\documentclass[11pt]{article}
\usepackage[margin = 1in]{geometry}

\usepackage{lipsum}
\usepackage{amsfonts}
\usepackage{graphicx}
\usepackage{epstopdf}

\usepackage{mathtools,mathrsfs,mathabx}
\usepackage{empheq}
\usepackage{amsfonts}
\usepackage{amsgen,amsthm,amsmath,amstext,amsbsy,amsopn,amssymb,subfigure,stmaryrd,bbm,dsfont}
\usepackage{comment}
\usepackage{enumitem}
\usepackage{booktabs}
\usepackage{blkarray}
\usepackage{array}
\usepackage{algorithm}
\usepackage{algpseudocode}
\usepackage{url}
\usepackage{longtable}
\usepackage{mathtools}
\usepackage{multirow}
\usepackage{color}

\ifpdf
  \DeclareGraphicsExtensions{.eps,.pdf,.png,.jpg}
\else
  \DeclareGraphicsExtensions{.eps}
\fi

\newtheorem{example}{Example}  

\newtheorem{theorem}{Theorem}

\newtheorem{lemma}{Lemma}
\newtheorem{remark}{Remark}
\newtheorem{corollary}{Corollary}

\newtheorem{proposition}{Proposition}

\DeclareMathAlphabet\mathbfcal{OMS}{cmsy}{b}{n}

\newcommand{\be}{\begin{equation}}
\newcommand{\ee}{\end{equation}}
\newcommand{\bea}{\begin{eqnarray}}
\newcommand{\eea}{\end{eqnarray}}
\newcommand{\beas}{\begin{eqnarray*}}
	\newcommand{\eeas}{\end{eqnarray*}}

\newcommand{\bbR}{\mathbb{R}}
\newcommand{\bbS}{\mathbb{S}}

\newcommand{\cU}{\mathcal{U}}
\newcommand{\cV}{\mathcal{V}}
\newcommand{\cH}{\mathcal{H}}
\newcommand{\cL}{\mathcal{L}}

\newcommand{\0}{{\mathbf{0}}}

\renewcommand{\u}{{\mathbf{u}}}

\newcommand{\X}{{\mathbf{X}}}
\newcommand{\B}{{\mathbf{B}}}

\newcommand{\D}{{\mathbf{D}}}

\newcommand{\I}{{\mathbf{I}}}
\newcommand{\M}{{\mathbf{M}}}

\renewcommand{\P}{{\mathbf{P}}}
\renewcommand{\O}{{\mathbf{O}}}

\newcommand{\R}{{\mathbf{R}}}
\renewcommand{\L}{{\mathbf{L}}}
\renewcommand{\S}{{\mathbf{S}}}

\newcommand{\U}{{\mathbf{U}}}
\newcommand{\V}{{\mathbf{V}}}
\newcommand{\W}{{\mathbf{W}}}

\newcommand{\A}{{\mathbf{A}}}
\newcommand{\Y}{{\mathbf{Y}}}
\newcommand{\Z}{{\mathbf{Z}}}

\newcommand{\cN}{{\cal N}}
\newcommand{\cM}{{\cal M}}

\newcommand{\scC}{{\mathscr{C}}}

\newcommand{\bOmega}{\boldsymbol{\Omega}}
\newcommand{\bSigma}{\boldsymbol{\Sigma}}
\newcommand{\bDelta}{\boldsymbol{\Delta}}

\newcommand{\rank}{{\rm rank}}

\newcommand{\tr}{{\rm tr}}

\newcommand{\Hess}{{\rm Hess}}
\newcommand{\grad}{{\rm grad}}

\newcommand{\F}{{\rm F}}
\newcommand{\rspan}{{\rm span}}
\newcommand{\st}{{\rm St}}
\newcommand{\sym}{{\rm Sym}}
\renewcommand{\skew}{{\rm Skew}}
\newcommand{\rmD}{{\rm D}}
\newcommand{\GL}{{\rm GL}}

\newcommand{\argmin}{\mathop{\rm arg\min}}

\newcommand{\bbO}{\mathbb{O}}

\makeatletter
\newcommand*{\rom}[1]{\expandafter\@slowromancap\romannumeral #1@}
\makeatother

\allowdisplaybreaks

\graphicspath{ {images/} }

\begin{document}
	\title{On Geometric Connections of Embedded and Quotient Geometries in Riemannian Fixed-rank Matrix Optimization}
	
	\author{Yuetian Luo$^1$, ~ Xudong Li$^2$, ~ and ~ Anru R. Zhang$^{3}$}
	
	\date{}
	\maketitle

	\footnotetext[1]{Data Science Institute, University of Chicago.  (\texttt{yuetian@uchicago.edu})}
	\footnotetext[2]{School of Data Science, Fudan University. (\texttt{lixudong@fudan.edu.cn})}
	\footnotetext[3]{Department of Biostatistics and Bioinformatics, Duke University. (\texttt{anru.zhang@duke.edu})}
	
	\bigskip

\begin{abstract}
In this paper, we propose a general procedure for establishing the geometric landscape connections of a Riemannian optimization problem under the embedded and quotient geometries. By applying the general procedure to the fixed-rank positive semidefinite (PSD) and general matrix optimization, we establish an exact Riemannian gradient connection under two geometries at every point on the manifold and sandwich inequalities between the spectra of Riemannian Hessians at Riemannian first-order stationary points (FOSPs). These results immediately imply an equivalence on the sets of Riemannian FOSPs, Riemannian second-order stationary points (SOSPs), and strict saddles of fixed-rank matrix optimization under the embedded and the quotient geometries. To the best of our knowledge, this is the first geometric landscape connection between the embedded and the quotient geometries for fixed-rank matrix optimization and it provides a concrete example of how these two geometries are connected in Riemannian optimization. In addition, the effects of the Riemannian metric and quotient structure on the landscape connection are discussed. We also observe an algorithmic connection between two geometries with some specific Riemannian metrics in fixed-rank matrix optimization: there is an equivalence between gradient flows under two geometries with shared spectra of Riemannian Hessians. A number of novel ideas and technical ingredients including a unified treatment for different Riemannian metrics, novel metrics for the Stiefel manifold, and new horizontal space representations under quotient geometries are developed to obtain our results. The results in this paper deepen our understanding of geometric and algorithmic connections of Riemannian optimization under different Riemannian geometries and provide a few new theoretical insights to unanswered questions in the literature.
\end{abstract}

\section{Introduction}\label{sec:intro}

Riemannian optimization is a powerful method to tackle a general class of optimization problems with geometric constraints so that the solution is constrained to a Riemannian manifold. One key component in Riemannian optimization is the Riemannian geometry. Over the past decades, numerous Riemannian geometries, including different manifold classes, quotient manifold structures, and Riemannian metrics, have been proposed in various problems for having either better geometric properties or faster algorithmic convergences \cite{absil2009geometric,bonnabel2010riemannian,edelman1998geometry,massart2020quotient,mishra2012riemannian,mishra2014fixed,mishra2016riemannian,vandereycken2013riemannian}. In the literature, two popular choices of manifold classes in Riemannian optimization are {\it embedded submanifold} and {\it quotient manifold}. The embedded geometry often allows computing and interpreting the geometric notions more straightforwardly; while the optimization methods via quotient geometry can be more versatile as quotient geometry provides more choices of quotient structures and Riemannian metrics. The readers are referred to  \cite{absil2009optimization,boumal2020introduction,meyer2011geometric, mishra2014fixed} for surveys on these topics.

The Riemannian optimization under embedded and quotient geometries are not obviously related (even commented to be fundamentally different in \cite{journee2010low}) and are often studied separately in the literature. It is also unclear how to choose between these two geometries in Riemannian optimization. On the other hand, a few empirical studies on the algorithmic comparisons between the embedded and the quotient geometries in matrix completion and graph-based clustering problems showed that for either the gradient or trust-region-based algorithms, these two geometries perform more or less the same in terms of the total computational time \cite{douik2021low,mishra2012riemannian,mishra2014fixed}. It has been asked by \cite{vandereycken2013low} on the reason behind. \cite{vandereycken2010riemannian} hinted that embedded and quotient approaches are probably related as the manifolds under these two geometries are diffeomorphic to each other. However, it remains elusive in the literature how they are exactly connected in specific Riemannian optimization problems from either an algorithmic or a geometric point of view.

In this work, we make the first attempt to answer these questions by proposing a general framework to investigate the first-order and second-order geometric landscape connections of an optimization problem under the embedded and quotient geometries. Throughout this paper, the concept of landscape refers to the objective surfaces or the global geometry of the objective function, i.e., a characterization of all stationary points and an explicit description of geometry on the entire parameter domain \cite{ge2017optimization,li2016symmetry,zhang2020symmetry}. The first-order geometric connection can often be easily established and the general procedure for connecting second-order geometries includes three steps: 1. compute the quadratic form of Riemannian Hessians under two geometries; 2. construct a carefully-designed mapping $\cL$ between the horizontal space under the quotient geometry to the tangent space under the embedded geometry at proper reference points to connect Riemannian Hessians; 3. establish the spectra connection between Riemannian Hessians via bounding the spectrum of $\cL$.

We then specifically consider the following fixed-rank matrix optimization problems:
\begin{equation} \label{eq: PSD-manifold-formulation}
	\text{PSD case}: \quad \quad  \min_{\X \in \bbS^{p \times p} \succcurlyeq 0, \rank(\X) = r} f(\X), \quad 0 < r \leq p,
	\end{equation}
\begin{equation} \label{eq: general prob}
	\text{general case}: \quad \min_{\X \in \bbR^{p_1 \times p_2}, \rank(\X) = r} f(\X), \quad 0 < r \leq \min\{p_1,p_2\}.
\end{equation} 
In the positive semidefinite (PSD) case, without loss of generality, we assume $f$ is symmetric in $\X$, i.e., $f(\X) = f(\X^\top)$; otherwise, we can set $\tilde{f}(\X) = \frac{1}{2}(f(\X) + f(\X^\top))$ and have $\tilde{f}(\X) = f(\X)$ for all $\X \succcurlyeq 0$ without changing the problem \cite{bhojanapalli2016dropping}. In both cases, we assume $f$ is twice continuously differentiable with respect to $\X$ and the Euclidean metric. Both embedded and quotient geometries have been studied for the sets of fixed-rank matrices and many algorithms have been proposed for \eqref{eq: PSD-manifold-formulation} and \eqref{eq: general prob} on each individual Riemannian geometry. See Section \ref{sec: related-literature} for a review of existing results. 

By applying the general procedure, we establish the geometric connections of \eqref{eq: PSD-manifold-formulation} and \eqref{eq: general prob} under the embedded and a variety of quotient geometries (Theorems \ref{th: embedded-quotient-connection-PSD1}--\ref{th: embedded-quotient-connection-general3}, Corollaries \ref{coro: landscape connection PSD}, \ref{coro: landscape connection general case}) informally summarized as follows.
\begin{theorem}[Informal results] Consider optimization problems \eqref{eq: PSD-manifold-formulation} and \eqref{eq: general prob} on fixed-rank PSD and general matrix manifolds. 
\begin{itemize}[leftmargin=*]
	\item There exists an equivalence relation on the sets of Riemannian first-order stationary points (FOSPs), Riemannian second-order stationary points (SOSPs), and strict saddles of \eqref{eq: PSD-manifold-formulation} or \eqref{eq: general prob} under embedded and quotient geometries.
	\item There exists a one-to-one relation between Riemannian gradients under embedded and quotient geometries at every point on the manifold; the spectra of Riemannian Hessians of \eqref{eq: PSD-manifold-formulation} or \eqref{eq: general prob} under two geometries are sandwiched by each other at Riemannian FOSPs.
\end{itemize}
\end{theorem} 
To the best of our knowledge, this is the first geometric landscape connection between the embedded and the quotient geometries for fixed-rank matrix optimization.

 In addition, the effects of Riemannian metric and quotient structure on the landscape connection, i.e., the relationship between the stationary points and surface geometries under the embedded and quotient geometries, are discussed. We also discover an intriguing algorithmic connection of \eqref{eq: PSD-manifold-formulation} and \eqref{eq: general prob} under the embedded and the quotient geometries with some specific Riemannian metrics. Specifically, we find the gradient flows under two geometries are closely related when the spectra of Riemannian Hessians under two geometries are close at Riemannian FOSPs, and they are {\it identical} when the spectra of two Riemannian Hessians coincide. Such results provide theoretical insights into empirical observations in \cite{mishra2012riemannian}. 

In a broad sense, embedded and quotient geometries are the most common two choices in Riemannian optimization. It is known that the manifolds under two geometries are diffeomorphic under proper assumptions \cite[Proposition 3.5.23]{abraham2012manifolds}. See details in Section \ref{sec: general-strategy-for-connection}.
However, it is unclear how the geometry-dependent key concepts of optimization problems are related. This paper bridges them from geometric and algorithmic points of view and illustrates explicitly how they are connected in solving fixed-rank matrix optimization problems.

\subsection{Related Literature} \label{sec: related-literature}
This work is related to a range of literature on low-rank matrix optimization, Riemannian/nonconvex optimization, and geometric landscape analysis of an optimization problem. 

First, choosing a proper Riemannian geometry is undoubtedly a central topic in Riemannian optimization and numerous geometries have been proposed under different considerations. For example, \cite{vandereycken2013riemannian} proposed a homogeneous space geometry on the set of fixed-rank PSD matrices such that complete geodesics can be obtained. \cite{mishra2012riemannian,mishra2016riemannian} proposed new Riemannian metrics under fixed-rank quotient geometries tailored to objective functions. Different Riemannian manifold structures have also been considered to have better geometric properties \cite{absil2014two,bonnabel2010riemannian}. In this work, we focus on the choice between embedded and quotient geometries and study its effect on the corresponding Riemannian optimization problem.

Second, for fixed-rank matrix optimization \eqref{eq: PSD-manifold-formulation} and \eqref{eq: general prob}, a number of Riemannian optimization methods, including (conjugate) gradient descent, (Gauss-)Newton, trust-region have been developed under either embedded geometry \cite{hou2020fast,luo2020recursive,shalit2012online,vandereycken2013low,vandereycken2010riemannian,wei2016guarantees} or quotient geometry \cite{absil2014two,absil2009geometric,boumal2011rtrmc,edelman1998geometry,huang2018blind,meyer2011linear,meyer2011regression,mishra2014fixed,ngo2012scaled}. We refer readers to \cite{absil2009optimization,boumal2020introduction,cai2018exploiting} for the recent algorithmic development in Riemannian matrix optimization. In addition to Riemannian optimization, a number of other methods including convex relaxation \cite{cai2013compressed,recht2010guaranteed}, non-convex factorization \cite{jain2013low,ma2019implicit,sun2015guaranteed,tu2016low}, projected gradient descent \cite{jain2010guaranteed}, and penalty method \cite{gao2010majorized} have been proposed to solve \eqref{eq: PSD-manifold-formulation} and \eqref{eq: general prob} as well. A comparison of these different approaches can be found in \cite{chen2018harnessing,chi2019nonconvex}. 

Third, a few attempts have been made to analyze the landscape of a Riemannian matrix optimization problem. For example, \cite{ahn2021riemannian,maunu2019well} provided landscape analyses for robust subspace recovery and matrix factorization, respectively, over the Grassmannian manifold. Under the embedded geometry, \cite{uschmajew2018critical} showed the landscape of \eqref{eq: general prob} is benign when $f$ is quadratic and satisfies certain restricted spectral bounds properties. Different from this line of work which focuses on the landscape of the problem under a single Riemannian geometry when $f$ is well-conditioned, here we study the geometric landscape connections of \eqref{eq: PSD-manifold-formulation} and \eqref{eq: general prob} under the embedded and the quotient geometries for a general $f$. 

Finally, there are a few recent studies on the geometric connections of different approaches for rank-constrained optimization. For example, \cite{ha2020equivalence} studied the relationship between Euclidean FOSPs/SOSPs under the factorization formulation and fixed points of the projected gradient descent (PGD) in the general low-rank matrix optimization with the bounded rank constraint. They showed while the sets of Euclidean FOSPs under the factorization formulation can be larger, the sets of Euclidean SOSPs are contained in the set of fixed points of the PGD with a small step size. Another related work is \cite{luo2021nonconvex}, where they studied the landscape connections of the Euclidean factorization formulation and the Riemannian formulation with embedded geometry for both \eqref{eq: PSD-manifold-formulation} and \eqref{eq: general prob} and identified the close connections between the landscapes under these two formulations, e.g., the set of FOSPs and SOSPs under these two formulations are exactly the same when constraining to rank-$r$ matrices. Different from the Euclidean factorization formulations considered in \cite{ha2020equivalence,luo2021nonconvex}, we consider the formulation under quotient geometry, where the invariance in matrix factorization is canceled out by quotient, so the regularization term in the factorization formulation is not needed here. Quotient geometry and Euclidean geometry for matrix factorization are different in general (see the forthcoming Remark \ref{rem: choise of W-L-R}). Moreover, we consider more quotient geometries induced from different matrix factorizations than \cite{ha2020equivalence,luo2021nonconvex} and develop a unified procedure for connecting landscape geometries under the embedded and the quotient geometries. Finally, we also discover a new algorithmic connection between the gradient flows under the corresponding embedded and quotient geometries (see Remark \ref{rem: algorithmic-connection}) and this is not covered in \cite{ha2020equivalence,luo2021nonconvex}. In addition, under the quotient geometry, a general equivalence relation on the sets of Riemannian FOSPs and SOSPs of objectives on the total space and the quotient space is given in \cite[Section 9.11]{boumal2020introduction}. Complementary to these results, here we consider the geometric landscape connections of \eqref{eq: PSD-manifold-formulation} and \eqref{eq: general prob} under two Riemannian formulations with different geometries, i.e., embedded and quotient geometries.  

\subsection{Organization of the Paper}\label{sec: organization}

The rest of this article is organized as follows. After a brief introduction of notation in Section \ref{sec: notation} and preliminaries for Riemannian optimization on embedded and quotient manifolds in Section \ref{sec: Riemannian-opt-background}, we introduce a general procedure for establishing landscape connections of an optimization problem under embedded and quotient geometries in Section \ref{sec: general-strategy-for-connection}. Then we discuss the embedded and quotient geometries for fixed-rank matrix manifolds in Section \ref{sec: embedded-quotient-fixed-rank-matrix}. By applying the general procedure, our main results on the geometric landscape connections of \eqref{eq: PSD-manifold-formulation} and \eqref{eq: general prob} under the embedded and the quotient geometries are presented in Sections \ref{sec: connection-PSD} and \ref{sec: connection-general}, respectively. Conclusion and future work are given in Section \ref{sec: conclusion}. The proof for the main results is provided in the main body of the paper. Additional preliminaries, proofs, and lemmas are collected in Appendices \ref{sec: additional-preliminaries}-\ref{sec: additional-lemmas}, respectively.

\subsection{Notation and Preliminaries} \label{sec: notation}

The following notation will be used throughout this article. For any positive integer $p$, denote $[p]=\{1,\ldots, p\}$. We use $\bbR^{p_1 \times p_2}$, $\bbS^{p \times p}$, $\bbS_+(r)$, $\bbR^{p \times r}_*$ and $\st(r,p)$ to denote the spaces of $p_1$-by-$p_2$ real matrices, $p$-by-$p$ real symmetric matrices, $r$-by-$r$ real symmetric positive definite matrices, $p$-by-$r$ real full column rank matrices and $p$-by-$r$ real matrices with orthonormal columns, respectively. We also let $\bbO_r$ be the set of $r$-by-$r$ orthogonal matrices, i.e., $\bbO_r = \st(r,r)$. Uppercase and lowercase letters (e.g., $A, B, a, b$), lowercase boldface letters (e.g., $\mathbf{u}, \mathbf{v}$), uppercase boldface letters (e.g., $\U, \V$) are often used to denote scalars, vectors, matrices, respectively. For any $a, b \in \bbR$, let $a \wedge b := \min\{a,b\}, a \vee b := \max\{a,b\}$. For any matrix $\X \in \mathbb{R}^{p_1\times p_2}$ with singular value decomposition (SVD) $\sum_{i=1}^{p_1 \land p_2} \sigma_i(\X) \mathbf{u}_i \mathbf{v}_i^\top$, where $\sigma_1(\X) \geq \sigma_2(\X) \geq \cdots \geq \sigma_{p_1 \wedge p_2} (\X)$, denote its Frobenius norm and spectral norm as $\|\X\|_\F = \sqrt{\sum_{i} \sigma^2_i(\X)}$ and $\|\X\| = \sigma_1(\X)$, respectively. Also, denote $\X^{-1}$, $\X^{-\top}$ and $\X^\dagger$ as the inverse, transpose inverse, and Moore-Penrose inverse of $\X$, respectively. For any $\X \in \bbR^{p \times p}$, let $\sym(\X) = (\X + \X^\top)/2$,  $\skew(\X) =(\X - \X^\top)/2$, and $\tr(\X)$ be the symmetric part, skew-symmetric part, and the trace of $\X$, respectively. For any $\X \in \bbS^{p \times p}$ having eigendecomposition $\U \bSigma \U^\top$ with non-increasing eigenvalues on the diagonal of $\bSigma$, let $\lambda_i(\X)$ be the $i$-th largest eigenvalue of $\X$, $\lambda_{\min}(\X)$ be the least eigenvalue of $\X$, and $\X^{1/2} = \U \bSigma^{1/2} \U^\top$. We note $\X \succcurlyeq 0$ if $\X$ is a symmetric positive semidefinite (PSD) matrix. Throughout the paper, the SVD (or eigendecomposition) of a rank $r$ matrix $\X$ (or symmetric matrix $\X$) refers to its economic version and we say a column orthonormal matrix $\U'$ spans the top $r$ left singular space or eigenspace of $\X$ if $\U' = \U \O$ for some $\O \in \bbO_r$, where $\U$ is formed by the top $r$ left singular vectors or eigenvectors of $\X$. For any $\U\in \st(r,p)$, $P_{\U} = \U\U^\top$ represents the orthogonal projector onto the column space of $\U$; we also note $\U_\perp\in \st(p-r, p)$ as an orthonormal complement of $\U$. We use bracket subscripts to denote sub-matrices. For example, $\X_{[i_1,i_2]}$ is the entry of $\X$ on the $i_1$-th row and $i_2$-th column. In addition, $\I_r$ is the $r$-by-$r$ identity matrix.  Finally, the dimension of a linear space $\cV$ is denoted as $\dim(\cV)$. For any two linear spaces $\cV_1, \cV_2$, the sum of $\cV_1$ and $\cV_2$ is denoted by $\cV_1+\cV_2 := \{\mathbf{v}_1 + \mathbf{v}_2| \mathbf{v}_1 \in \cV_1,\mathbf{v}_2 \in \cV_2\}$. If every vector in $\cV_1 + \cV_2$ can be uniquely decomposed into $\mathbf{v}_1 + \mathbf{v}_2$, where $\mathbf{v}_1 \in \cV_1,\mathbf{v}_2 \in \cV_2$, then we call the sum of $\cV_1$ and $\cV_2$ the direct sum, denoted by $\cV_1 \oplus \cV_2$. The direct sum satisfies a key property: $\dim(\cV_1 \oplus \cV_2) = \dim(\cV_1) + \dim(\cV_2)$. For any two Euclidean spaces $\cV_1$ and $\cV_2$ endowed with inner product $g(\cdot, \cdot)$, we say $\cV_1$ is orthogonal to $\cV_2$ with respect to $g$ and note $\cV_1 \perp \cV_2$, if and only if $g(\mathbf{v}_1, \mathbf{v}_2) =0$ for any $\mathbf{v}_1 \in \cV_1,\mathbf{v}_2 \in \cV_2$.

Suppose $f: \bbR^{p_1 \times p_2} \to \bbR$ is a differentiable scalar function and $\phi: \bbR^{p_1 \times p_2} \to \bbR^{q_1 \times q_2}$ is a differentiable matrix-valued function. Let the Euclidean gradient of $f$ at $\X$ be $\nabla f(\X)$, i.e., $(\nabla f(\X))_{[i,j]} = \frac{\partial f(\X)}{\partial \X_{[i,j]}}$ for $i\in [p_1], j\in [p_2]$. The Euclidean gradient of $\phi$ is a linear operator from $\bbR^{p_1 \times p_2}$ to $\bbR^{q_1 \times q_2}$ such that $(\nabla \phi (\X) [\Z] )_{[i,j]} = \sum_{k \in [p_1],l \in [p_2]} \frac{\partial ( \phi (\X) )_{[i,j]} }{\partial \X_{[k,l]}} \Z_{[k,l]}$ for any $\Z \in \bbR^{p_1 \times p_2}, i \in [q_1], j \in [q_2]$. For a twice continuously differentiable function $f$, let $\nabla^2 f(\X)[\cdot]$ be its Euclidean Hessian, which is the gradient of $\nabla f(\X)$ and can be viewed as a linear operator from $\bbR^{p_1 \times p_2}$ to $\bbR^{p_1 \times p_2}$ satisfying
\begin{equation*}
	( \nabla^2 f(\X)[\Z] )_{[i,j]} =  \sum_{k \in [p_1],l \in [p_2]} \frac{\partial (\nabla f (\X) )_{[i,j]} }{\partial \X_{[k,l]}} \Z_{[k,l]} = \sum_{k \in [p_1],l \in [p_2]} \frac{\partial^2 f (\X) }{\partial \X_{[k,l]} \partial \X_{[i,j]} } \Z_{[k,l]}.
\end{equation*} 
Define the bilinear form of the Hessian of $f$ as $\nabla^2 f(\X)[\Z_1, \Z_2]:= \langle \nabla^2 f(\X)[\Z_1], \Z_2 \rangle $ for any $\Z_1, \Z_2 \in \bbR^{p_1 \times p_2}$, where $\langle \cdot, \cdot \rangle$ is the standard Euclidean inner product.

\section{Riemannian Optimization Under Embedded and Quotient Geometries} \label{sec: Riemannian-opt-background}
In this section, we first give a brief introduction to Riemannian optimization and then discuss how to perform Riemannian optimization under embedded and quotient geometries. 

Riemannian optimization concerns optimizing a real-valued function $f$ defined on a Riemannian manifold $\cM$. The readers are referred to \cite{absil2009optimization,boumal2020introduction,hu2020brief} for more details. The calculations of Riemannian gradients and Riemannian Hessians are key ingredients to perform continuous optimization over the Riemannian manifold. Suppose $\X \in \cM$, $g_{\X}( \cdot, \cdot )$ is the Riemannian metric, and $T_\X \cM$ is the tangent space of $\cM$ at $\X$. Then a {\it vector field} $\xi$ on a manifold $\cM$ is a mapping that assigns each point $\X \in \cM$ a tangent vector $\xi_\X \in T_\X \cM$. Suppose $\cM'$ is another smooth manifold equipped with the metric $g'$. For a smooth map between two manifolds $F: \cM \to \cM'$, the {\it differential} of $F$ at $\X$ is a linear operator $\rmD F(\X): T_\X \cM \to T_{F(\X)} \cM'$ defined by $\rmD F(\X)[\mathbf{v}]:= \frac{d F( c(t) )}{dt}|_{t=0} $, where $c$ is a smooth curve on $\cM$ passing through $\X$ at $t = 0$ with velocity $\mathbf{v}$. 

 The {\it Riemannian gradient} of a smooth function $f:\cM \to \bbR$ at $\X$ is defined as the unique tangent vector, ${\rm grad}\, f(\X) \in T_\X \cM$, such that $g_\X( \grad \, f(\X),\xi_\X ) = {\rm D} \, f(\X)[\xi_\X], \forall\, \xi_\X \in T_\X \cM$, where ${\rm D} f(\X)[\xi_\X]$ is the differential of $f$ at point $\X$ along the direction $\xi_\X$. The {\it Riemannian Hessian} of $f$ at $\X\in\cM$ is a linear mapping ${\rm Hess} \,f(\X): T_\X\cM \to T_\X\cM$ defined as
\begin{equation} \label{def: Riemannain-Hessian}
	{\rm Hess}\, f(\X)[\xi_\X] = \widebar{\nabla}_{\xi_\X} {\rm grad}\, f \in T_\X \cM, \,\quad \forall \xi_\X \in T_\X \cM,
\end{equation}
where $\widebar{\nabla}$ is the {\it Riemannian connection} on $\cM$, which is a generalization of the directional derivative on a vector field to Riemannian manifolds \cite[Section 5.3]{absil2009optimization} (See Appendix \ref{sec: additional-preliminaries} for more details). The bilinear form of Riemannian Hessian is defined as $\Hess f(\X)[\xi_\X, \theta_\X]:= g_\X ( \Hess f(\X)[\xi_\X], \theta_\X )$ for any $\xi_\X, \theta_\X \in T_{\X}\cM$. We say $\X \in \cM$ is a {\it Riemannian FOSP} of $f$ if $\grad f(\X) = \0$ and a {\it Riemannian SOSP} of $f$ if $\grad f(\X) = \0$ and $\Hess f(\X) \succcurlyeq 0$. Moreover, $\X \in \cM$ is a {\it local minimizer} of $f$ if there exists a neighborhood $\cN$ of $\X$ in $\cM$ such that $f(\X) \leq f(\X')$ for all $\X' \in \cN$. Finally, we call a Riemannian FOSP a {\it strict saddle} if the Riemannian Hessian evaluated at this point has a strict negative eigenvalue. 

In this work, we mainly focus on two classes of manifolds: embedded submanifold and quotient manifold. The embedded submanifold can be viewed as a generalization of the notion of surface in $\bbR^{d}$ and the Riemannian gradients and Hessians under the embedded geometry can often be concretely written out as every geometric object lies in the embedding space. For example, suppose $\cM$ is a Riemannian embedded submanifold of the Riemannian manifold $\widetilde{\cM}$ and the objective function $f: \cM \to \bbR$ is the restriction of $\tilde{f}: \widetilde{\cM} \to \bbR$ to the embedded submanifold $\cM$. Then we have the following simple expressions for the Riemannian gradient of $f$ and the Riemannian connection \cite[Eq. (3.37) and Proposition 5.3.2]{absil2009optimization}: $\grad f(\X) = P_{T_\X \cM}( \grad \tilde{f} (\X) ) $ and $\widebar{\nabla}_{\xi_\X} \eta = P_{T_\X \cM}( \widebar{\nabla}'_{\xi_\X} \eta  ),$
where $P_{T_\X \cM}(\cdot)$ is the projection operator onto the tangent space $T_\X \cM$, $\xi, \eta$ are two vector fields on $\cM$, $\grad \tilde{f} (\X)$ and $\widebar{\nabla}'$ are the Riemannian gradient of $\tilde{f}$ and the Riemannian connection on $\widetilde{\cM}$, respectively. On the other hand, the geometric objects under quotient manifolds are more abstract. The following Section \ref{sec: Riemannian-opt-quotient} aims to provide more details on how to perform Riemannian optimization on quotient manifolds. 

\subsection{Riemannian Optimization on Quotient Manifolds} \label{sec: Riemannian-opt-quotient}

Quotient manifolds are often defined via an equivalence relation ``$\sim$'' that satisfies symmetric, reflexive, and transitive properties \cite[Section 3.4.1]{absil2009optimization}. The equivalence classes are often abstract objects and cannot be directly applied in numerical computations. Riemannian optimization on quotient manifolds works on representatives of these equivalence classes instead. To be specific, suppose $\widebar{\cM}$ is an embedded submanifold equipped with an equivalence relation $\sim$. The {\it equivalence class} (or {\it fiber}) of $\widebar{\cM}$ at a given point $\X$ is defined by the set $[\X] = \{\X_1 \in \widebar{\cM}: \X_1 \sim \X \}$. The set $ \cM := \widebar{\cM}/\sim = \{[\X]: \X \in \widebar{\cM} \}$ is called a {\it quotient} of $\widebar{\cM}$ by $\sim$. The mapping $\pi: \widebar{\cM} \to \widebar{\cM}/\sim$, $ \X \mapsto [\X]$ is called the {\it quotient map} or {\it canonical projection} and the set $\widebar{\cM}$ is called the {\it total space} of the quotient $\widebar{\cM}/\sim$. If $\cM$ further admits a smooth manifold structure and $\pi$ is a smooth submersion, then we call $\cM$ a {\it quotient manifold} of $\widebar{\cM}$. Here, a smooth map $F: \cM \to \cM'$ is called a {\it smooth submersion} if its differential $\rmD F$ is surjective at each point of $\cM$ \cite[Chapter 4]{lee2013smooth}.

Due to the abstractness, the tangent space $T_{[\X]} \cM$ of $\cM$ at $[\X]$ calls for a representation in the tangent space $T_{\X}\widebar{\cM}$ of the total space $\widebar{\cM}$. By the equivalence relation $\sim$, the representation of elements in $T_{[\X]} \cM$ should be restricted to the directions in $T_\X \widebar{\cM}$ without inducing displacement along the equivalence class $[\X]$. This can be achieved by decomposing $T_{\X}\widebar{\cM}$ into complementary spaces $T_{\X}\widebar{\cM} = \cV_{\X} \widebar{\cM} \oplus \cH_{\X} \widebar{\cM}$. Here, $\cV_{\X} \widebar{\cM}$ is called the {\it vertical space} that contains tangent vectors of the equivalence class $[\X]$. $\cH_{\X} \widebar{\cM}$ is called the {\it horizontal space} of $T_{\X}\widebar{\cM}$, which is complementary to $\cV_\X \widebar{\cM}$ and provides a proper representation of the abstract tangent space $T_{[\X]} \cM$ \cite[Section 3.5.8]{absil2009optimization}. Once $\widebar{\cM}$ is endowed with $\cH_{\X}\widebar{\cM}$, a given tangent vector $\eta_{[\X]} \in T_{[\X]} \cM$ at $[\X]$ is uniquely represented by a horizontal tangent vector $\eta_\X \in \cH_\X \widebar{\cM}$ that satisfies ${\rm D} \pi (\X)[\eta_\X] = \eta_{[\X]}$ \cite[Section 3.5.8]{absil2009optimization}. The tangent vector $\eta_\X \in \cH_\X \widebar{\cM}$ is also called the {\it horizontal lift} of $\eta_{[\X]}$ at $\X$.

Next, we introduce the notion of {\it Riemannian quotient manifolds}. Suppose the total space $\widebar{\cM}$ is endowed with a Riemannian metric $\bar{g}_\X$, and for every $[\X] \in \cM$ and every $\eta_{[\X]}, \theta_{[\X]} \in T_{[\X]} \cM$, the expression $\bar{g}_\X(\eta_\X, \theta_\X )$, i.e., the inner product of the horizontal lifts of $\eta_{[\X]}, \theta_{[\X]}$ at $\X$, does not depend on the choice of the representative $\X$. Then the metric $\bar{g}_\X$ in the total space induces a metric $g_{[\X]}$ on the quotient space, i.e., $g_{[\X]}(\eta_{[\X]}, \theta_{[\X]}):= \bar{g}_\X(\eta_\X, \theta_\X)$. The quotient manifold $\cM$ endowed with $g_{[\X]}$ is called a {\it Riemannian quotient manifold} of $\widebar{\cM}$ and the quotient mapping $\pi: \widebar{\cM} \to \cM$ is called a {\it Riemannian submersion} \cite[Section 3.6.2]{absil2009optimization}. 
Optimization on Riemannian quotient manifolds is particularly convenient because the computation of representatives of Riemannian gradients and Hessians in the abstract quotient space can be directly performed by means of their analogous in the total space. To be specific, suppose $\bar{f}: \widebar{\cM} \to \bbR$ is an objective function in the total space and is invariant along the fiber, i.e., $\bar{f}(\X_1) = \bar{f}(\X_2)$ whenever $\X_1 \sim \X_2$. Then $\bar{f}$ induces a function $f: \cM \to \bbR$ on the quotient space and the horizontal lift of the Riemannian gradient of $f$ can be obtained as follows \cite[Section 3.4.2]{meyer2011geometric}:
\begin{equation} \label{eq: quotient-gradient-general}
	\overline{\grad f([\X])} = P_\X^\cH (\grad \bar{f} (\X) ).
\end{equation} 
Here $P_\X^\cH(\cdot)$ denotes the projection operator onto the horizontal space at $\X$ ($P_\X^\cH(\cdot)$ depends on the metric $\bar{g}_\X$) and $\grad \bar{f}(\X)$ denotes the Riemannian gradient of $\bar{f}$ at $\X$ in the total space. The following Lemma \ref{lm: quotient-ortho-hori-gradient-express} shows that if the horizontal space is canonically chosen, i.e., $\cH_\X \widebar{\cM}$ is the orthogonal complement of $\cV_\X \widebar{\cM}$ in $T_\X \widebar{\cM}$ with respect to $\bar{g}_\X$, then $\grad \bar{f} (\X)$ automatically lies in the horizontal space at $\X$. 
\begin{lemma}{\rm (\cite[Section 3.6.2]{absil2009optimization})} \label{lm: quotient-ortho-hori-gradient-express}
	Suppose $\cH_\X \widebar{\cM}:= \{ \eta_\X \in T_\X \widebar{\cM}: \bar{g}_\X(\eta_\X, \theta_\X) = 0 \text{ for all } \theta_\X \in \cV_\X \widebar{\cM} \}$. Then $\overline{\grad f([\X])} = \grad \bar{f} (\X)$.
\end{lemma}

Finally, the Riemannian connection on the quotient manifold $\cM$ can also be uniquely represented by the Riemannian connection in the total space $\widebar{\cM}$. Suppose $\eta, \theta$ are two vector fields on $\cM$ and $\eta_\X$ and $\theta_\X$ are the horizontal lifts of $\eta_{[\X]}$ and $\theta_{[\X]}$ in $\cH_\X\widebar{\cM}$. Then the horizontal lift of $\widebar{\nabla}_{\theta_{[\X]}} \eta$ on the quotient manifold is given by $\overline{ \widebar{\nabla}_{\theta_{[\X]}} \eta } = P_\X^{\cH} (\widebar{\nabla}_{\theta_\X} \bar{\eta})$, where $\bar{\eta}$ denotes the horizontal lift of the vector field $\eta$ and $\widebar{\nabla}_{\theta_\X} \bar{\eta}$ is the Riemannian connection in the total space. Combining \eqref{def: Riemannain-Hessian}, we have the horizontal lift of the Riemannian Hessian of $f$ on $\cM$ satisfies
\begin{equation} \label{eq: quotient-hessian-linear-from}
	\overline{\Hess f([\X])[\theta_{[\X]}]} = P_\X^{\cH}\left( \widebar{\nabla}_{\theta_\X} \overline{\grad f} \right)
\end{equation} for any $\theta_{[\X]} \in T_{[\X]} \cM$ and its horizontal lift $\theta_\X$. We also define the bilinear form of the horizontal lift of the Riemannian Hessian as $\overline{\Hess f([\X])} [\theta_{\X}, \eta_\X ] := \bar{g}_\X\left( \overline{\Hess f([\X])[\theta_{[\X]}]}, \eta_\X \right) $ for any $\theta_\X, \eta_\X \in \cH_\X \widebar{\cM}$. Then, by recalling the definition of the Riemannian metric $g_{[\X]}$ in the quotient space, we have
\begin{equation*}
	\begin{split}
		\overline{\Hess f([\X])} [\theta_{\X}, \eta_\X ] &= \bar{g}_\X\left(\overline{\Hess f([\X])[\theta_{[\X]}]}, \eta_\X \right)\\
		& = g_{[\X]} \left( \Hess f([\X])[\theta_{[\X]}], \eta_{[\X]}  \right) = \Hess f([\X])[\theta_{[\X]}, \eta_{[\X]}].
	\end{split}
\end{equation*} 
So $\Hess f([\X])$ is completely characterized by $\overline{\Hess f([\X])}$ in the lifted horizontal space.

\section{A General Procedure for Establishing Geometric Connections in Riemannian Optimization} \label{sec: general-strategy-for-connection}

In this section, we present a general procedure for connecting landscape properties of an optimization problem (not necessarily restricted to fixed-rank matrix optimization) under different Riemannian geometries. For the convenience of presentation, we focus on the connection between embedded and quotient geometries, while this procedure can be applied in broader settings. 

Suppose $\cM^e$, endowed with the Riemannian metric $g_\X$, is a $d$-dimensional Riemannian embedded submanifold of the Riemannian manifold $\widebar{\cM}^e$. Here the superscript $e$ in $\cM^e$ stands for the ``embedded''. Consider the following optimization problem under $\cM^e$
\begin{equation}\label{eq: general-strategy-example-embedded}
    \text{(Opt. under Embedded Geometry)} \quad \min_{\X \in \cM^e} f(\X),
\end{equation}
where $f: \cM^e \to \bbR$ is twice continuously differentiable and is the restriction of $\bar{f}: \widebar{\cM}^e \to \bbR$ to the embedded submanifold $\cM^e$. Suppose $\widebar{\cM}^q$ is another $(d+m)$-dimensional smooth manifold and there exists a surjective submersion $\bar{\ell}: \Z \in \widebar{\cM}^q \to \bar{\ell}(\Z) \in \cM^e$. Then \eqref{eq: general-strategy-example-embedded} can be equivalently reformulated on $\widebar{\cM}^q$ as
\begin{equation} \label{eq: general-strategy-example}
	\min_{\Z \in \widebar{\cM}^q } \bar{h}(\Z) := f(\bar{\ell}(\Z)).
\end{equation} 
Here, the equivalence means (i) if $\Z^* \in \argmin_{\Z \in \widebar{\cM}^q } \bar{h}(\Z)$, then $\bar{\ell}(\Z^*) \in  \argmin_{\X \in \cM^e} f(\X)$; (ii) if $\X^*  \in  \argmin_{\X \in \cM^e} f(\X)$, then there exists an $\Z^* \in \widebar{\cM}^q$ such that $\Z^* \in \argmin_{\Z \in \widebar{\cM}^q } \bar{h}(\Z)$ and $\X^* = \bar{\ell}(\Z^*)$.
We will see concrete examples of $\cM^e$, $\widebar{\cM}^q$ and the mapping $\bar{\ell}$ later in the context of fixed-rank matrix optimization. The transformation between \eqref{eq: general-strategy-example-embedded} and \eqref{eq: general-strategy-example} can be regarded as a generalization of the classic technique of changing variables in the Riemannian optimization setting. Since $\bar{\ell}$ is a submersion, we can deduce from the rank theorem \cite[Theorem 4.12]{lee2013smooth} that by picking certain special coordinate charts $(\cU^q, \varphi^q)$ and $(\cU^e, \varphi^e)$ on $\widebar{\cM}^q$ and $\cM^e$, respectively, $\bar{h}$ has the coordinate representation, denoted by $\hat{\bar{h}}$, $ \hat{\bar{h}}(z_1,\ldots,z_{d+m}) = \bar{h} \circ (\varphi^q)^{-1} (z_1,\ldots,z_{d+m}) \overset{ \eqref{eq: general-strategy-example} } = f \circ (\varphi^e)^{-1} \varphi^e ( \bar{\ell}((\varphi^q)^{-1} (z_1,\ldots,z_{d+m})) )   = f \circ (\varphi^e)^{-1}(z_1,\ldots,z_d) = \hat{f}(z_1,\ldots,z_d)$ where $\hat{f}$ is the coordinate representation of $f$. This implies the Riemannian FOSPs and SOSPs of \eqref{eq: general-strategy-example-embedded} and \eqref{eq: general-strategy-example} 
are the same, i.e., $(\varphi^{q})^{-1}(z_1,\ldots,z_{d+m})$ is a FOSP or SOSP on $\widebar{\cM}^q$ if and only if $(\varphi^e)^{-1}(z_1,\ldots,z_d)$ is a FOSP or SOSP on $\cM^e$. However, the relationship between the spectra of Riemannian Hessians under two formulations is less clear as $\varphi^e$ and $\varphi^q$ are difficult to characterize in general. More importantly, since $\bar{\ell}$ is not necessarily a bijection, there is no one-to-one correspondence between points on $\widebar{\cM}^q$ and $\cM^e$. This motivates us to apply one additional quotient step. Suppose ``$\sim$'' is an equivalence relation on $\widebar{\cM}^q$ defined by $\bar{\ell}$, i.e., $\Z_1 \sim \Z_2$ whenever $\bar{\ell}(\Z_1) = \bar{\ell}(\Z_2)$, and $\cM^q:=\widebar{\cM}^q/\sim$ with metric $g_{[\Z]}$ is a Riemannian quotient manifold, where the superscript $q$ in $\cM^q$ stands for the ``quotient'', then $\bar{h}(\Z)$ induces a function $h([\Z])$ on the quotient manifold $\cM^q$ and \eqref{eq: general-strategy-example-embedded} can be transformed to an optimization on the quotient manifold $\mathcal{M}^q$: 
\begin{equation} \label{eq: general-strategy-example-quotient}
\text{(Opt. under Quotient Geometry)} \quad \min_{[\Z] \in \cM^q}h([\Z]).
\end{equation}

First, since $\bar{\ell}$ induces a diffeomorphism $\ell$ between $\cM^q$ and $\cM^e$ \cite[Proposition 3.5.23]{abraham2012manifolds}, i.e., $\cM^q$ and $\cM^e$ are diffeomorphic, a simple fact stated in the following lemma shows there is an equivalence relationship between the sets of local minimizers of \eqref{eq: general-strategy-example-embedded} and  \eqref{eq: general-strategy-example-quotient}.
\begin{lemma} \label{lm: local-minimizer-connection}
If $\X$ is a local minimizer  of \eqref{eq: general-strategy-example-embedded}, then $\ell^{-1}(\X)$ is a local minimizer  of \eqref{eq: general-strategy-example-quotient}; if $[\Z]$ is a local minimizer   of \eqref{eq: general-strategy-example-quotient}, then $\X = \ell([\Z])$ is a local minimizer of \eqref{eq: general-strategy-example-embedded}.
\end{lemma}

We show next that with a careful treatment with Riemannian gradients and Hessians of \eqref{eq: general-strategy-example-embedded} and \eqref{eq: general-strategy-example-quotient}, we are able to obtain a much richer geometric connection between the landscapes of these two problems.

\subsection{Outline of the Procedure} \label{sec: general-strategy-outline}

First, we find it is often relatively easy to connect the first-order geometries between \eqref{eq: general-strategy-example-embedded} and \eqref{eq: general-strategy-example-quotient}. For example, by taking derivative on both sides of $\bar{h}(\Z) = f(\bar{\ell}(\Z))$ along the direction $\theta_\Z \in \cH_\Z \widebar{\cM}^q$ and the chain rule, we have
\begin{equation} \label{eq: general-gradient-connection}
\begin{split}
	g_{[\Z]} (\grad \, h([\Z]), \theta_{[\Z]}) & \overset{ \eqref{eq: quotient-gradient-general} }  = \bar{g}_\Z \left( P_\Z^{\cH} (\grad\, \bar{h}(\Z)), \theta_\Z \right) \overset{(a)}  =  \bar{g}_\Z \left( \grad\, \bar{h}(\Z), \theta_\Z \right)\\
	& = \rmD \bar{h}(\Z)[\theta_\Z] = \rmD f(\bar{\ell}(\Z))[ \rmD \bar{\ell}(\Z)[\theta_\Z] ] \overset{(b)}= g_{\bar{\ell}(\Z)} \left( \grad f(\bar{\ell}(\Z)),  \rmD \bar{\ell}(\Z)[\theta_\Z] \right).
\end{split}
\end{equation} 
Here $g_{[\Z]}$ and $\bar{g}_\Z$ are Riemannian metrics on $\cM^q$ and $\widebar{\cM}^q$, respectively; (a) is because $\theta_\Z \in \cH_\Z \widebar{\cM}^q$; (b) is because $\rmD \bar{\ell}(\Z)[\theta_\Z]\in  T_{\bar{\ell}(\Z)} \cM^e$ \cite[Section 3.5.1]{absil2009optimization}. Thus, for reasonable choices of $\bar{\ell}$, we hope to find a connection between $\overline{\grad\, h([\Z])}$ and $\grad f(\bar{\ell}(\Z))$ based on \eqref{eq: general-gradient-connection}. 

Next, we present a general three-step procedure to connect the second-order geometries between \eqref{eq: general-strategy-example-embedded} and \eqref{eq: general-strategy-example-quotient}. 
\begin{itemize}[leftmargin=*]
	\item {\bf Step 1: Compute the quadratic forms of Riemannian Hessians.} We first compute the Riemannian Hessians from their definitions. For fixed-rank matrix optimization, we will derive the explicit expressions for the quadratic form of the Riemannian Hessians in Sections \ref{sec: connection-PSD} and \ref{sec: connection-general} and will see the quadratic forms of the Riemannian Hessians of $f(\X)$ and $h([\Z])$ always involve the quadratic form of the Riemannian Hessian of $\bar{f}$ and the Riemannian gradient of $\bar{f}$ or $\bar{h}$. For concreteness, we assume
	\begin{equation} \label{eq: general-setup-Hessian-expression}
		\begin{split}
			\Hess f(\X)[\xi_\X, \xi_\X] &= \Hess \bar{f} (\X)[\xi_\X, \xi_\X] +  \Psi_1  , \quad \forall \xi_\X \in T_\X \cM^e,\\
			\overline{\Hess\, h([\Z])}[\theta_\Z, \theta_\Z] &= \Hess \bar{f} (\bar{\ell}(\Z))[\varphi(\theta_\Z), \varphi(\theta_\Z)] + \Psi_2  ,\quad \forall \theta_\Z \in \cH_\Z \widebar{\cM}^q.
		\end{split}
	\end{equation} Here $\varphi: \cH_\Z \widebar{\cM}^q \to T_{\bar{\ell}(\Z)} \widebar{\cM}^e$ and $\cH_\Z \widebar{\cM}^q$ is the horizontal space of $T_\Z \widebar{\cM}^q$; $\Psi_1$ and $\Psi_2$ incorporate remaining terms of $\Hess f(\X)[\xi_\X, \xi_\X]$ and $\overline{\Hess\, h([\Z])}[\theta_\Z, \theta_\Z]$ other than the quadratic form of the Riemannian Hessian of $\bar{f}$. In the embedded submanifold setting, $\Psi_1$ has a closed form expression: $\Psi_1 = g_\X( \mathfrak{A}_\X ( \xi_\X, P_{T_\X \cM^e \perp} \grad \bar{f}(\X)),\xi_\X ) $, where $\mathfrak{A}_\X$ is the Weingarten map and $P_{T_\X \cM^e \perp}$ is the orthogonal complement of the projector $P_{T_\X \cM^e}$. See details in \cite{absil2013extrinsic}. In the quotient manifold setting, as we will see in Propositions \ref{prop: gradient-hessian-exp-PSD} and \ref{prop: gradient-hessian-exp-general}, we often have $\varphi(\theta_\Z) = \rmD \bar{\ell}(\Z)[\theta_\Z]$. 
	
	\item {\bf Step 2: Find a proper mapping $\cL$ between $\cH_\Z \widebar{\cM}^q$ and $T_\X \cM^e$ to connect Riemannian Hessians.
	} To establish the second-order geometric connection, i.e., the connection between $\Hess f(\X)$ and $\overline{\Hess\, h([\Z])}$ with $\X = \bar{\ell}(\Z)$, a natural idea is to first connect $\Hess \bar{f} (\X)[\xi_\X, \xi_\X]$ with $\Hess \bar{f} (\bar{\ell}(\Z))[\varphi(\theta_\Z), \varphi(\theta_\Z)]$. To do this, we would like to find a mapping $\cL$ between $\cH_\Z \widebar{\cM}^q$ and $T_\X \cM^e$ such that $\cL(\theta_\Z)=\varphi(\theta_\Z)$ and set $\xi_\X = \cL(\theta_\Z)$. 
	Moreover, $\cL$ is further constrained to be a bijection so that we can connect the whole spectra of $\Hess f(\X)$ with the spectra of $\overline{\Hess\, h([\Z])}$ as we will see in Step 3. We will see later that $\cL$ is often taken to be the restriction of $ \rmD \bar{\ell}(\Z)$ to $\cH_\Z \widebar{\cM}^q$. 
	
	On the other hand, such a mapping $\cL$ alone seems not enough for connecting $\Hess f(\X)[\xi_\X, \xi_\X]$ with $\overline{\Hess\, h([\Z])}[\theta_\Z, \theta_\Z]$ as $\Psi_1$ and $\Psi_2$ in \eqref{eq: general-setup-Hessian-expression} can be complex and distinct from each other (see the forthcoming Propositions \ref{prop: gradient-hessian-exp-PSD} and \ref{prop: gradient-hessian-exp-general}). This motivates us to put some assumptions on $\X$ or $[\Z]$ in order to proceed. In the following, we consider a simple setting to illustrate what assumptions may be needed. 
	\begin{example}\label{ex: example-1}
	Consider the optimization problem: $ \min_{x > 0} f(x) $, where $f$ is a scalar function defined on the positive part of the real line. If we consider the factorization $x = z^2$ and let $\bar{h}(z) = f(z^2)$, it is easy to see $\bar{h}(z) = \bar{h}(-z)$. So we can consider the equivalence classes $[z] = \{ z, -z \}$ and take the quotient manifold $\cM^q = (\bbR \setminus \{0\})/(+-)$, i.e., quotienting out the sign of the non-zero real number, as the search space. 
	Suppose the Euclidean inner product is adopted as the Riemannian metric. Then in both cases, the Riemannian Hessians of $f(x)$ and $h([z])$ are equal to the corresponding Euclidean Hessians (see also the forthcoming Proposition \ref{prop: gradient-hessian-exp-PSD}):
	\begin{equation} \label{eq: scalar-motivating-example}
		\begin{split}
		\Hess f(x)[\xi_x, \xi_x] = \xi^2_x f''(x); \quad \overline{\Hess\, h([z])}[\theta_z, \theta_z]	= \theta_z^2 \bar{h}''(z) = 2\theta_z^2f'(z^2) + 4 z^2\theta_z^2 f''(z^2),
		\end{split}
	\end{equation} where $f'$ and $f''$ denote the first and second derivatives of $f$. Given $z_1 \in \bbR, z_1 \neq 0$, $x_1 = z_1^2$ and \eqref{eq: scalar-motivating-example}, we see a natural assumption to connect $\Hess f(x_1)[\xi_{x_1}, \xi_{x_1}]$ with $\overline{\Hess\, h([z_1])}[\theta_{z_1}, \theta_{z_1}]$ is by assuming $z_1^2$ is a FOSP of $\min_{x > 0}f(x)$, i.e., $f'(z_1^2) = 0$, otherwise it is even not guaranteed that $\Hess f(x_1)[\xi_{x_1}, \xi_{x_1}]$ and $\overline{\Hess\, h([z_1])}[\theta_{z_1}, \theta_{z_1}]$ will have the same sign. At the same time, a more encouraging fact is that if we further let the bijection $\cL$ to be $\cL(\theta_z) = 2z \theta_z$ for any $\theta_z \in \bbR$, then we have $\overline{\Hess\, h([z_1])}[\theta_{z_1}, \theta_{z_1}]	= \Hess f(x_1)[\cL(\theta_{z_1}), \cL(\theta_{z_1})]$. 
	\end{example}
	Motivated by Example \ref{ex: example-1}, we find by putting some proper first-order assumptions on $\X$ or $[\Z]$, we could hope for a nice connection between $\Hess f(\X)$ and $\overline{\Hess\, h([\Z])}$. In fact, as we will see in Sections \ref{sec: connection-PSD} and \ref{sec: connection-general}, this intuition applies to all fixed-rank matrix optimization examples we consider.
	\item {\bf Step 3: Establish the spectra connection between Riemannian Hessians via bounding spectrum of $\cL$.} Suppose one has successfully worked through Steps 1 and 2 and come to the stage that $\overline{\Hess\, h([\Z])}[\theta_\Z, \theta_\Z] = \Hess f(\X)[\cL(\theta_\Z), \cL(\theta_\Z)]$ is shown to hold for any $\theta_\Z \in \cH_\Z \widebar{\cM}^q$ at properly chosen $[\Z]$ and $\X = \bar{\ell}(\Z)$. The following Theorem \ref{th: hessian-sandwich} shows a sandwich inequality between the spectra of $\Hess f(\X)$ and $\overline{\Hess\, h([\Z])}$ based on spectrum bounds of $\cL$ and its proof is provided in Appendix \ref{sec: additional-proof-general-strategy-in-connection}.
\end{itemize}

\begin{theorem}[Sandwich Inequalities for Spectra of Hessians] \label{th: hessian-sandwich}
	Suppose $\X \in \cM^e$, $\Z \in \widebar{\cM}^q$, $\dim(T_\X \cM^e) = \dim(\cH_\Z \widebar{\cM}^q) = d$. Then both $\Hess f(\X)$ and $\overline{\Hess\, h([\Z])}$ have $d$ eigenvalues. Moreover, if $\cL: \cH_\Z \widebar{\cM}^q \to T_\X \cM^e$ is a bijection satisfying
	\begin{equation} \label{eq: general-setting-hessian-connection}
		\overline{\Hess\, h([\Z])}[\theta_\Z, \theta_\Z] = \Hess f(\X)[\cL(\theta_\Z), \cL(\theta_\Z)], \quad \forall \theta_\Z \in \cH_\X \widebar{\cM}^q
	\end{equation} and
	\begin{equation} \label{eq: spectrum-bound-general-setting-L}
		\alpha \bar{g}_\Z (\theta_\Z, \theta_\Z) \leq g_\X (\cL(\theta_\Z),\cL(\theta_\Z)) \leq \beta \bar{g}_\Z (\theta_\Z, \theta_\Z), \quad \forall \theta_\Z \in \cH_\X \widebar{\cM}^q
	\end{equation} holds for some $0\leq \alpha \leq \beta$, then for $k = 1,\ldots,d$, we have
	$$\lambda_k( \overline{\Hess\, h([\Z])} ) \text{ is sandwiched between } \alpha \lambda_k( \Hess f(\X) ) \text{ and } \beta \lambda_k( \Hess f(\X)),$$
	where $\lambda_k( \Hess f(\X) )$ and $\lambda_k( \overline{\Hess\, h([\Z])} )$ are the $k$-th largest eigenvalues of $\Hess f(\X)$ and $\overline{\Hess\, h([\Z])}$, respectively.
\end{theorem}

\section{Embedded and Quotient Geometries on Fixed-rank Matrices} \label{sec: embedded-quotient-fixed-rank-matrix}

The set of $p$-by-$p$ rank $r$ PSD matrices $\cM_{r+}:=\left\{ 
\X\in \bbS^{p\times p} \mid {\rm rank}(\X) = r, \X \succcurlyeq 0
\right\}$ and the set of $p_1$-by-$p_2$ rank $r$ matrices $\cM_r:=\left\{ 
\X\in \bbR^{p_1\times p_2}\mid {\rm rank}(\X) = r
\right\}$ are two manifolds of particular interest. In Sections \ref{sec: embeddded-fixed-rank-matrix} and \ref{sec: quotient-fixed-rank-matrix}, we introduce embedded and quotient geometries on $\cM_{r+}$ and $\cM_r$, respectively.

\subsection{Embedded Geometries for $\cM_{r+}$ and $\cM_r$} \label{sec: embeddded-fixed-rank-matrix}
The following lemma shows that $\cM_{r+}$ and $\cM_r$ are smooth embedded submanifolds of $\bbR^{p \times p}$ and $\bbR^{p_1 \times p_2}$, respectively, and it also summarizes commonly used algebraic representations of the corresponding tangent spaces. To emphasizing the embedding natural of $\cM_{r+}$ and $\cM_r$, we write them as $\cM^e_{r+}$ and $\cM^e_{r}$, respectively.

\begin{lemma}{\rm(\cite[Chapter 5]{helmke2012optimization}, \cite[Proposition 5.2]{vandereycken2010riemannian}, \cite[Example 8.14]{lee2013smooth})}\label{lm: Mr+ and Mr manifold}
	$\cM^e_{r+}$ and $\cM^e_{r}$ are smooth embedded submanifolds of $\bbR^{p \times p}$ and $\bbR^{p_1 \times p_2}$ with dimensions $(pr - r(r-1)/2)$ and $(p_1 + p_2 -r)r$, respectively. The tangent space $T_{\X}\cM^e_{r+}$ at $\X \in \cM^e_{r+}$ is
	\begin{equation}
	\label{eq:tangent Mr+}
	T_\X \cM_{r+}^e = \left\{
	[\U\quad \U_{\perp}] \begin{bmatrix}
	\S & \D^\top \\[2pt]
	\D & \0
	\end{bmatrix}
	[\U\quad \U_{\perp}]^\top: \S \in \bbS^{r \times r}, \D \in \bbR^{(p-r) \times r}
	\right\},
	\end{equation} where $\U \in \st(r,p)$ spans the top $r$ eigenspace of $\X$.
	
	 The tangent space $T_{\X}\cM^e_{r}$ at $\X \in \cM^e_{r}$ is
	\begin{equation}
	\label{eq:tangent Mr}
	T_\X \cM_r^e = \left\{
	[\U\quad \U_{\perp}] \begin{bmatrix}
	\S & \D_2^\top \\[2pt]
	\D_1 & \0
	\end{bmatrix}
	[\V\quad \V_{\perp}]^\top: \S \in \bbR^{r \times r}, \D_1 \in \bbR^{(p_1 - r)\times r},  \D_2 \in \bbR^{(p_2 - r)\times r}
	\right\},
	\end{equation} where $\U \in \st(r,p_1)$ and $\V\in \st(r,p_2)$ span the left and right singular subspaces of $\X$, respectively.
\end{lemma} 
 In addition, we assume the embedded submanifolds $\mathcal{M}^e_{r+}$ and $\mathcal{M}^e_r$ are endowed with the natural metric induced by the Euclidean inner product, i.e., $\langle \U, \V \rangle = \mathrm{trace}(\U^\top \V)$.

\subsection{Quotient Geometries for $\cM_{r+}$ and $\cM_r$} \label{sec: quotient-fixed-rank-matrix}
The versatile choices of fixed-rank matrix factorization yield various Riemannian quotient structures and metrics, which have been explored in the literature on both $\cM_{r+}$ \cite{bonnabel2010riemannian,journee2010low,massart2020quotient,meyer2011regression} and $\cM_r$ \cite{absil2014two,meyer2011linear,mishra2012riemannian,mishra2014fixed}. Due to the factorization, the total space of fixed-rank matrices under the quotient geometry (i.e., focus of this subsection) can often be written as a product space of some simple smooth manifolds, including the following three examples to be used later:
\begin{itemize}
\item[(1)] $\bbR^{p \times r}_*$: the set of real $p$-by-$r$ full column rank matrices; 
\item[(2)] $\st(r,p)$: the set of real $p$-by-$r$ matrices with orthonormal columns, i.e., the Stiefel manifold; 
\item[(3)] $\bbS_{+}(r)$: the set of $r$-by-$r$ real symmetric positive definite matrices. 
\end{itemize}
All these three manifolds are smooth homogeneous spaces and there exists a smooth structure on their product space \cite[Section 3.1.6]{absil2009optimization}. In Table \ref{tab: basic-prop-simple-manifold}, we summarize several basic properties of these simple manifolds. 

\begin{table}[ht]
	\centering
	\begin{tabular}{c | c | c | c}
	\hline
	& $\bbR^{p \times r}_*$ & $\st(r,p)$ & $\bbS_{+}(r)$ \\
	\hline
	Dimension & $pr$ & $pr-(r^2 + r)/2$ & $(r^2 + r)/2$ \\	
	\hline
	Matrix & \multirow{2}{1em}{$\Y$} & \multirow{2}{1em}{$\U$} & \multirow{2}{1em}{$\B$} \\
	representation & & & \\ 
	\hline
	\multirow{3}{5em}{Tangent space} & \multirow{3}{7em}{$T_{\Y} \bbR^{p \times r}_* = \bbR^{p \times r}$} & \multirow{3}{12em}{$T_{\U}\st(r,p) = \{\U\bOmega + \U_\perp \D: \bOmega = - \bOmega^\top \in \bbR^{r \times r}, \D \in \bbR^{(p-r) \times r} \}$} & \multirow{3}{8em}{$T_{\B} \bbS_{+} (r) =  \bbS^{r \times r}$}\\
	& &  & \\
	& & & \\
	\hline
	\multirow{3}{7em}{Projection onto tangent space} & \multirow{3}{8em}{$P_{T_{\Y} \bbR^{p \times r}_*}(\eta_\Y) = \eta_\Y$,$\forall \eta_\Y \in \bbR^{p \times r}$} & \multirow{3}{12em}{$P_{T_{\U}\st(r,p)}(\eta_\U) = P_{\U_\perp}(\eta_\U) + \U \skew(\U^\top \eta_\U)$, $\forall \eta_\U \in \bbR^{p \times r}$} & \multirow{3}{8em}{$P_{T_{\B} \bbS_{+} (r)}(\eta_\B) = \sym(\eta_\B)$, $\forall \eta_\B \in \bbR^{r \times r}$}\\
	& & &\\
	& & &\\
	\hline
	Metric $g$ on & \multirow{2}{12em}{$g_\Y(\theta_\Y, \eta_\Y) = \tr(\W_\Y \theta_\Y^\top \eta_\Y)$} & \multirow{2}{12em}{$g_\U (\theta_\U, \eta_\U) =  \tr( \V_{\blacklozenge} \theta_\U^\top \eta_\U)$} & \multirow{2}{8em}{$g_\B (\theta_\B, \eta_\B) = \tr(\W_\B \theta_\B \W_\B \eta_\B)$}\\
	tangent space & &  & \\
	\hline
	\end{tabular}\caption{Basic Riemannian Geometric Properties for $\bbR^{p \times r}_*$, $\st(r,p)$ and $\bbS_{+}(r)$ \cite{absil2009optimization,edelman1998geometry,mishra2014fixed}. Here for any square matrix $\X$, $\sym(\X) = (\X + \X^\top)/2$ and $\skew(\X) =(\X - \X^\top)/2$; $\W_\Y$, $\V_\blacklozenge$ and $\W_\B$ are $r$-by-$r$ symmetric positive definite weight matrices that specifies the Riemannian metric $g$ (see Remark \ref{rem: W_Y-choice-in-metric} for discussions).}\label{tab: basic-prop-simple-manifold}
\end{table}

\begin{remark} \label{rem: W_Y-choice-in-metric}
	We introduce weight matrices $\W_\Y$, $\V_\blacklozenge$ and $\W_\B$ based on $\Y$, $\blacklozenge$ and $\B$, respectively, while defining $g_\Y$, $g_\U$ and $g_\B$ so that various metrics considered in literature are covered. $\W_\Y$, $\V_\blacklozenge$ and $\W_\B$ are required to be in  $\mathbb{S}_+(r)$ so that $g_\Y$, $g_\U$ and $g_\B$ are genuine Riemannian metrics \cite[Section 3.6]{absil2009optimization}, e.g., $g_\Y(\theta_\Y, \theta_\Y) \geq 0,  \forall \theta_\Y \in T_{\Y} \bbR^{p \times r}_*$. Common choices of $\W_\Y$ include $\W_\Y = \I_r$ (flat metric) \cite{journee2010low}, $\W_\Y = (\Y^\top \Y)^{-1}$ (right-invariant metric) \cite{meyer2011linear}, and $\W_\Y = \Y^\top \Y$ \cite{mishra2012riemannian,huang2017solving}. Common choice for $\V_\blacklozenge$ and $\W_\B$ is $\V_\blacklozenge = \I$ and $\W_\B = \B^{-1}$ \cite{mishra2013low}. Finally, we note $\V_{\blacklozenge}$ is introduced with the to-be-specified matrix ``$\blacklozenge$'', this ``$\blacklozenge$'' matrix can often be other matrix factors in the fixed-rank matrix factorization; see the following Section \ref{sec: quotient-PSD} on $\cM^{q_2}_{r+}$, Section \ref{sec: quotient-general} on $\cM^{q_3}_{r}$ and Tables \ref{tab: illustraion-gap-coefficient-PSD} and \ref{tab: illustraion-gap-coefficient-general} for specific examples. To the best of our knowledge, such metric for $\st(r,p)$ is non-standard and has not been explored in the literature and we will see in Sections \ref{sec: connection-PSD}, \ref{sec: connection-general} and Remark \ref{rem: algorithmic-connection} that it is important to have this flexibility to reveal geometric and algorithmic connections between embedded and quotient geometries in fixed-rank matrix optimization. 
\end{remark}

Next, we present two quotient geometries for fixed-rank PSD matrices (based on full-rank factorization and polar factorization) and three quotient geometries for fixed-rank general matrices (based on full-rank factorization, polar factorization, and subspace-projection factorization) in Sections \ref{sec: quotient-PSD} and \ref{sec: quotient-general}, respectively. These quotient geometries have been explored in \cite{bonnabel2010riemannian,journee2010low,meyer2011linear,meyer2011regression,mishra2012riemannian,mishra2014fixed}. 
Here, we provide a unified way to characterize these quotient geometries, e.g., the introduction of $\W_\Y, \V_\blacklozenge$ and $\W_\B$ in the Riemannian metric, in both PSD and general cases, propose several new representations for the horizontal spaces (see the forthcoming Lemmas \ref{lm: general-quotient-manifold2-prop}, \ref{lm: general-quotient-manifold3-prop} and the discussion afterward), derive explicit formulas for the Riemannian Hessians under these quotient geometries (see the forthcoming Remark \ref{rem: closed-form-riemannian-hessian-quotient}) and show their geometric connections to the embedded geometry in fixed-rank matrix optimization.

\subsubsection{Quotient Geometries for $\cM_{r+}$} \label{sec: quotient-PSD}
Suppose $\X \in \bbS^{p \times p}$ is a rank $r$ PSD matrix with economic eigendecomposition $\X = \U' \bSigma' \U^{'\top} $. 
\vskip.1cm
{\bf (1) Full-rank Factorization $\mathcal{M}_{r+}^{q_1}$.} In this factorization, we view $\X$ as $\X = \Y \Y^\top$ for $\Y \in \bbR^{p \times r}_*$. Such a factorization exists, e.g., $\Y = \U' \bSigma'^{1/2}$, but is not unique because of the invariance mapping $\Y \mapsto \Y \O$ for any $\O \in \bbO_r$. To cope with it, we encode the invariance mapping in an abstract search space by defining the equivalence classes $[\Y] = \{ \Y \O: \O \in \bbO_r \}$. Since the invariance mapping is performed via the Lie group $\bbO_r$ smoothly, freely and properly, we have $\cM_{r+}^{q_1} := \widebar{\cM}_{r+}^{q_1}/\bbO_r$, where $\widebar{\cM}_{r+}^{q_1} = \bbR^{p \times r}_*$, is a quotient manifold of $\widebar{\cM}_{r+}^{q_1}$ \cite[Theorem 21.10]{lee2013smooth}. We equip $T_\Y \widebar{\cM}_{r+}^{q_1}$ with the metric $\bar{g}^{r+}_\Y(\eta_\Y, \theta_\Y)= \tr(\W_\Y \eta_\Y^\top \theta_\Y)$ as given in Table \ref{tab: basic-prop-simple-manifold}. In the following Lemma \ref{lm: psd-quotient-manifold1-prop}, we provide the corresponding vertical and horizontal spaces of $T_\Y \widebar{\cM}_{r+}^{q_1}$ and show under some proper assumptions on $\W_\Y$, $\cM_{r+}^{q_1}$ is a Riemannian quotient manifold endowed with the Riemannian metric $g^{r+}_{[\Y]}$ induced from $\bar{g}^{r+}_\Y$.
\begin{lemma} \label{lm: psd-quotient-manifold1-prop} (i) Given $\U \in \st(r,p)$ that spans the top $r$ eigenspace of $\Y\Y^\top$ and $\P = \U^\top \Y$, the vertical and horizontal spaces of $T_\Y \widebar{\cM}_{r+}^{q_1}$ are given as follows:
	\begin{equation*}
		\begin{split}
			\cV_\Y \widebar{\cM}_{r+}^{q_1} &= \{ \theta_\Y: \theta_\Y = \Y \bOmega, \bOmega = - \bOmega^\top \in \bbR^{r \times r} \} = \{ \theta_\Y: \theta_\Y = \U \bOmega \P^{-\top}, \bOmega = - \bOmega^\top \in \bbR^{r \times r} \},\\
			\cH_\Y \widebar{\cM}_{r+}^{q_1} &= \{ \theta_\Y: \theta_\Y = (\U \S + \U_\perp \D) \P^{-\top}, \S\P^{-\top}\W_\Y \P^{-1} \in \bbS^{r \times r}, \D \in \bbR^{(p-r) \times r} \},
		\end{split}
	\end{equation*} with $\dim(\cV_\Y \widebar{\cM}_{r+}^{q_1}) = (r^2-r)/2$, $\dim(\cH_\Y \widebar{\cM}_{r+}^{q_1}) = pr-(r^2-r)/2$ and $\cV_\Y \widebar{\cM}_{r+}^{q_1} \perp \cH_\Y \widebar{\cM}_{r+}^{q_1}$ with respect to $\bar{g}^{r+}_\Y$.
	
	(ii) Moreover, $\cM_{r+}^{q_1}$ is a Riemannian quotient manifold endowed with the metric $g_{[\Y]}^{r+}$ induced from $\bar{g}_\Y^{r+}$ if and only if $\W_\Y = \O \W_{\Y\O} \O^\top$ holds for any $\O \in \bbO_r$.
\end{lemma}

\vskip.1cm

{\bf (2) Polar Factorization $\mathcal{M}_{r+}^{q_2}$.} We factorize $\X = \U'\O (\O^\top \bSigma' \O) (\U'\O)^\top = \U \B \U^\top$ with $\O \in \bbO_r, \U \in \st(r,p)$, $\B \in \bbS_+(r)$ \cite{bonnabel2010riemannian}. Due to the rotational invariance $(\U, \B) \mapsto (\U\O, \O^\top \B \O)$ for $\O \in \bbO_r$, we define the search space as the following equivalence classes $[\U, \B] = \{(\U \O, \O^\top \B \O): \O \in \bbO_r \}$. This results in the second quotient manifold: $\cM_{r+}^{q_2}:= \widebar{\cM}_{r+}^{q_2}/\bbO_r$, where $\widebar{\cM}_{r+}^{q_2} = \st(r,p) \times \bbS_+(r)$. By taking the metrics for $\st(r,p)$ and $\bbS_+(r)$ given in Table \ref{tab: basic-prop-simple-manifold} with the notation $\V_\blacklozenge$ replaced by $\V_\B$, we endow $\widebar{\cM}_{r+}^{q_2}$ with the metric $\bar{g}_{(\U,\B)}^{r+}(\eta_{(\U, \B)}, \theta_{(\U, \B)}) = \tr(\V_\B \eta_U^\top \theta_U) + \tr(\W_\B \eta_B \W_\B \theta_B)$ for $\eta_{(\U, \B)} = [\eta_U^\top \quad \eta_B^\top]^\top, \theta_{(\U, \B)} = [\theta_U^\top \quad \theta_B^\top]^\top \in T_{(\U,\B)} \widebar{\cM}_{r+}^{q_2}$. 
\begin{lemma} \label{lm: psd-quotient-manifold2-prop} 
\label{lm:BS10-Thm1} (i) The vertical and horizontal spaces of $T_{(\U,\B)} \widebar{\cM}_{r+}^{q_2}$ are given as:
\begin{equation*}
\begin{split}
	\cV_{(\U,\B)} \widebar{\cM}_{r+}^{q_2} &= \{ \theta_{(\U,\B)} =[\theta_U^\top \quad \theta_B^\top]^\top: \theta_U = \U \bOmega, \theta_B = \B \bOmega - \bOmega \B, \bOmega = - \bOmega^\top \in \bbR^{r \times r} \},\\
			\cH_{(\U,\B)} \widebar{\cM}_{r+}^{q_2} &= \{ \theta_{(\U,\B)} =[\theta_U^\top \quad \theta_B^\top]^\top: \theta_U =  \U_\perp \D, \theta_B \in \bbS^{r \times r}, \D \in \bbR^{(p-r) \times r}\}.
\end{split}		
\end{equation*} 
Here, $\dim(\cV_{(\U,\B)} \widebar{\cM}_{r+}^{q_2}) = (r^2-r)/2$, $\dim(\cH_{(\U,\B)} \widebar{\cM}_{r+}^{q_2}) = pr-(r^2-r)/2$.

(ii) Moreover, $\cM_{r+}^{q_2}$ is a Riemannian quotient manifold endowed with the metric $g^{r+}_{[\U,\B]}$ induced from $\bar{g}^{r+}_{(\U,\B)}$ if and only if $\V_\B = \O \V_{\O^\top \B \O} \O^\top$ and $\W_\B = \O \W_{\O^\top \B \O} \O^\top$ for any $\O \in \bbO_r$.
\end{lemma}

\begin{remark}({\bf $\cV_{(\U,\B)} \widebar{\cM}_{r+}^{q_2}$ and $\cH_{(\U,\B)} \widebar{\cM}_{r+}^{q_2}$ are not orthogonal})
\label{rem: choice-horizontal-space}
In the context of Lemma \ref{lm:BS10-Thm1}, $\cH_{(\U,\B)} \widebar{\cM}_{r+}^{q_2}$ is complementary but not orthogonal to $\cV_{(\U,\B)} \widebar{\cM}_{r+}^{q_2}$, which means $\cH_{(\U,\B)} \widebar{\cM}_{r+}^{q_2}$ is not a canonical horizontal space for the quotient manifold $\mathcal{M}_{r+}^{q_2}$. In fact, it is easier to find a correspondence between this non-canonical horizontal space $\cH_{(\U,\B)} \widebar{\cM}_{r+}^{q_2}$ and the tangent space $T_{\X}\cM_{r+}^e$ under the embedded geometry (see later Proposition \ref{prop: psd-bijection2}). Such a correspondence is critical in establishing the landscape connections of embedded and quotient geometries in Riemannian optimization as we have illustrated in Step 2 of the general procedure in Section \ref{sec: general-strategy-for-connection}.	
\end{remark}

Finally, the following Lemma \ref{lm: completeness-quotient-PSD} shows the full-rank and polar factorizations fully take into account the invariance in the corresponding matrix factorizations and its proof is provided in Appendix \ref{proof-sec: embedded-quotient-fixed-rank-matrix}. 
\begin{lemma}\label{lm: completeness-quotient-PSD}~
\begin{itemize}[leftmargin=*]
	 \item (Full-rank Factorization) Let $\Y_1, \Y_2 \in \bbR_*^{p \times r}$. Then $\Y_1 \Y_1^\top = \Y_2 \Y_2^\top$ if and only if $\Y_2 = \Y_1 \O$ for some $\O \in \bbO_r$.
	 \item (Polar Factorization) Let $\U_1, \U_2 \in \st(r,p)$ and $\B_1, \B_2 \in \bbS_+(r)$. Then $\U_1 \B_1 \U_1^\top = \U_2 \B_2 \U_2^\top$ if and only if $\U_2 = \U_1 \O $ and $\B_2 = \O^\top \B_1 \O$ for some $\O \in \bbO_r$.
\end{itemize}Moreover, $\cM_{r+} = \{\Y \Y^\top: \Y \in \bbR^{p \times r}_* \} = \{\U \B \U^\top: \U \in \st(r,p), \B \in \bbS_+(r) \}$.
\end{lemma}

In Table \ref{tab: basic-prop-quotient-psd}, we summarize the basic properties of the full-rank factorization and polar factorization based quotient manifolds on fixed-rank PSD matrices.
\begin{table}[ht]
	\centering
	\begin{tabular}{c | c | c }
	\hline
	& $\cM_{r+}^{q_1}$ & $\cM_{r+}^{q_2}$ \\
	\hline
	Matrix representation & $\Y$ & $(\U,\B)$ \\
	\hline
		Equivalence classes & $[\Y] = \{\Y \O: \O \in \bbO_r \}$ & $[\U, \B]=\{(\U\O, \O^\top \B \O), \O \in \bbO_r \} $ \\
	\hline
	Total space $\widebar{\cM}_{r+}$ &$ \bbR^{p \times r}_* $ & $\st(r,p)\times \bbS_{+} (r)$ \\
	\hline
	Tangent space in & \multirow{2}{3em}{$T_\Y \bbR_*^{p \times r}$} & \multirow{2}{10em}{ $T_\U \st(r,p) \times T_\B \bbS_+(r)$ } \\
	total space & & \\
	\hline
	Metric $\bar{g}^{r+}$ on & \multirow{2}{7em}{$\tr(\W_\Y \eta_\Y^\top \theta_\Y)$, $\W_\Y \in \bbS_+(r)$} & \multirow{2}{15em}{$  \tr(\V_\B\eta_U^\top \theta_U) + \tr(\W_\B \eta_B \W_\B \theta_B)$, $\V_\B,\W_\B \in \bbS_+(r)$} \\
	total space & &  \\
	\hline
	\end{tabular}\caption{Basic Properties of Quotient Manifolds $\cM_{r+}^{q_1}$, $\cM_{r+}^{q_2}$.} \label{tab: basic-prop-quotient-psd}
\end{table}

\subsubsection{Quotient Geometries for $\cM_r$} \label{sec: quotient-general}
In this section, we introduce three quotient structures for $\cM_r$ based on three fixed-rank matrix factorizations. For $\X \in \cM_r$, denote the SVD as $\X = \U' \bSigma' \V^{'\top}$. 
\vskip.1cm
{\bf (1) Full-rank Factorization $\mathcal{M}_{r}^{q_1}$.} In this factorization, we rearrange the SVD of $\X$ as $\X = (\U' \bSigma^{'1/2}) (\bSigma^{'1/2} \V^{'\top}) = \L \R^\top $, where $\L \in \bbR^{p_1 \times r}_*$, $\R \in \bbR^{p_2 \times r}_*$. The first quotient geometry for $\cM_r$ results from the invariance mapping $(\L, \R) \mapsto (\L \M, \R\M^{-\top})$ for $\M \in \GL(r)$, where $\GL(r) := \{ \M \in \bbR^{r \times r}: \M \text{ is invertible}\}$ denotes the degree $r$ general linear group. It is thus straightforward to consider the equivalence classes $[\L, \R]=\{ (\L\M, \R \M^{-\top}): \M \in \GL(r) \}$ as the search space. The set of equivalence classes form the quotient manifold $\cM_{r}^{q_1}:= \widebar{\cM}_r^{q_1}/\GL(r)$, where $\widebar{\cM}_r^{q_1} = \bbR_*^{p_1 \times r} \times \bbR_*^{p_2 \times r}$, as $\GL(r)$ is a Lie group \cite[Theorem 21.10]{lee2013smooth}. Suppose $\W_{\L,\R}, \V_{\L,\R}$ are two $r$-by-$r$ positive definite matrices depending on $\L$ and $\R$, the metric we endow on $T_{(\L, \R)} \widebar{\cM}_{r}^{q_1}$ is $\bar{g}_{(\L,\R)}^r ( \eta_{(\L,\R)}, \theta_{(\L,\R)} ) = \tr( \W_{\L,\R} \eta_L^\top \theta_R ) + \tr(\V_{\L, \R} \eta_R^\top \theta_R)$ for $\eta_{(\L, \R)} = [ \eta_L^\top \quad \eta_R^\top ]^\top, \theta_{(\L, \R)} = [ \theta_L^\top \quad \theta_R^\top ]^\top \in T_{(\L, \R)} \widebar{\cM}_{r}^{q_1}$. In the following Lemma \ref{lm: general-quotient-manifold1-prop}, we show with some proper assumptions on $\W_{\L,\R}$ and $ \V_{\L, \R}$, $\cM_{r}^{q_1}$ is a Riemannian quotient manifold.
\begin{lemma} \label{lm: general-quotient-manifold1-prop}
 (i) Suppose $\U \in \st(r,p_1)$ and $ \V \in \st(r,p_2)$ span the top $r$ left and right singular subspaces of $\L\R^\top$, respectively and $\P_1 = \U^\top \L, \P_2 = \V^\top \R$. Then the vertical and horizontal spaces of $T_{(\L,\R)} \widebar{\cM}_{r}^{q_1}$ are given as follows:
	\begin{equation*}
		\begin{split}
			\cV_{(\L,\R)} \widebar{\cM}_{r}^{q_1} &= \left\{ \theta_{(\L,\R)}=[\theta_L^\top \quad \theta_R^\top]^\top : \theta_L = \U \S \P_2^{-\top}, \theta_R = -\V \S^\top \P_1^{-\top}, \S \in \bbR^{r \times r} \right\},\\
			\cH_{(\L,\R)} \widebar{\cM}_{r}^{q_1} &= \left\{ \theta_{(\L,\R)}= \begin{bmatrix}
				\theta_L\\
				\theta_R
			\end{bmatrix}: \begin{array}{l}
				\theta_L = (\U \S \P_2\W_{\L, \R}^{-1} \P_2^\top + \U_\perp \D_1 ) \P_2^{-\top}, \D_1 \in \bbR^{(p_1-r)\times r}\\
				\theta_R = (\V \S^\top \P_1 \V_{\L,\R}^{-1} \P_1^\top + \V_\perp \D_2  ) \P_1^{-\top}, \D_2 \in \bbR^{(p_2-r)\times r}
			\end{array}, \S \in \bbR^{r \times r}   \right\},
		\end{split}
	\end{equation*} with $\dim(\cV_{(\L,\R)} \widebar{\cM}_{r}^{q_1}) = r^2$, $\dim(\cH_{(\L,\R)} \widebar{\cM}_{r}^{q_1}) = (p_1 + p_2-r)r$ and $\cV_{(\L,\R)} \widebar{\cM}_{r}^{q_1} \perp \cH_{(\L,\R)} \widebar{\cM}_{r}^{q_1}$ with respect to $\bar{g}^{r}_{(\L,\R)}$. 
	
	(ii) Moveover, $\cM_{r}^{q_1}$ is a Riemannian quotient manifold endowed with metric $g_{[\L,\R]}^r$ induced from $\bar{g}_{(\L,\R)}^r$ if and only if $\W_{\L,\R} = \M \W_{\L\M,\R \M^{-\top}} \M^\top$ and $\V_{\L,\R} = \M^{-\top } \V_{\L\M,\R \M^{-\top}} \M^{-1}$ hold for any $\M \in \GL(r)$.
\end{lemma}
\begin{remark} \label{rem: choise of W-L-R}
	Similarly to $\W_\Y$ for the PSD case, we introduce $\W_{\L,\R}$ and $ \V_{\L,\R}$ in $\bar{g}_{(\L,\R)}^r$ to accommodate various metric choices considered in literature. Common choices include: $\W_{\L,\R} = (\L^\top \L)^{-1}, \V_{\L,\R} = (\R^\top \R)^{-1}$ \cite{meyer2011linear,mishra2014fixed} and $\W_{\L,\R} = \R^\top \R, \V_{\L,\R} = \L^\top \L$ \cite{mishra2012riemannian}. Distinct from the quotient geometry on $\cM^{q_1}_{r+}$, the flat metric, i.e., the Euclidean one with $\W_{\L,\R} = \I_r$ and $\V_{\L,\R} = \I_r$, is no longer proper as it does not yield a valid Riemannian metric in the quotient space.
\end{remark}

\vskip.1cm
{\bf (2) Polar Factorization $\mathcal{M}_{r}^{q_2}$.} We consider another factorization: $\X = \U'\O (\O^\top \bSigma' \O  ) (\V' \O)^\top = \U \B \V^\top$, where $\O \in \bbO_r, \U \in \st(r,p_1)$, $\B \in \bbS_+(r), \V \in \st(r,p_2)$. The rotational invariance mapping here is $(\U, \B, \V) \mapsto (\U \O, \O^\top \B \O, \V \O)$ for $\O \in \bbO_r$. This gives us the equivalence classes $[\U, \B, \V] = \{ (\U \O, \O^\top \B \O, \V \O): \O \in \bbO_r \}$ and the second quotient manifold $\cM_r^{q_2} = \widebar{\cM}_{r}^{q_2}/\bbO_r$, where $\widebar{\cM}_{r}^{q_2} = \st(r,p_1) \times \bbS_+(r) \times \st(r,p_2)$. We pick metrics for $\st(r,p)$ and $\bbS_+(r)$ given in Table \ref{tab: basic-prop-simple-manifold} with ``$\blacklozenge$'' as $\B$, and endow $T_{(\U,\B, \V)} \widebar{\cM}_{r}^{q_2}$ with the product metric $\bar{g}_{(\U,\B,\V)}^r ( \eta_{(\U,\B,\V)}, \theta_{(\U,\B,\V)} ) = \tr(\V_\B \eta_U^\top \theta_U ) + \tr(\W_\B \eta_B \W_\B \theta_B ) + \tr(\V_\B \eta_V^\top \theta_V) $ for $\eta_{(\U,\B,\V)} = [ \eta_U^\top \quad \eta_B^\top \quad \eta_V^\top ]^\top, \theta_{(\U,\B,\V)} = [ \theta_U^\top \quad \theta_B^\top \quad \theta_V^\top ]^\top \in T_{(\U,\B,\V)} \widebar{\cM}_{r}^{q_2}$. For technical convenience, we focus on the case $\V_\B = \I_r, \W_\B = \B^{-1}$, which is the most standard metric under this geometry \cite{mishra2013low,mishra2014fixed}. In the following Lemma \ref{lm: general-quotient-manifold2-prop}, we show $\cM_{r}^{q_2}$ is a Riemannian quotient manifold and provide expressions for its vertical and horizontal spaces.
\begin{lemma}  \label{lm: general-quotient-manifold2-prop}
	$\cM_{r}^{q_2}$ endowed with the metric $g^{r}_{[\U,\B,\V]}$ induced from $\bar{g}^{r}_{(\U,\B,\V)}$ is a Riemannian quotient manifold with the following vertical and horizontal spaces:
	\begin{equation} \label{eq: vertical-horizontal-quotient-general-manifold2}
		\begin{split}
				\cV_{(\U,\B,\V)} \widebar{\cM}_{r}^{q_2} &= \left\{ \theta_{(\U,\B,\V)} =[\theta_U^\top \quad \theta_B^\top \quad \theta_V^\top]^\top: \begin{array}{c}
					\theta_U = \U \bOmega,\theta_B = \B \bOmega - \bOmega \B,\\
					\theta_V = \V \bOmega,\bOmega = - \bOmega^\top \in \bbR^{r \times r}
				\end{array} \right\},\\
			\cH_{(\U,\B,\V)} \widebar{\cM}_{r}^{q_2} &= \left\{ \theta_{(\U,\B,\V)} =[\theta_U^\top \quad \theta_B^\top \quad \theta_V^\top]^\top: \begin{array}{c}
				\theta_U =\U_\perp \D_1 + \U \bOmega, \theta_B \in \bbS^{r \times r},\theta_V =\V_\perp \D_2 - \V \bOmega, \\
				\bOmega = - \bOmega^\top \in \bbR^{r \times r} , \D_1 \in \bbR^{(p_1-r) \times r}, \D_2 \in \bbR^{(p_2-r) \times r}
			\end{array} \right\},
		\end{split}
	\end{equation}with $\dim(\cV_{(\U,\B,\V)} \widebar{\cM}_{r}^{q_2}) = (r^2-r)/2$, $\dim(\cH_{(\U,\B,\V)} \widebar{\cM}_{r}^{q_2}) = (p_1 + p_2 - r)r$.
\end{lemma}

\vskip.1cm
{\bf (3) Subspace-projection Factorization $\mathcal{M}_{r}^{q_3}$.} The third quotient geometry is based on the factorization $\X = \U' (\bSigma' \V^{'\top}) = \U \Y^\top$, where $\U \in \st(r,p_1)$ and $\Y \in \bbR^{p_2 \times r}_*$. This factorization is called the subspace-projection factorization \cite{mishra2014fixed} as $\U$ represents the column space of $\X$ and $\Y$ is the left projection coefficient matrix of $\X$ on $\U$. Here the rotational invariance mapping is $(\U, \Y) \mapsto (\U\O, \Y\O)$ for $\O \in \bbO_r$ and the equivalence classes are $[\U, \Y] = \{ (\U\O, \Y \O): \O \in \bbO_r \}$. This results in the third quotient manifold we are interested in: $\cM_r^{q_3}:= \widebar{\cM}_{r}^{q_3}/\bbO_r$, where $\widebar{\cM}_{r}^{q_3} = \st(r,p_1) \times \bbR_*^{p_2 \times r}$. By taking the metric on $\st(r,p_1)$ with ``$\blacklozenge$'' in $\V_\blacklozenge$ as $\Y$, we endow $\widebar{\cM}_{r}^{q_3}$ with the metric $\bar{g}_{(\U,\Y)}^{r}(\eta_{(\U, \Y)}, \theta_{(\U, \Y)}) = \tr(\V_\Y \eta_U^\top \theta_U) + \tr(\W_\Y \eta_Y^\top \theta_Y)$ for $\eta_{(\U, \Y)} = [\eta_U^\top \quad \eta_Y^\top]^\top, \theta_{(\U, \Y)} = [\theta_U^\top \quad \theta_Y^\top]^\top \in T_{(\U,\Y)} \widebar{\cM}_{r}^{q_3}$. In the following Lemma \ref{lm: general-quotient-manifold3-prop}, we provide the vertical and horizontal spaces of $T_{(\U,\Y)} \widebar{\cM}_{r}^{q_3}$ and show with some proper assumptions on $\V_\Y$ and $\W_{\Y}$, $\cM_{r}^{q_3}$ is a Riemannian quotient manifold.
\begin{lemma} \label{lm: general-quotient-manifold3-prop} (i)
	The vertical and horizontal spaces of $T_{(\U,\Y)} \widebar{\cM}_{r}^{q_3}$ are
	\begin{equation} \label{eq: vertical-horizontal-quotient-general-manifold3}
		\begin{split}
			\cV_{(\U,\Y)} \widebar{\cM}_{r}^{q_3} &= \{ \theta_{(\U,\Y)} = [\theta_U^\top \quad \theta_Y^\top]^\top: \theta_U = \U \bOmega, \theta_Y = \Y \bOmega, \bOmega = - \bOmega^\top \in \bbR^{r \times r} \},\\
			\cH_{(\U,\Y)} \widebar{\cM}_{r}^{q_3} &= \{ \theta_{(\U,\Y)} =[\theta_U^\top \quad \theta_Y^\top]^\top: \theta_U =  \U_\perp \D, \theta_Y \in \bbR^{p_2 \times r}, \D \in \bbR^{(p_1-r) \times r} \}.
		\end{split}
	\end{equation}
	Here, $\dim(\cV_{(\U,\Y)} \widebar{\cM}_{r}^{q_3}) = (r^2-r)/2$, $\dim(\cH_{(\U,\Y)} \widebar{\cM}_{r}^{q_3}) = (p_1 + p_2 - r)r$.
	
	(ii) Moreover, $\cM_{r}^{q_3}$ is a Riemannian quotient manifold endowed with metric $g_{[\U,\Y]}^r$ induced from $\bar{g}_{(\U,\Y)}^r$ if and only if $\V_\Y = \O \V_{\Y\O} \O^\top$ and $\W_\Y = \O \W_{\Y\O} \O^\top$ hold for any $\O \in \bbO_r$.
\end{lemma}
\begin{remark}
In Lemmas \ref{lm: general-quotient-manifold2-prop} and \ref{lm: general-quotient-manifold3-prop}, we introduce new horizontal spaces that are distinct from the canonical horizontal spaces in the literature in $\cM_r^{q_2}$ and $\cM_r^{q_3}$ \cite{absil2014two,meyer2011linear,mishra2013low,mishra2014fixed}. These new horizontal spaces admit closed-form expressions, which makes developing a correspondence between the embedded manifold tangent space $T_\X \cM^e_r$ and the non-canonical horizontal spaces $\cH_{(\U,\B,\V)} \widebar{\cM}_{r}^{q_2}$, $\cH_{(\U,\Y)} \widebar{\cM}_{r}^{q_3}$ easier, and facilitates the later landscape analysis (see later in Propositions \ref{prop: general-bijection2} and \ref{prop: general-bijection3}). 
\end{remark}

 The following Lemma \ref{lm: completeness-quotient-general} shows the above full-rank, polar, and subspace-projection factorizations fully take into account the invariance in the corresponding quotient geometries and its proof is provided in Appendix \ref{proof-sec: embedded-quotient-fixed-rank-matrix}.
\begin{lemma} \label{lm: completeness-quotient-general}~
\begin{itemize}[leftmargin=*]
	 \item (Full-rank Factorization) Let $\L_1, \L_2 \in \bbR_*^{p_1 \times r}$ and $\R_1, \R_2 \in \bbR_*^{p_2 \times r}$. Then $\L_1 \R_1^\top = \L_2 \R_2^\top$ if and only if $\L_2 = \L_1 \M$ and $\R_2 = \R_1 \M^{-\top} $ for some $\M \in \GL(r)$.
	 \item (Polar Factorization) Let $\U_1, \U_2 \in \st(r,p_1)$, $\B_1, \B_2 \in \bbS_+(r)$ and $\V_1, \V_2 \in \st(r,p_2)$. Then $\U_1 \B_1 \V_1^\top = \U_2 \B_2 \V_2^\top$ if and only if $\U_2 = \U_1 \O $, $\B_2 = \O^\top \B_1 \O$ and $\V_2 = \V_1 \O$ for some $\O \in \bbO_r$.
	 \item (Subspace-projection Factorization) Let $\U_1, \U_2 \in \st(r,p_1)$ and $\Y_1, \Y_2 \in \bbR_*^{p_2 \times r}$. Then $\U_1 \Y_1^\top = \U_2 \Y_2^\top$ if and only if $\U_2 = \U_1 \O $ and $\Y_2 = \Y_1 \O$ for some $\O \in \bbO_r$.
\end{itemize}Finally, $\cM_{r} = \{\L \R^\top: \L \in \bbR^{p_1 \times r}_*, \R \in \bbR_*^{p_2 \times r} \} = \{\U \B \V^\top: \U \in \st(r,p_1), \B \in \bbS_+(r), \V \in \st(r, p_2) \} =\{\U \Y^\top: \U \in \st(r,p_1), \Y \in \bbR_*^{p_2 \times r} \} $.
\end{lemma}

In Table \ref{tab: basic-prop-quotient-general}, we summarize the basic properties of the full-rank factorization, polar factorization, and subspace-projection factorization based quotient manifolds of the general fixed-rank matrices.
\begin{table}[ht]
\centering
\begin{tabular}{c | c | c| c}
		\hline
	& $\cM_{r}^{q_1}$ & $\cM_{r}^{q_2}$ & $\cM_{r}^{q_3}$\\
	\hline
	Matrix & \multirow{2}{3em}{$(\L,\R)$} & \multirow{2}{4em}{$(\U, \B,\V)$} & \multirow{2}{3em}{$(\U, \Y)$}\\
	representation & & & \\
	\hline
	\multirow{3}{4em}{Equivalence classes} & \multirow{3}{9em}{$[\L,\R] = \{ (\L\M, \R \M^{-\top}): \M \in \GL(r) \}$} & \multirow{3}{14em}{$[\U, \B, \V] = \{ (\U \O, \O^\top \B \O, \V \O): \O \in \bbO_r \}$} & \multirow{3}{10em}{$[\U, \Y] = \{ (\U\O, \Y \O): \O \in \bbO_r \}$} \\
	& & &\\
	& & & \\
	\hline 
	Total space & \multirow{2}{9em}{$\bbR_*^{p_1 \times r} \times \bbR_*^{p_2 \times r}$} & \multirow{2}{14em}{$ \st(r,p_1) \times \bbS_+(r) \times \st(r,p_2)$} & \multirow{2}{9em}{$ \st(r,p_1) \times \bbR_*^{p_2 \times r}$} \\
	$\widebar{\cM}_r$& & &\\
	\hline
	Tangent space & \multirow{2}{9em}{$T_{\L} \bbR^{p_1 \times r}_* \times T_{\R} \bbR^{p_2 \times r}_*$} & \multirow{2}{14em}{$ T_\U \st(r,p_1) \times T_\B \bbS_+(r) \times T_\V \st(r,p_2)$} & \multirow{2}{10em}{$T_\U \st(r,p_1) \times T_\Y \bbR^{p_2 \times r}_*$} \\
	in total space & & & \\
	\hline 
	\multirow{4}{5em}{Metric $\bar{g}^r$ on total space} & \multirow{4}{10em}{$ \tr( \W_{\L,\R} \eta_L^\top \theta_L ) + \tr(\V_{\L, \R} \eta_R^\top \theta_R)$, $\W_{\L,\R} \in \bbS_+(r), \V_{\L, \R} \in \bbS_+(r)$} & \multirow{4}{14em}{$ \tr( \eta_U^\top \theta_U ) + \tr(\B^{-1} \eta_B \B^{-1} \theta_B ) + \tr(\eta_V^\top \theta_V) $} & \multirow{4}{8em}{$\tr(\V_\Y \eta_U^\top \theta_U) + \tr(\W_\Y \eta_Y^\top \theta_Y)$, $\V_\Y ,\W_\Y \in \bbS_+(r)$}\\
	 & & & \\
	& & & \\
	& & & \\
	\hline
\end{tabular}
\caption{Basic Properties for Quotient Manifolds $\cM_{r}^{q_1}$, $\cM_{r}^{q_2}$ and $\cM_r^{q_3}$.} \label{tab: basic-prop-quotient-general}
\end{table}

\section{Geometric Connections of Embedded and Quotient Geometries in Fixed-rank PSD Matrix Optimization} \label{sec: connection-PSD}

In this section, we apply the general procedure in Section \ref{sec: general-strategy-for-connection} to connect the landscape properties of optimization \eqref{eq: PSD-manifold-formulation} under the embedded and the quotient geometries. First, under the two quotient geometries introduced in Section \ref{sec: quotient-PSD}, the optimization problem \eqref{eq: PSD-manifold-formulation} can be reformulated as follows,
\begin{subequations}\label{eq: PSD-opt-problem-quotient}
	\begin{align}
		\text{on } \widebar{\cM}_{r+}^{q_1}: &\quad   \min_{\Y \in \bbR^{p \times r}_*}\bar{h}_{r+}(\Y):= f(\Y \Y^\top), \label{eq: PSD-opt-problem-quotient-sub1}\\
		\text{on } \widebar{\cM}_{r+}^{q_2}: &\quad \min_{\U \in \st(r,p), \B \in \bbS_+(r)}\bar{h}_{r+}(\U, \B):= f(\U \B \U^\top) \label{eq: PSD-opt-problem-quotient-sub2}.
	\end{align}
\end{subequations}
 Since $\bar{h}_{r+}(\Y)$ and $\bar{h}_{r+}(\U, \B)$ are invariant along the fibers of $\widebar{\cM}_{r+}^{q_1}$ and $\widebar{\cM}_{r+}^{q_2}$, they induce functions $h_{r+}([\Y])$ and $h_{r+}([\U,\B])$ on quotient manifolds $\cM_{r+}^{q_1}$ and $\cM_{r+}^{q_2}$, respectively. Next, we provide the expressions for Riemannian gradients and Hessians of \eqref{eq: PSD-manifold-formulation} under both geometries (Step 1), construct the bijective maps between $T_\X \cM_{r+}^e$ and $\cH_\Y \widebar{\cM}_{r+}^{q_1}$, $\cH_{(\U,\B)} \widebar{\cM}_{r+}^{q_2}$ (Step 2), and give their spectrum bounds (Step 3).

\begin{proposition}[Riemannian Gradients and Hessians of \eqref{eq: PSD-manifold-formulation}] \label{prop: gradient-hessian-exp-PSD}
	The Riemannian gradients and Hessians of \eqref{eq: PSD-manifold-formulation} under the embedded and the quotient geometries introduced in Section \ref{sec: embedded-quotient-fixed-rank-matrix} are: 
	\begin{itemize}[leftmargin=*]
		\item On $\cM^e_{r+}$: Suppose $\X \in \cM^e_{r+}$, $\U$ spans the top $r$ eigenspace of $\X$, $\xi_\X = [\U \quad \U_\perp] \begin{bmatrix}
			\S & \D^\top\\
			\D & \0
		\end{bmatrix} [\U \quad \U_\perp]^\top \in T_{\X}\cM^e_{r+}$. Then  
		\begin{equation} \label{eq: embedded-gd-hessian-psd}
	\begin{split}
		\grad f(\X) &= P_{\U} \nabla f(\X)P_{\U} + P_{\U_\perp} \nabla f(\X)P_{\U} + P_{\U} \nabla f(\X)P_{\U_\perp},\\
		\Hess f(\X)[\xi_\X, \xi_\X] &= \nabla^2 f(\X)[\xi_\X, \xi_\X] + 2\langle \nabla f(\X), \U_\perp \D \bSigma^{-1} \D^\top \U_\perp^\top \rangle,
	\end{split} 
	\end{equation} where $\bSigma = \U^\top \X \U$. 
		\item On $\cM_{r+}^{q_1}$: Suppose $\Y \in \bbR^{p \times r}_*$ and $\theta_{\Y} \in \cH_{\Y} \widebar{\cM}_{r+}^{q_1}$. Then
		\begin{equation}  \label{eq: quotient-gradient-Hessian-PSD1}
		\begin{split}
			\overline{\grad\, h_{r+}([\Y])} &= 2\nabla f(\Y \Y^\top) \Y \W_\Y^{-1}, \\
			\overline{\Hess \, h_{r+}([\Y])}[\theta_\Y, \theta_\Y] 
		&= \nabla^2 f(\Y \Y^\top)[\Y\theta_\Y^\top + \theta_\Y \Y^\top , \Y\theta_\Y^\top + \theta_\Y \Y^\top ] + 2\langle \nabla f(\Y \Y^\top ), \theta_\Y \theta_\Y^\top \rangle\\
			& +2\langle \nabla f(\Y \Y^\top) \Y  \rmD \W_\Y^{-1} [\theta_\Y], \theta_\Y\W_\Y \rangle + \langle \rmD \W_\Y\left[\overline{ \grad \, h_{r+}([\Y])}\right], \theta_\Y^\top \theta_\Y  \rangle/2.
		\end{split}
		\end{equation}
		\item On $\cM_{r+}^{q_2}$: Suppose $\U \in \st(r,p)$, $\B \in \bbS_+(r)$ and $\theta_{(\U,\B)} = [\theta_U^\top \quad \theta_B^\top]^\top \in \cH_{(\U,\B)} \widebar{\cM}_{r+}^{q_2}$. Then
		\begin{equation}\label{eq: quotient-gradient-Hessian-PSD2}
			\begin{split}
				&\overline{\grad\, h_{r+}([\U,\B])} = \begin{bmatrix}
				\overline{\grad_\U\, h_{r+}([\U,\B])}\\
				\overline{\grad_\B\, h_{r+}([\U,\B])}
			\end{bmatrix} = \begin{bmatrix}
				2 P_{\U_\perp} \nabla f(\U \B \U^\top) \U \B \V_\B^{-1} \\
				\W_\B^{-1} \U^\top \nabla f(\U \B \U^\top) \U \W_\B^{-1}
			\end{bmatrix},\\
			&\overline{\Hess \, h_{r+}([\U,\B])}[\theta_{(\U,\B)}, \theta_{(\U,\B)}]\\
			 &=\nabla^2 f(\U \B \U^\top)[\U \B \theta_U^\top + \U \theta_B \U^\top + \theta_U \B \U^\top, \U \B \theta_U^\top + \U \theta_B \U^\top + \theta_U \B \U^\top]  + 2 \langle  \nabla f(\U \B \U^\top), \theta_U \B \theta_U^\top \rangle\\
			& \,+2\langle \nabla f(\U \B \U^\top) \U, 2 \theta_U \theta_B + \U \rmD \W_\B^{-1}[ \theta_B] \W_\B \theta_B + \theta_U \V_\B \rmD \V_\B^{-1}[\theta_B] \B - \theta_U \U^\top \theta_U \B - \U \theta_U^\top \theta_U \B \rangle \\
			& \, + \tr( \rmD \V_\B [ \overline{ \grad_\B \, h_{r+}([\U, \B])}  ] \theta_U^\top \theta_U   )/2 + \tr( \sym(  \W_\B \theta_B \rmD \W_\B [  \overline{ \grad_\B \, h_{r+}([\U, \B])}  ]  ) \theta_B ).
			\end{split}
		\end{equation}
	\end{itemize}
\end{proposition}

\begin{remark}({\bf Quadratic Form of Riemannian Hessians}) \label{rem: closed-form-riemannian-hessian-quotient} In Proposition \ref{prop: gradient-hessian-exp-PSD}, we only give the quadratic expressions of the Hessians as we use them exclusively throughout the paper. It is easy to obtain general bilinear expressions by noting that $\Hess f(\X)[\xi_\X,\theta_\X] =(\Hess f(\X)[\xi_\X+\theta_\X,\xi_\X+\theta_\X] - \Hess f(\X)[\xi_\X-\theta_\X,\xi_\X-\theta_\X]  )/4 $.

We also note the Riemannian Hessian expressions under the fixed-rank quotient geometries have been explicitly or implicitly developed in \cite{journee2010low,meyer2011geometric,mishra2014fixed,mishra2011low} for some specific problems. Most of these works only provide the linear form of the Riemannian Hessians, which often do not admit a closed-form expression due to the horizontal projection involved in \eqref{eq: quotient-hessian-linear-from}. Here we provide explicit formulas for the quadratic form of the Riemannian Hessians and the closed-form expressions are critical in establishing the landscape connections of embedded and quotient geometries in fixed-rank matrix optimization.
\end{remark}

\begin{proposition}[Bijection Between $\cH_\Y \widebar{\cM}_{r+}^{q_1}$ and $T_\X \cM_{r+}^e$] \label{prop: psd-bijection1} Suppose $\Y \in \bbR^{p \times r}_*$, $\X = \Y \Y^\top$, the eigenspace of $\X$ is $\U$, and $\P = \U^\top \Y$. For any $\theta_\Y \in \cH_\Y \widebar{\cM}_{r+}^{q_1}$ and $\xi_\X =  [\U \quad \U_\perp] \begin{bmatrix}
			\S & \D^\top\\
			\D & \0
		\end{bmatrix} [\U \quad \U_\perp]^\top \in T_{\X}\cM^e_{r+}$, define
\begin{equation} \label{def: xi-theta-correspondence-PSD1}
	\begin{split}
		\xi^{\theta_\Y}_{\X}&:= [\U \quad \U_\perp] \begin{bmatrix}
		\P \theta_\Y^\top \U + \U^\top \theta_\Y \P^\top & \P \theta_\Y^\top \U_\perp \\
		\U_\perp^\top \theta_\Y \P^\top & \0
	\end{bmatrix}[\U \quad \U_\perp]^\top  \in T_{\X}\cM^e_{r+}, \\
	\theta_{\Y}^{\xi_\X} &:= (\U \S' + \U_\perp \D )\P^{-\top} \in \cH_\Y \widebar{\cM}_{r+}^{q_1}, 
	\end{split}
\end{equation} where $\S'$ in $\theta_{\Y}^{\xi_\X}$ is uniquely determined by the linear equation system . Then we can find a linear bijective mapping $\cL_\Y^{r+}$ between $\cH_\Y \widebar{\cM}_{r+}^{q_1}$ and $T_\X \cM_{r+}^e$,
\begin{equation*}
	\cL_\Y^{r+}: \theta_\Y \in  \cH_\Y \widebar{\cM}_{r+}^{q_1} \longrightarrow \xi_{\X}^{\theta_\Y} \in T_\X \cM_{r+}^e \quad \text{and} \quad (\cL_\Y^{r+})^{-1}: \xi_\X \in T_\X \cM_{r+}^e \to \theta_\Y^{\xi_\X} \in \cH_\Y \widebar{\cM}_{r+}^{q_1},
\end{equation*} such that $\cL_\Y^{r+}(\theta_\Y) = \Y \theta_\Y^\top + \theta_\Y \Y^\top$ holds for any $\theta_\Y \in \cH_\Y \widebar{\cM}_{r+}^{q_1}$.

	 Finally, $\cL_{\Y}^{r+}$ satisfies the following spectrum bound:
		\begin{equation} \label{ineq: bijection-spectrum-psd1}
		 2\sigma_r(\P \W_\Y^{-1} \P^\top) \bar{g}^{r+}_\Y(\theta_\Y, \theta_\Y) \leq \|\cL^{r+}_\Y(\theta_\Y)\|_\F^2 \leq 4\sigma_1(\P \W_\Y^{-1} \P^\top) \bar{g}^{r+}_\Y(\theta_\Y, \theta_\Y), \quad \forall \theta_\Y \in \cH_\Y \widebar{\cM}_{r+}^{q_1}.
	\end{equation}	
\end{proposition}
{\noindent \bf Proof of Proposition \ref{prop: psd-bijection1}.} The proof is divided into two steps: in Step 1, we show $\xi_\X^{\theta_\Y}$ and $\theta_{\Y}^{\xi_\X}$ are well defined for any $\theta_\Y$ and $\xi_\X$; in Step 2, we show $\cL_{\Y}^{r+}$ is a bijection and prove its spectrum bounds.

{\bf Step 1.} First, it is clear for any $\theta_\Y \in \cH_\Y \widebar{\cM}_{r+}^{q_1}$, $\xi_\X^{\theta_\Y}$ is well defined. To show $\theta_{\Y}^{\xi_\X}$ is well defined for any $\xi_\X \in T_{\X}\cM^e_{r+}$, we need to show the linear system $\widebar{\S} + \widebar{\S}^{\top} = \S$, $\widebar{\S} \P^{-\top} \W_\Y \P^{-1} = \P^{-\top} \W_\Y \P^{-1} \widebar{\S}^{\top} $ has a unique solution with respect to $\widebar{\S}$. Simple calculations assert that
\begin{equation}
\label{eq: barS subseq Sp}
\begin{split}
	&\{\widebar{\S}: \widebar{\S} + \widebar{\S}^{\top} = \S, \widebar{\S} \P^{-\top} \W_\Y \P^{-1} = \P^{-\top} \W_\Y \P^{-1} \widebar{\S}^{\top}  \} \\
    \subseteq& \{\widetilde{\S}: \P^{-\top} \W_\Y \P^{-1} \widetilde{\S} + \widetilde{\S} \P^{-\top} \W_\Y \P^{-1}  = \P^{-\top} \W_\Y \P^{-1} \S \}.
\end{split}
\end{equation} 
Observing facts that \begin{equation}
\label{eq:SpSylvester}
     \P^{-\top} \W_\Y \P^{-1} \widetilde{\S} + \widetilde{\S} \P^{-\top} \W_\Y \P^{-1}  = \P^{-\top} \W_\Y \P^{-1} \S
\end{equation} 
is a Sylvester equation with respect to $\widetilde{\S}$, and $\P^{-\top} \W_\Y \P^{-1}$ and $-\P^{-\top} \W_\Y \P^{-1}$ have disjoint spectrums, we know from \cite[Theorem VII.2.1]{bhatia2013matrix} that \eqref{eq:SpSylvester} has a unique solution, denoted $\S'$. Since $\S$ is symmetric and $\P^{-\top} \W_\Y \P^{-1}$ is a PSD matrix, we have 
\begin{equation} \label{eq: two-symmetric-equation}
\begin{split}
	\P^{-\top} \W_\Y \P^{-1}  \S' +  \S' \P^{-\top} \W_\Y \P^{-1} - \P^{-\top} \W_\Y \P^{-1} \S &= \0 \\
 \S^{' \top} \P^{-\top} \W_\Y \P^{-1}+\P^{-\top} \W_\Y \P^{-1}  \S^{'\top} - \S\P^{-\top} \W_\Y \P^{-1} &= \0. 
\end{split}
\end{equation} By summing two equations in \eqref{eq: two-symmetric-equation}, we get $\P^{-\top} \W_\Y \P^{-1} ( \S' +  \S^{'\top} - \S) + ( \S' +  \S^{'\top} - \S) \P^{-\top} \W_\Y \P^{-1} = \0 $, i.e., $ \S' +  \S^{'\top} - \S$ is a solution to the new Sylvester equation $\P^{-\top} \W_\Y \P^{-1} \S_1 + \S_1 \P^{-\top} \W_\Y \P^{-1} = \0$ with respect to $\S_1$. Now, we know again by  \cite[Theorem VII.2.1]{bhatia2013matrix} that $\0$ is the unique solution to the system $\P^{-\top} \W_\Y \P^{-1} \S_1 + \S_1 \P^{-\top} \W_\Y \P^{-1} = \0$. Thus, $ \S' +  \S^{'\top} = \S$, and from \eqref{eq: two-symmetric-equation}, it further holds $ \S' \P^{-\top} \W_\Y \P^{-1} =\P^{-\top} \W_\Y \P^{-1}  \S^{'\top} $. This, together with \eqref{eq: barS subseq Sp} and the uniqueness of $ \S'$, asserts that $\S'$ is the unique solution to the linear  system $\widebar{\S} \P^{-\top} \W_\Y \P^{-1} = \P^{-\top} \W_\Y \P^{-1}\widebar{\S}^\top $, $\widebar{\S} + \widebar{\S}^\top = \S$. This finishes the proof of this part.

{\bf Step 2.} Note that both $\cH_\Y \widebar{\cM}_{r+}^{q_1}$ and $T_\X \cM_{r+}^e$ are of dimension $(pr - (r^2-r)/2)$. Let $\cL_\Y^{r+'} : \xi_\X \in T_\X\cM^e_{r+} \longrightarrow \theta_\Y^{\xi_\X} \in \cH_\Y \widebar{\cM}_{r+}^{q_1}$. Moreover, for any $\xi_\X = [\U \quad \U_\perp] \begin{bmatrix}
			\S & \D^\top\\
			\D & \0
		\end{bmatrix} [\U \quad \U_\perp]^\top\in T_\X\cM^e_{r+}$, we have
\begin{equation} \label{eq: bijection-PSD1}
	\cL^{r+}_\Y ( \cL_\Y^{r+'}(\xi_\X) ) = \cL^{r+}_\Y (\theta_\Y^{\xi_\X}) =  [\U \quad \U_\perp] \begin{bmatrix}
		\P  \theta_\Y^{\xi_\X \top} \U + \U^\top  \theta_\Y^{\xi_\X} \P^\top & \P  \theta_\Y^{\xi_\X\top} \U_\perp \\
		\U_\perp^\top \theta_\Y^{\xi_\X} \P^\top & \0
	\end{bmatrix}[\U \quad \U_\perp]^\top = \xi_\X.
\end{equation}
Therefore, $\cL^{r+}_\Y$ is a bijection and $\cL_\Y^{r+'} = (\cL^{r+}_\Y)^{-1}$. At the same time, it is easy to check $ \cL_\Y^{r+}$ satisfies $\cL_\Y^{r+}(\theta_\Y)=\Y \theta_\Y^\top + \theta_\Y \Y^\top$ by observing $\Y = \U \P$.

Next, we provide the spectrum bounds for $\cL^{r+}_\Y$. For any $\theta_\Y = (\U\S + \U_\perp \D )\P^{-\top} \in \cH_\Y \widebar{\cM}_{r+}^{q_1}$, we have
	\begin{equation} \label{ineq: upper-bound-theta-norm}
	\begin{split}
		\bar{g}^{r+}_\Y(\theta_\Y,\theta_\Y) = \tr(\W_\Y \theta_\Y^\top \theta_\Y) = \| \theta_\Y\W_\Y^{1/2}\|_\F^2& \leq \sigma^2_1(\P^{-\top}\W_\Y^{1/2})  ( \|\S\|_\F^2 + \|\D\|_\F^2 ) \\
		& = \sigma_1(\P^{-\top}\W_\Y \P^{-1}) ( \|\S\|_\F^2 + \|\D\|_\F^2 ) \\
		& = ( \|\S\|_\F^2 + \|\D\|_\F^2 )/\sigma_r(\P \W_\Y^{-1} \P^\top),
	\end{split}
	\end{equation} and 
	\begin{equation} \label{ineq: quadratic-S'-bound-PSD1}
		\begin{split}
			\langle \S^{\top}, \S \rangle &\overset{(a)}= \langle \P \W_\Y^{-1} \P^\top \S \P^{-\top} \W_\Y \P^{-1}, \S \rangle \\
			&= \langle ( \P \W_\Y^{-1} \P^\top)^{1/2} \S ( \P \W_\Y^{-1} \P^\top)^{-1/2}, ( \P \W_\Y^{-1} \P^\top)^{1/2} \S ( \P \W_\Y^{-1} \P^\top)^{-1/2} \rangle \geq 0,
		\end{split}
	\end{equation} here (a) is because $ \S \P^{-\top} \W_\Y \P^{-1} = \P^{-\top} \W_\Y \P^{-1} \S^\top $ by the construction of $\cH_\Y \widebar{\cM}_{r+}^{q_1}$. Thus
\begin{equation} \label{ineq: L_Y^r+-bound}
\begin{split}
	\|\cL^{r+}_\Y(\theta_\Y)\|_\F^2 = \|\xi^{\theta_\Y}_{\X}\|_\F^2 & \overset{ \eqref{def: xi-theta-correspondence-PSD1}  }= \|\P \theta_\Y^\top \U + \U^\top \theta_\Y \P^\top \|_\F^2 + 2 \|\U_\perp^\top \theta_\Y \P^\top \|_\F^2\\
	& = \|\S^\top + \S\|_\F^2 + 2\|\D\|_\F^2\\
	& \overset{\eqref{ineq: quadratic-S'-bound-PSD1} } \geq 2 (\|\S\|_\F^2 + \|\D\|_\F^2) \overset{\eqref{ineq: upper-bound-theta-norm} } \geq 2\sigma_r(\P \W_\Y^{-1} \P^\top)\bar{g}^{r+}_\Y(\theta_\Y,\theta_\Y),
\end{split}
\end{equation} and
\begin{equation*}
		\begin{split}
		\|\cL^{r+}_\Y(\theta_\Y)\|_\F^2 = \|\xi^{\theta_\Y}_{\X}\|_\F^2 & \overset{ \eqref{def: xi-theta-correspondence-PSD1}  }= \|\P \theta_\Y^\top \U + \U^\top \theta_\Y \P^\top \|_\F^2 + 2 \|\U_\perp^\top \theta_\Y \P^\top \|_\F^2 \\
		& \leq 4 \|\U^\top \theta_\Y \P^\top\|_\F^2 + 2 \|\U_\perp^\top \theta_\Y \P^\top \|_\F^2\\
		& = 2 \|\U^\top \theta_\Y \P^\top\|_\F^2 + 2 \|\theta_\Y \P^\top\|_\F^2 \\
		& \leq 4  \|\theta_\Y \P^{\top}\|_\F^2 = 4  \|\theta_\Y \W_\Y^{1/2} \W_\Y^{-1/2} \P^{\top}\|_\F^2 \\
		&\overset{\eqref{ineq: upper-bound-theta-norm} }\leq 4\sigma_1(\P \W_\Y^{-1} \P^\top ) \bar{g}^{r+}_\Y(\theta_\Y,\theta_\Y).
	\end{split}
\end{equation*} This finishes the proof of this proposition. \quad $\blacksquare$

\begin{proposition}[Bijection Between $\cH_{(\U,\B)} \widebar{\cM}_{r+}^{q_2}$ and $T_\X \cM_{r+}^e$] \label{prop: psd-bijection2} Suppose $\U \in \st(r,p)$, $\B \in \bbS_+(r)$ and $\X = \U \B \U^\top$. For any $\theta_{(\U,\B)} = [\theta_U^\top \quad \theta_B^\top]^\top \in \cH_{(\U,\B)} \widebar{\cM}_{r+}^{q_2}$ and $\xi_\X =  [\U \quad \U_\perp] \begin{bmatrix}
			\S & \D^\top\\
			\D & \0
		\end{bmatrix} [\U \quad \U_\perp]^\top \in T_{\X}\cM^e_{r+}$, define
\begin{equation} \label{def: xi-theta-correspondence-PSD2}
	\begin{split}
		\xi^{\theta_{(\U,\B)}}_{\X}&:= [\U \quad \U_\perp] \begin{bmatrix}
		\theta_B & \B \theta_U^\top \U_\perp \\
		\U_\perp^\top \theta_U \B & \0
	\end{bmatrix}[\U \quad \U_\perp]^\top  \in T_{\X}\cM^e_{r+}, \\
	\theta_{(\U,\B)}^{\xi_\X} &:= [(\U_\perp \D \B^{-1})^\top \quad \S]^\top \in \cH_{(\U,\B)} \widebar{\cM}_{r+}^{q_2}.
	\end{split}
\end{equation} Then we can find a linear bijective mapping $\cL_{\U,\B}^{r+}$ between $\cH_{(\U,\B)} \widebar{\cM}_{r+}^{q_2}$ and $T_\X \cM_{r+}^e$,
\begin{equation*}
\begin{split}
	&\cL_{\U,\B}^{r+}: \theta_{(\U,\B)} \in  \cH_{(\U,\B)} \widebar{\cM}_{r+}^{q_2} \longrightarrow \xi_{\X}^{\theta_{(\U,\B)}} \in T_\X \cM_{r+}^e, \\
	&(\cL_{\U,\B}^{r+})^{-1}: \xi_\X \in T_\X \cM_{r+}^e \to \theta_{(\U,\B)}^{\xi_\X} \in \cH_{(\U,\B)} \widebar{\cM}_{r+}^{q_2}
\end{split}
\end{equation*}such that $\cL_{\U,\B}^{r+}(\theta_{(\U,\B)}) = \U \B \theta_U^\top + \U \theta_B \U^\top + \theta_U \B \U^\top$ holds for any $\theta_{(\U,\B)} \in \cH_{(\U,\B)} \widebar{\cM}_{r+}^{q_2}$.	

Finally, $\cL_{\U,\B}^{r+}$ satisfies the following spectrum bound: $\forall \theta_{(\U,\B)} \in \cH_{(\U,\B)} \widebar{\cM}_{r+}^{q_2}$,
		\begin{equation} \label{ineq: bijection-spectrum-psd2}
		\begin{split}
		(\sigma^2_r(\W_\B^{-1}) \wedge 2 \sigma_r^2(\V_\B^{-1/2} \B) ) \bar{g}^{r+}_{(\U,\B)}(\theta_{(\U,\B)}, \theta_{(\U,\B)}) &\leq \|\cL^{r+}_{\U,\B}(\theta_{(\U,\B)})\|_\F^2\\
		 &\quad \leq  (\sigma^2_1(\W^{-1}_\B) \vee 2 \sigma^2_1( \V_\B^{-1/2} \B ) )\bar{g}^{r+}_{(\U,\B)}(\theta_{(\U,\B)}, \theta_{(\U,\B)}).	
		\end{split}
	\end{equation}	
\end{proposition}
{\noindent \bf Proof of Proposition \ref{prop: psd-bijection2}.} First, it is easy to see $\xi^{\theta_{(\U,\B)}}_{\X}$ and $\theta_{(\U,\B)}^{\xi_\X}$ are well defined given any $\theta_{(\U,\B)} \in \cH_{(\U,\B)} \widebar{\cM}_{r+}^{q_2}$ and $\xi_\X \in T_{\X}\cM^e_{r+}$. 

Next, we show $\cL_{\U,\B}^{r+}$ is a bijection. Notice $\cH_{(\U,\B)} \widebar{\cM}_{r+}^{q_2}$ is of dimension $(pr - (r^2-r)/2)$, which is the same with $T_\X \cM_{r+}^e$. Suppose $\cL_{\U,\B}^{r+'} : \xi_\X \in T_\X\cM^e_{r+} \longrightarrow \theta_{(\U,\B)}^{\xi_\X} \in \cH_{(\U,\B)} \widebar{\cM}_{r+}^{q_2}$. For any $\xi_\X = [\U \quad \U_\perp] \begin{bmatrix}
			\S & \D^\top\\
			\D & \0
		\end{bmatrix} [\U \quad \U_\perp]^\top \in T_\X\cM^e_{r+}$, we have
\begin{equation} \label{eq: bijection-PSD2}
	\cL^{r+}_{\U,\B} ( \cL_{\U,\B}^{r+'}(\xi_\X) ) = \cL^{r+}_{\U,\B} (\theta_{(\U,\B)}^{\xi_\X}) =  [\U \quad \U_\perp] \begin{bmatrix}
		\theta_B^{\xi_\X} & \B \theta^{\xi_\X \top}_U \U_\perp \\
		\U_\perp^\top \theta_U^{\xi_\X} \B & \0
	\end{bmatrix}[\U \quad \U_\perp]^\top = \xi_\X.
\end{equation} Since $\cL^{r+}_{\U,\B}$ and $\cL_{\U,\B}^{r+'}$ are linear maps, \eqref{eq: bijection-PSD2} implies $\cL^{r+}_{\U,\B}$ is a bijection and $\cL_{\U,\B}^{r+'} = (\cL^{r+}_{\U,\B})^{-1}$. At the same time, it is easy to check $\cL_{\U,\B}^{r+}(\theta_{(\U,\B)}) = \U \B \theta_U^\top + \U \theta_B \U^\top + \theta_U \B \U^\top$ holds for any $\theta_{(\U,\B)} \in \cH_{(\U,\B)} \widebar{\cM}_{r+}^{q_2}$.

Next, we provide the spectrum bounds for $\cL^{r+}_{\U,\B}$. For any $\theta_{(\U,\B)} = [(\U_\perp \D)^\top \quad \theta_B]^\top \in \cH_{(\U,\B)} \widebar{\cM}_{r+}^{q_2}$, we have $\bar{g}^{r+}_{(\U,\B)}(\theta_{(\U,\B)},\theta_{(\U,\B)}) = \|\D\V_\B^{1/2}\|_\F^2 + \tr(\W_\B \theta_B \W_\B \theta_B) = \|\D\V_\B^{1/2}\|_\F^2 + \|\W_\B^{1/2} \theta_B \W_\B^{1/2}\|_\F^2$. Thus, we have
\begin{equation*}
	\begin{split}
		\|\cL^{r+}_{\U,\B}(\theta_{(\U,\B)})\|_\F^2 \overset{ \eqref{def: xi-theta-correspondence-PSD2} } = \|\theta_\B\|_\F^2 + 2 \|\B \theta_\U^\top \U_\perp\|_\F^2 &= \|\theta_\B\|_\F^2 + 2\|\B \D^\top\|_\F^2\\
		& = \|\W_\B^{-1/2}\W_\B^{1/2} \theta_\B \W_\B^{1/2} \W_\B^{-1/2}\|_\F^2 + 2\|\D\V_\B^{1/2} \V_\B^{-1/2}  \B \|_\F^2 \\
		 &\leq (\sigma^2_1(\W^{-1}_\B) \vee 2 \sigma^2_1( \V_\B^{-1/2} \B ) ) \left( \|\W_\B^{1/2} \theta_B \W_\B^{1/2}\|_\F^2  +  \|\D\V_\B^{1/2}\|_\F^2\right)\\
		& = (\sigma^2_1(\W^{-1}_\B) \vee 2 \sigma^2_1( \V_\B^{-1/2} \B ) ) \bar{g}^{r+}_{(\U,\B)}(\theta_{(\U,\B)},\theta_{(\U,\B)}),
	\end{split}
\end{equation*} and
\begin{equation*}
	\begin{split}
		\|\cL^{r+}_{\U,\B}(\theta_{(\U,\B)})\|_\F^2  \overset{ \eqref{def: xi-theta-correspondence-PSD2} }= \|\theta_\B\|_\F^2 + 2\|\B \D^\top\|_\F^2 &\geq (\sigma^2_r(\W_\B^{-1}) \wedge 2 \sigma_r^2(\V_\B^{-1/2} \B) ) \left( \|\W_\B^{1/2} \theta_\B \W_\B^{1/2}\|_\F^2  +  \|\D\V_\B^{1/2}\|_\F^2\right)\\
		& = (\sigma^2_r(\W_\B^{-1}) \wedge 2 \sigma_r^2(\V_\B^{-1/2} \B) ) \bar{g}^{r+}_{(\U,\B)}(\theta_{(\U,\B)},\theta_{(\U,\B)}). \quad \quad \quad  \blacksquare
	\end{split}
\end{equation*}

Next, we present our first main result on the geometric landscape connections of Riemannian optimization \eqref{eq: PSD-manifold-formulation} under the embedded and the quotient geometries. 

\begin{theorem}[Geometric Landscape Connections of \eqref{eq: PSD-manifold-formulation} on $\cM_{r+}^e$ and $\cM_{r+}^{q_1}$] \label{th: embedded-quotient-connection-PSD1}
	 Suppose the conditions in Proposition \ref{prop: psd-bijection1} hold and the $\W_\Y$ in $\bar{g}_\Y^{r+}$ satisfies $\W_\Y = \O \W_{\Y\O} \O^\top$ for any $\O \in \bbO_r$. Then
	\begin{equation} \label{eq: gradient-connect-PSD1}
	\begin{split}
		\grad f(\X) &= \left( \overline{\grad\, h_{r+}([\Y])} \W_\Y \Y^\dagger + \left(\overline{\grad\, h_{r+}([\Y])} \W_\Y \Y^\dagger\right)^\top(\I_p - \Y \Y^\dagger) \right)/2,\\
		 \overline{\grad\, h_{r+}([\Y])} &= 2 \grad f(\X) \Y \W_\Y^{-1}.
	\end{split}
	\end{equation}
Furthermore, if $[\Y]$ is a Riemannian FOSP of \eqref{eq: PSD-opt-problem-quotient-sub1}, we have: 
\begin{equation} \label{eq: Hessian-connection-PSD1}
	\begin{split}
		\overline{\Hess \, h_{r+}([\Y])}[\theta_\Y, \theta_\Y]  &= \Hess f(\X)[\cL_\Y^{r+}(\theta_\Y),\cL_\Y^{r+}(\theta_\Y)], \quad \forall \theta_\Y \in \cH_\Y \widebar{\cM}_{r+}^{q_1}.
	\end{split}
	\end{equation}
	Finally, $\overline{\Hess \, h_{r+}([\Y])}$ and $\Hess f(\X)$ have $(pr-(r^2-r)/2)$ eigenvalues; for $i = 1,\ldots, pr-(r^2-r)/2$, 
	 \begin{equation*}
	    \begin{split}
	        & \lambda_i(\overline{\Hess \, h_{r+}([\Y])}) \text{ is sandwiched between } \\
	        & 2\sigma_r(\P \W_\Y^{-1} \P^\top) \lambda_i(\Hess f(\X))  \text{ and } 4\sigma_1(\P \W_\Y^{-1} \P^\top) \lambda_i(\Hess f(\X)).
	    \end{split}
	 \end{equation*}
\end{theorem}
{\noindent \bf Proof of Theorem \ref{th: embedded-quotient-connection-PSD1}.} First, notice $\Y$ lies in the column space spanned by $\U$ and $\Y \Y^\dagger = P_\U$. So \eqref{eq: gradient-connect-PSD1} is by direct calculation from the gradient expressions in Proposition \ref{prop: gradient-hessian-exp-PSD}. 

Next, we prove \eqref{eq: Hessian-connection-PSD1}. Since $[\Y]$ is a Riemannian FOSP of \eqref{eq: PSD-opt-problem-quotient-sub1}, we have
\begin{equation} \label{eq: FOSP-condition-PSD1}
	\overline{ \grad \, h_{r+}([\Y])} = \0, \quad \nabla f(\Y \Y^\top) \Y = \0, 
\end{equation} 
and have $\grad f(\X) = \0$ by \eqref{eq: gradient-connect-PSD1}. So $ \nabla f(\X) = P_{\U_\perp} \nabla f(\X) P_{\U_\perp}$. Recall $\P = \U^\top \Y$, $\X = \Y \Y^\top$ and let $\bSigma = \U^\top \X \U$. Given any $\theta_\Y \in \cH_{\Y} \widebar{\cM}_{r+}^{q_1}$, we have
\begin{equation} \label{eq: Hessian-con-gradient-1}
	\langle \nabla f(\X), P_{\U_\perp} \theta_\Y \P^{\top} \bSigma^{-1} \P \theta_\Y^\top  P_{\U_\perp}    \rangle = \langle \nabla f(\X), \theta_\Y \theta_\Y^\top \rangle,
\end{equation} where the equality is because $\P$ is nonsingular, $\P \P^\top = \bSigma$ and $ \nabla f(\X) = P_{\U_\perp} \nabla f(\X) P_{\U_\perp} $.

Then by Proposition \ref{prop: gradient-hessian-exp-PSD}:
\begin{equation*} 
\begin{split}
	& \quad \overline{\Hess \, h_{r+}([\Y])}[\theta_\Y, \theta_\Y] \\
	&=   \nabla^2 f(\Y \Y^\top)[\Y\theta_\Y^\top + \theta_\Y \Y^\top , \Y\theta_\Y^\top + \theta_\Y \Y^\top ] + 2\langle \nabla f(\Y \Y^\top ), \theta_\Y \theta_\Y^\top \rangle\\
	&\quad  + 2\langle \nabla f(\Y \Y^\top) \Y  \rmD \W_\Y^{-1} [\theta_\Y], \theta_\Y\W_\Y \rangle + \langle \rmD \W_\Y\left[\overline{ \grad \, h_{r+}([\Y])}\right], \theta_\Y^\top \theta_\Y  \rangle/2\\
	& \overset{ \eqref{eq: FOSP-condition-PSD1} }=  \nabla^2 f(\X)[\Y\theta_\Y^\top + \theta_\Y \Y^\top , \Y\theta_\Y^\top + \theta_\Y \Y^\top ] + 2\langle \nabla f(\Y \Y^\top ), \theta_\Y \theta_\Y^\top \rangle \\
	& \overset{ \text{Proposition } \ref{prop: psd-bijection1}, \eqref{eq: Hessian-con-gradient-1} }= \nabla^2 f(\X)[\cL_\Y^{r+}(\theta_\Y),\cL_\Y^{r+}(\theta_\Y)] + 2\langle \nabla f(\X), P_{\U_\perp} \theta_\Y \P^{\top} \bSigma^{-1} \P \theta_\Y^\top  P_{\U_\perp}    \rangle\\
	& = \Hess f(\X)[\cL_\Y^{r+}(\theta_\Y),\cL_\Y^{r+}(\theta_\Y)],
\end{split}
\end{equation*} 
where the last equality follows from the expression of  $\Hess f(\X)$ in \eqref{eq: embedded-gd-hessian-psd} and the definition of $\cL_\Y^{r+}$.

Then, by \eqref{ineq: bijection-spectrum-psd1}, \eqref{eq: Hessian-connection-PSD1} and Theorem \ref{th: hessian-sandwich}, we have $\overline{\Hess \, h_{r+}([\Y])}$ and $\Hess f(\X)$ have $(pr-(r^2-r)/2)$ eigenvalues and $\widebar{\lambda}_i(\Hess \, h_{r+}([\Y]))$ is sandwiched between $2\sigma_r(\P \W_\Y^{-1} \P^\top) \lambda_i(\Hess f(\X)) $ and $4\sigma_1(\P \W_\Y^{-1} \P^\top) \lambda_i(\Hess f(\X))$ for $i = 1,\ldots,pr-(r^2-r)/2$. \quad $\blacksquare$

\begin{theorem}[Geometric Landscape Connections of \eqref{eq: PSD-manifold-formulation} on $\cM_{r+}^e$ and $\cM_{r+}^{q_2}$] \label{th: embedded-quotient-connection-PSD2}
	 Suppose $\U \in \st(r,p)$, $\B \in \bbS_+(r)$, $\X = \U \B \U^\top$ and $\V_\B$, $\W_\B$ in $\bar{g}^{r+}_{(\U,\B)}$ satisfies $\V_\B = \O \V_{\O^\top \B \O} \O^\top$ and $\W_\B = \O \W_{\O^\top \B \O} \O^\top$ for any $\O \in \bbO_r$. Then
	\begin{equation} \label{eq: gradient-connect-PSD2}
	\begin{split}
		\grad f(\X) &= ( \overline{\grad_\U\, h_{r+}([\U,\B])} \V_\B \B^{-1} \U^\top )/2 + ( \overline{\grad_\U\, h_{r+}([\U,\B])} \V_\B  \B^{-1} \U^\top )^\top/2 \\
		& \quad + \U \W_\B \overline{\grad_\B\, h_{r+}([\U,\B])} \W_\B \U^\top,\\
		 \overline{\grad\, h_{r+}([\U,\B])} &= \begin{bmatrix}
		 	2 P_{\U_\perp} \grad f(\X) \U \B \V^{-1}_\B \\
		 	\W^{-1}_\B \U^\top \grad f(\X) \U \W^{-1}_\B
		 \end{bmatrix}.
	\end{split}
	\end{equation}
Furthermore, if $[\U,\B]$ is a Riemannian FOSP of \eqref{eq: PSD-opt-problem-quotient-sub2}, we have: 
\begin{equation} \label{eq: Hessian-connection-PSD2}
	\begin{split}
		\overline{\Hess \, h_{r+}([\U,\B])}[\theta_{(\U,\B)}, \theta_{(\U,\B)}]  &= \Hess f(\X)[\cL_{\U,\B}^{r+}(\theta_{(\U,\B)}),\cL_{\U,\B}^{r+}(\theta_{(\U,\B)})], \quad \forall \theta_{(\U,\B)} \in \cH_{(\U,\B)} \widebar{\cM}_{r+}^{q_2}.
	\end{split}
	\end{equation}
	Finally, $\overline{\Hess \, h_{r+}([\U,\B])}$ has $(pr-(r^2-r)/2)$ eigenvalues and for $i = 1,\ldots, pr-(r^2-r)/2$,
	\begin{equation*}
		\begin{split}
			& \lambda_i(\overline{\Hess \, h_{r+}([\U,\B])}) \text{ is sandwiched between }\\
			&(\sigma^2_r(\W_\B^{-1}) \wedge 2 \sigma_r^2(\V_\B^{-1/2} \B) ) \lambda_i(\Hess f(\X)) \text{ and }  (\sigma^2_1(\W^{-1}_\B) \vee 2 \sigma^2_1( \V_\B^{-1/2} \B ) ) \lambda_i(\Hess f(\X)).
		\end{split}
	\end{equation*}
\end{theorem}
{\noindent \bf Proof of Theorem \ref{th: embedded-quotient-connection-PSD2}.} The proof of this theorem is similar to the proof of Theorem \ref{th: embedded-quotient-connection-PSD1} and is postponed to Appendix \ref{sec: addition-proof-psd}. \quad $\blacksquare$

In Table \ref{tab: illustraion-gap-coefficient-PSD}, we compute the explicit gap coefficients in the sandwich inequalities in Theorems \ref{th: embedded-quotient-connection-PSD1} and \ref{th: embedded-quotient-connection-PSD2} under several Riemannian metrics specified by $\W_\Y, \V_\B$ and $\W_\B$ in $\bar{g}^{r+}_\Y$ and $\bar{g}^{r+}_{(\U,\B)}$.
\begin{table}[ht]
	\centering
	\begin{tabular}{c | c | c | c}
		\hline
		 & \multirow{2}{10em}{Choices of $\W_\Y, \V_\B$ and $\W_\B$ in $\bar{g}^{r+}$} & \multirow{2}{8em}{Gap Coefficient Lower Bound} & \multirow{2}{8em}{Gap Coefficient Upper Bound} \\
		 & & & \\
		 \hline
		 \multirow{2}{7em}{$\cM_{r+}^e$ v.s. $\cM^{q_1}_{r+}$} & $\W_\Y = \I_r$  & $2\sigma_r(\X)$  & $4 \sigma_1(\X)$ \\
		 \cline{2-4} 
		 & $\W_\Y = 2 \Y^\top \Y$ & $1$ & $2$ \\
		 \hline
		 \multirow{2}{7em}{$\cM_{r+}^e$ v.s. $\cM^{q_2}_{r+}$} & $\V_\B = \I_r$, $\W_\B = \B^{-1}$  & $\sigma^2_r(\X)$  & $2 \sigma^2_1(\X)$ \\
		  \cline{2-4} 
		 & $\V_\B = 2\B^2, \W_\B = \I_r$  & $1$  & $1$\\
		 \hline
	\end{tabular}
	\caption{Gap coefficients in the sandwich inequalities in Theorems \ref{th: embedded-quotient-connection-PSD1} and \ref{th: embedded-quotient-connection-PSD2} under different metrics.
	} \label{tab: illustraion-gap-coefficient-PSD}
\end{table}

We note in Theorems \ref{th: embedded-quotient-connection-PSD1} and \ref{th: embedded-quotient-connection-PSD2}, the connection of Riemannian gradients under two geometries is established at every point on the manifold while the connection between Riemaniann Hessians is established at Riemannian FOSPs. Theorems \ref{th: embedded-quotient-connection-PSD1} and \ref{th: embedded-quotient-connection-PSD2} also show the following equivalence of Riemannian FOSPs, SOSPs, and strict saddles of optimization \eqref{eq: PSD-manifold-formulation} under the embedded and the quotient geometries.
\begin{corollary} ({\bf Equivalence on Riemannian FOSPs, SOSPs and strict saddles of \eqref{eq: PSD-manifold-formulation} Under Embedded and Quotient Geometries})\label{coro: landscape connection PSD} Suppose $\W_\Y = \O \W_{\Y\O} \O^\top$, $\V_\B = \O \V_{\O^\top \B \O} \O^\top$ and $\W_\B = \O \W_{\O^\top \B \O} \O^\top$ hold for any $\O \in \bbO_r$. Then we have
\begin{itemize}
\item[(a)] given $\Y \in \bbR^{p \times r}_*, \U \in \st(r,p)$ and $\B \in \bbS_{+}(r)$, if $[\Y]$ ($[\U,\B]$) is a Riemannian FOSP or SOSP or strict saddle of \eqref{eq: PSD-opt-problem-quotient-sub1} (\eqref{eq: PSD-opt-problem-quotient-sub2}), then $\X = \Y \Y^\top$ ($\X = \U \B \U^\top$) is a Riemannian FOSP or SOSP or strict saddle of \eqref{eq: PSD-manifold-formulation} under the embedded geometry; 
\item[(b)] if $\X$ is a Riemannian FOSP or SOSP or strict saddle of \eqref{eq: PSD-manifold-formulation} under the embedded geometry, then there is a unique $[\Y]$ ($[\U,\B]$) such that $\Y \Y^\top = \X$ ($\U \B \U^\top = \X$) and it is a Riemannian FOSP or SOSP or strict saddle of \eqref{eq: PSD-opt-problem-quotient-sub1} (\eqref{eq: PSD-opt-problem-quotient-sub2}).
\end{itemize}
\end{corollary}

Corollary \ref{coro: landscape connection PSD} is proved based on Theorems \ref{th: embedded-quotient-connection-PSD1} and \ref{th: embedded-quotient-connection-PSD2} in Appendix \ref{sec: addition-proof-psd}. We also note that on the same manifold, the equivalence of Riemannian SOSPs under two geometries with different metrics can be deduced more directly from the fact that the quadratic form of the Riemannian Hessian does not depend on the metric at Riemannian FOSPs from \cite[Eq. (5.29)]{absil2009optimization}. To derive the connection result between two different manifolds considered here, one needs to further explore the relationship between manifolds.

\begin{remark}({\bf Spectrum Relation of Riemannian Hessians}) The sandwich inequalities in Theorems \ref{th: embedded-quotient-connection-PSD1} and \ref{th: embedded-quotient-connection-PSD2} provide a finer connection on the spectrum of the Riemannian Hessians under the embedded and the quotient geometries at Riemannian FOSPs. These results are useful in transferring a few common geometric landscape properties from one geometry formulation to another. One such example is the so-called strict saddle property \cite{ge2015escaping,lee2019first}, which states that the function has a strict negative curvature at all stationary points but local minima. With this strict saddle property, various Riemannian gradient descent and trust-region methods are guaranteed to escape all strict saddles and converge to an SOSP \cite{boumal2019global,criscitiello2019efficiently,lee2019first,sun2018geometric,sun2019escaping}.
\end{remark}

\begin{remark}({\bf Effects of Riemannian Metric and Quotient Structure on Landscape Connection}) \label{rem: effiect-of-metric-on-landscape} Corollary \ref{coro: landscape connection PSD} shows the choices of the quotient structure and metric do not affect the landscape connection of Riemannian FOSPs and SOSPs under two geometries. However, they do affect the gap of the sandwich inequalities in Theorems \ref{th: embedded-quotient-connection-PSD1} and \ref{th: embedded-quotient-connection-PSD2}. For example, in the scenario of Theorem \ref{th: embedded-quotient-connection-PSD1} with $\W_\Y = \I_r$, the coefficients in lower and upper bounds of the sandwich inequality are $2\sigma_r(\X)$ and $4\sigma_1(\X)$, and their ratio relies on the condition number of $\X$. \cite{zhang2021preconditioned,zheng2022riemannian,zhuo2021computational} showed recently that such dependence on the condition number of $\X$ yields slow convergence of factorization-based algorithms, especially when $\X$ is ill-conditioned or the rank of the input matrix $\Y$ is overspecified. Moreover, the ratio of upper and lower bounds of the sandwich inequality becomes $2$ if we choose $\W_\Y = 2\Y^\top \Y$, meaning that the algorithms under this metric are more robust to the ill-condition of $\X$ \cite{tong2020accelerating,zhang2021preconditioned}. More surprisingly, the spectrum of two Riemannian Hessians can be exactly the same as the upper and lower bounds of sandwich inequality match in some special setting, e.g., $\cM_{r+}^e$ versus $\cM_{r+}^{q_2}$ with $\V_\B = 2\B^2$ and $\W_\B = \I_r$ (see Table \ref{tab: illustraion-gap-coefficient-PSD}). To make this happen, the chosen metrics for positive definite matrices and the Stiefel manifold in $\cM_{r+}^{q_2}$ are non-standard. See Remark \ref{rem: W_Y-choice-in-metric} for a detailed discussion. As we will discuss in Remark \ref{rem: algorithmic-connection}, when the gap coefficients are universal constants, there is an algorithmic connection between adopting embedded and quotient geometries in Riemannian fixed-rank matrix optimization, and this connection is promoted to an equivalence when the spectrum of two Riemannian Hessians coincide.
	
\end{remark}

\begin{remark}({\bf Implications on Landscape Connections of Different Geometries for Riemannian Optimization}) \label{rem: implication-on-connection-diff-approaches}
	Generally speaking, embedded and quotient geometries are the most common two choices in Riemannian optimization. Compared to quotient geometry, embedded geometry allows computing and interpreting many geometric objects straightforwardly. Theorems \ref{th: embedded-quotient-connection-PSD1} and \ref{th: embedded-quotient-connection-PSD2} establish a strong geometric landscape connection between two geometries in fixed-rank PSD matrix optimization and this provides an example under which two different geometries are indeed connected in treating the same constraint in Riemannian optimization. Finally, we note although we focus on the geometric connections of \eqref{eq: PSD-manifold-formulation} under the embedded and the quotient geometries, it is also relatively easy to obtain the geometric connections under different quotient geometries based on our results.
\end{remark}

\section{Geometric Connections of Embedded and Quotient Geometries in Fixed-rank General Matrix Optimization} \label{sec: connection-general}
In this section, we present the geometric landscape connections of optimization problem \eqref{eq: general prob} under the embedded and the quotient geometries. First, the problem \eqref{eq: general prob} can be reformulated under each of the three quotient geometries in Section \ref{sec: quotient-general} as follows,
\begin{subequations}\label{eq: general-opt-problem-quotient}
	\begin{align}
		\text{on } \widebar{\cM}_{r}^{q_1}:& \quad  \min_{\L \in \bbR^{p_1 \times r}_*,\R \in \bbR^{p_2 \times r}_*}\bar{h}_{r}(\L,\R):=  f(\L \R^\top), \label{eq: general-opt-problem-quotient-sub1}\\
		\text{on } \widebar{\cM}_{r}^{q_2}:& \quad \min_{\U \in \st(r,p_1), \B \in \bbS_+(r), \V \in \st(r,p_2)}\bar{h}_{r}(\U, \B, \V):=  f(\U \B \V^\top), \label{eq: general-opt-problem-quotient-sub2}\\
		\text{on } \widebar{\cM}_{r}^{q_3}:& \quad \min_{\U \in \st(r,p_1), \Y \in \bbR^{p_2 \times r}_*}\bar{h}_{r}(\U, \Y):= f(\U \Y^\top).\label{eq: general-opt-problem-quotient-sub3}
	\end{align}
\end{subequations}
Since $\bar{h}_{r}(\L,\R)$, $\bar{h}_{r}(\U, \B, \V)$ and $\bar{h}_{r}(\U, \Y)$ are invariant along the fibers of $\widebar{\cM}_{r}^{q_1}$, $\widebar{\cM}_{r}^{q_2}$ and $\widebar{\cM}_{r}^{q_3}$, they induce functions $h_{r}([\L,\R])$, $h_{r}([\U,\B,\V])$ and $h_r([\U,\Y])$ on quotient manifolds $\cM_{r}^{q_1}$, $\cM_{r}^{q_2}$ and $\cM_{r}^{q_3}$, respectively. In the following Proposition \ref{prop: gradient-hessian-exp-general}, we provide Riemannian gradients and Hessians of \eqref{eq: general prob} under embedded and quotient geometries. 

\begin{proposition}[Riemannian Gradients and Hessians of \eqref{eq: general prob}] \label{prop: gradient-hessian-exp-general}
	The Riemannian gradients and Hessians of \eqref{eq: general prob} under the embedded and the quotient geometries introduced in Section \ref{sec: embedded-quotient-fixed-rank-matrix} are: 
	\begin{itemize}[leftmargin=*]
		\item On $\cM^e_{r}$: Suppose $\X \in \cM_{r}^e$, $\U,\V$ span the top $r$ left and right singular subspaces of $\X$, respectively and $\xi_\X = [\U \quad \U_\perp] \begin{bmatrix}
			\S & \D_2^\top\\
			\D_1 & \0
		\end{bmatrix} [\V \quad \V_\perp]^\top \in T_{\X}\cM^e_{r}$. Then 
		\begin{equation}\label{eq: embedded-gd-hessian-general}
	\begin{split}
		\grad f(\X) &= P_{\U} \nabla f(\X)P_{\V} + P_{\U_\perp} \nabla f(\X)P_{\V} + P_{\U} \nabla f(\X)P_{\V_\perp},\\
		\Hess f(\X)[\xi_\X, \xi_\X] &= \nabla^2 f(\X)[\xi_\X, \xi_\X] + 2\langle \nabla f(\X), \U_\perp \D_1 \bSigma^{-1} \D_2^\top \V_\perp^\top \rangle,
	\end{split}
	\end{equation} where $\bSigma = \U^\top \X \V$.
		\item On $\cM_{r}^{q_1}$: Suppose $\L \in \bbR^{p_1 \times r}_*$, $\R \in \bbR^{p_2 \times r}_*$ and $\theta_{(\L,\R)} = [\theta_L^\top \quad \theta_R^\top]^\top \in \cH_{(\L,\R)} \widebar{\cM}_{r}^{q_1}$. Then
		\begin{equation}  \label{eq: quotient-gradient-Hessian-general1}
		\begin{split}
				&\overline{\grad\, h_{r}([\L ,\R])} = \begin{bmatrix}
				\overline{\grad_{\L}\, h_{r}([\L ,\R])}\\
				\overline{\grad_{\R}\, h_{r}([\L ,\R])}
			\end{bmatrix} =  \begin{bmatrix}
				\nabla f(\L \R^\top) \R \W_{\L ,\R}^{-1} \\
				(\nabla f(\L \R^\top))^\top \L \V_{\L ,\R}^{-1} 
			\end{bmatrix}, \\
						&\overline{\Hess \, h_{r}([\L,\R])}[\theta_{(\L,\R)}, \theta_{(\L,\R)}] \\
		 &=  \nabla^2 f(\L \R^\top)[\L \theta_R^\top + \theta_L \R^\top,\L \theta_R^\top + \theta_L \R^\top] + 2 \langle \nabla f(\L \R^\top), \theta_L \theta_R^\top \rangle \\
		 & \quad + \langle \nabla f(\L \R^\top) \R \rmD \W_{\L,\R}^{-1} [\theta_{(\L,\R)}] , \theta_L \W_{\L,\R} \rangle + \langle  \left(\nabla f(\L \R^\top) \right)^\top \L \rmD \V_{\L,\R}^{-1} [\theta_{(\L,\R)}], \theta_R \V_{\L,\R} \rangle \\
		 & \quad +  \langle \rmD \W_{\L,\R}[ \overline{ \grad \, h_{r}([\L,\R])} ], \theta_L^\top \theta_L \rangle /2 +  \langle \rmD \V_{\L,\R} [ \overline{ \grad \, h_{r}([\L,\R])} ], \theta_R^\top \theta_R \rangle /2.
		\end{split}
		\end{equation}
		\item On $\cM_{r}^{q_2}$: Suppose $\U \in \st(r,p_1)$, $\B \in \bbS_+(r)$, $\V \in \st(r,p_2)$ and $\theta_{(\U,\B,\V)} = [\theta_U^\top \quad \theta_B^\top \quad \theta_V^\top]^\top \in \cH_{(\U,\B,\V)} \widebar{\cM}_r^{q_2}$. Then
		\begin{equation}\label{eq: quotient-gradient-Hessian-general2}
			\begin{split}
				&\overline{\grad\, h_{r}([\U,\B,\V])} = \begin{bmatrix}
				\overline{\grad_\U\, h_{r}([\U,\B,\V])}\\
				\overline{\grad_\B\, h_{r}([\U,\B,\V])} \\
				\overline{\grad_\V\, h_{r}([\U,\B,\V])}
			\end{bmatrix} \\
			& \quad \quad \quad \quad \quad \quad \quad \quad = \begin{bmatrix}
				P_{\U_\perp} \nabla f(\U \B \V^\top) \V \B + \U\left(\skew(\bDelta )\B + \B  \skew(\bDelta ) \right)/2 \\
				\B \sym(\bDelta) \B \\
				P_{\V_\perp} (\nabla f(\U \B \V^\top))^\top \U \B - \V \left( \skew(\bDelta )\B + \B  \skew(\bDelta ) \right)/2
			\end{bmatrix},\\
			& \overline{\Hess \, h_{r}([\U,\B,\V])}[\theta_{(\U,\B,\V)}, \theta_{(\U,\B,\V)}] \\
		=&\nabla^2 f(\U \B \V^\top)[\theta_U \B \V^\top + \U \theta_B \V^\top + \U \B \theta_V^\top, \theta_U \B \V^\top + \U \theta_B \V^\top + \U \B \theta_V^\top] + 2\langle \nabla f(\U\B\V^\top), \theta_U \B \theta_V^\top \rangle\\
		& +   \left\langle \bDelta,  \sym(\U^\top \theta_U \U^\top\theta_U) \B + \B \sym(\V^\top \theta_V \U^\top \theta_U) -2\theta_U^\top \theta_U \B \right\rangle/2\\
		& + \left\langle \bDelta,  \B\sym(\V^\top \theta_V \V^\top\theta_V) + \sym(\U^\top \theta_U \V^\top \theta_V) \B -2\B\theta_V^\top \theta_V +  2\theta_B \B^{-1} \theta_B \right\rangle/2\\
		& + \langle \bDelta', 2\theta_B -  \U^\top \theta_U \B - \theta_U^\top \U \B/2 -\V^\top \theta_V \B/2 \rangle + \langle \bDelta'', 2\theta_B -  \B\theta_V^\top \V - \B\V^\top \theta_V/2 -\B \theta_U^\top \U /2 \rangle,
			\end{split}
		\end{equation} where $\bDelta =  \U^\top\nabla f(\U \B \V^\top) \V$, $\bDelta' =  \theta_U^\top\nabla f(\U \B \V^\top) \V$ and $\bDelta'' =  \U^\top\nabla f(\U \B \V^\top) \theta_V$.
		\item On $\cM_r^{q_3}$: Suppose $\U \in \st(r,p_1)$, $\Y \in \bbR^{p_2 \times r}_*$ and  $\theta_{(\U,\Y)} = [\theta_U^\top \quad \theta_Y^\top]^\top \in \cH_{(\U,\Y)} \widebar{\cM}_{r}^{q_3}$. Then
		\begin{equation}\label{eq: quotient-gradient-Hessian-general3}
			\begin{split}
					&\overline{\grad\, h_{r}([\U,\Y])} = \begin{bmatrix}
				\overline{\grad_\U\, h_{r}([\U,\Y])}\\
				\overline{\grad_\Y\, h_{r}([\U,\Y])}
			\end{bmatrix} = \begin{bmatrix}
				P_{\U_\perp} \nabla f(\U \Y^\top)\Y \V_\Y^{-1} \\
				(\nabla f(\U \Y^\top))^\top \U \W_\Y^{-1}
			\end{bmatrix},\\
			&\overline{\Hess \, h_{r}([\U,\Y])}[\theta_{(\U,\Y)}, \theta_{(\U,\Y)}] \\
		=& \nabla^2 f(\U \Y^\top)[\U \theta_Y^\top + \theta_U \Y^\top, \U \theta_Y^\top + \theta_U \Y^\top] + 2 \langle  \nabla f(\U \Y^\top), \theta_U \theta_Y^\top \rangle - \langle  \U^\top \nabla f(\U \Y^\top) \Y, \theta_U^\top \theta_U \rangle \\
		& + \langle (\nabla f(\U \Y^\top))^\top \U \rmD \W_\Y^{-1} [\theta_Y], \theta_Y \W_\Y \rangle + \langle \nabla f(\U \Y^\top) \Y \rmD \V_\Y^{-1} [\theta_Y], \theta_U \V_\Y \rangle\\
		& + \langle \rmD \W_\Y\left[\overline{ \grad_\Y \, h_{r}([\U,\Y])}\right], \theta_Y^\top \theta_Y  \rangle /2 + \langle \rmD \V_\Y\left[\overline{ \grad_\Y \, h_{r}([\U,\Y])}\right], \theta_U^\top \theta_U  \rangle /2.
			\end{split}
		\end{equation}
	\end{itemize}
\end{proposition}

 Next, we construct bijective maps between $T_\X \cM_{r}^e$ and $\cH_{(\L,\R)} \widebar{\cM}_{r}^{q_1}$, $\cH_{(\U,\B,\V)} \widebar{\cM}_{r}^{q_2}$, $\cH_{(\U,\Y)} \widebar{\cM}_{r}^{q_3}$, and give their spectrum bounds.
\begin{proposition}[Bijection Between $\cH_{(\L,\R)} \widebar{\cM}_{r}^{q_1}$ and $T_\X \cM_{r}^e$] \label{prop: general-bijection1} Suppose $\L \in \bbR^{p_1 \times r}_*$, $\R \in \bbR^{p_2 \times r}_*$, $\X = \L \R^\top$ with its top $r$ left and right singular subspaces spanned by $\U$ and $\V$, respectively and $\P_1 = \U^\top \L$, $\P_2 = \V^\top \R$. For any $\theta_{(\L,\R)} = [\theta_L^\top \quad \theta_R^\top]^\top \in \cH_{(\L,\R)} \widebar{\cM}_{r}^{q_1}$ and $\xi_\X =  [\U \quad \U_\perp] \begin{bmatrix}
			\S & \D_2^\top\\
			\D_1 & \0
		\end{bmatrix} [\V \quad \V_\perp]^\top \in T_{\X}\cM^e_{r}$, define
\begin{equation} \label{def: xi-theta-correspondence-general1}
	\begin{split}
		\xi^{\theta_{(\L,\R)}}_{\X}&:= [\U \quad \U_\perp] \begin{bmatrix}
		\P_1 \theta_R^\top \V + \U^\top \theta_L \P_2^\top & \P_1 \theta_R^\top \V_\perp \\
		\U_\perp^\top \theta_L \P_2^\top & \0
	\end{bmatrix}[\V \quad \V_\perp]^\top  \in T_{\X}\cM^e_{r}, \\
	\theta_{{(\L,\R)}}^{\xi_\X} &:= \begin{bmatrix}
		(\U \S' \P_2 \W_{\L,\R}^{-1} \P_2^\top + \U_\perp \D_1 ) \P_2^{-\top}\\
		(\V \S^{'\top} \P_1 \V_{\L,\R}^{-1} \P_1^\top + \V_\perp \D_2 )\P_1^{-\top}
	\end{bmatrix} \in \cH_{(\L,\R)} \widebar{\cM}_{r}^{q_1},
	\end{split}
\end{equation} where $\S'$ in $\theta_{(\L,\R)}^{\xi_\X}$ is uniquely determined by the Sylvester equation $\P_1 \V_{\L,\R}^{-1} \P_1^\top \S' + \S' \P_2 \W_{\L,\R}^{-1} \P_2^\top   = \S$. Then we can find a linear bijective mapping $\cL_{\L,\R}^{r}$ between $\cH_{(\L,\R)} \widebar{\cM}_{r}^{q_1}$ and $T_\X \cM_{r}^e$,
\begin{equation*}
	\cL_{\L,\R}^{r}: \theta_{(\L,\R)} \in  \cH_{(\L,\R)} \widebar{\cM}_{r}^{q_1} \longrightarrow \xi_{\X}^{\theta_{(\L,\R)}} \in T_\X \cM_{r}^e \quad \text{and} \quad (\cL_{\L,\R}^{r})^{-1}: \xi_\X \in T_\X \cM_{r}^e \to \theta_{(\L,\R)}^{\xi_\X} \in \cH_{(\L,\R)} \widebar{\cM}_{r}^{q_1},
\end{equation*} such that $\cL_{\L,\R}^{r}(\theta_{(\L,\R)}) = \L \theta_R^\top + \theta_L \R^\top$ holds for any $\theta_{(\L,\R)} \in \cH_{(\L,\R)} \widebar{\cM}_{r}^{q_1}$.

	 Finally, we have the following spectrum bounds for $\cL_{\L,\R}^{r}$:
		\begin{equation} \label{ineq: bijection-spectrum-general1}
		 \gamma_{\L,\R} \cdot \bar{g}^{r}_{(\L,\R)}(\theta_{(\L,\R)}, \theta_{(\L,\R)}) \leq \|\cL^{r}_{\L,\R}(\theta_{(\L,\R)})\|_\F^2 \leq 2 \Gamma_{\L,\R} \cdot \bar{g}^{r}_{(\L,\R)}(\theta_{(\L,\R)}, \theta_{(\L,\R)}), \quad \forall \theta_{(\L,\R)} \in \cH_{(\L,\R)} \widebar{\cM}_{r}^{q_1},
	\end{equation} where $\gamma_{\L,\R} :=  \sigma_r(\P_2 \W_{\L,\R}^{-1} \P_2^\top) \wedge \sigma_r(\P_1 \V_{\L,\R}^{-1} \P_1^\top) $ and $\Gamma_{\L,\R} :=  \sigma_1(\P_2 \W_{\L,\R}^{-1} \P_2^\top) \vee \sigma_1(\P_1 \V_{\L,\R}^{-1} \P_1^\top)$.
\end{proposition}
{\noindent \bf Proof of Proposition \ref{prop: general-bijection1}.} First, the uniqueness of $\S'$ in $\theta_{(\L,\R)}^{\xi_\X}$ is guaranteed by the fact $\P_1 \V_{\L,\R}^{-1} \P_1^\top$ and $-\P_2 \W_{\L,\R}^{-1} \P_2^\top$ have disjoint spectra and \cite[Theorem VII.2.1]{bhatia2013matrix}. Thus, $\xi^{\theta_{(\L,\R)}}_{\X}$ and $\theta_{{(\L,\R)}}^{\xi_\X}$ are well defined given any $\theta_{(\L,\R)}$ and $\xi_\X$.

Next, we show $\cL^r_{\L,\R}$ is a bijection. Notice both $\cH_{(\L,\R)} \widebar{\cM}_{r}^{q_1}$ and $T_\X \cM_{r}^e$ are of dimension $(p_1 + p_2-r)r$. Suppose $\cL_{\L,\R}^{r'} : \xi_\X \in T_\X\cM^e_{r} \longrightarrow \theta_{(\L,\R)}^{\xi_\X} \in \cH_{(\L,\R)} \widebar{\cM}_{r}^{q_1}$. For any $\xi_\X = [\U \quad \U_\perp] \begin{bmatrix}
			\S & \D_2^\top\\
			\D_1 & \0
		\end{bmatrix} [\V \quad \V_\perp]^\top\in T_\X\cM^e_{r}$, we have
\begin{equation} \label{eq: bijection-general1}
	\cL^{r}_{\L,\R} ( \cL_\Y^{r'}(\xi_\X) ) = \cL^{r}_{\L,\R} (\theta_{(\L,\R)}^{\xi_\X}) =  [\U \quad \U_\perp] \begin{bmatrix}
		\P_1  \theta_R^{\xi_\X \top} \V + \U^\top  \theta_L^{\xi_\X} \P_2^\top & \P_1  \theta_R^{\xi_\X\top} \V_\perp \\
		\U_\perp^\top \theta_L^{\xi_\X} \P_2^\top & \0
	\end{bmatrix}[\V \quad \V_\perp]^\top = \xi_\X.
\end{equation}
Since $\cL^{r}_{\L,\R}$ and $\cL_{\L,\R}^{r'}$ are linear maps, \eqref{eq: bijection-general1} implies $\cL^{r}_{\L,\R}$ is a bijection and $\cL_{\L,\R}^{r'} = (\cL^{r}_{\L,\R})^{-1}$. At the same time, it is easy to check $ \cL_{\L,\R}^{r}$ satisfies $\cL_{\L,\R}^{r}(\theta_{(\L,\R)})=\L \theta_R^\top + \theta_L \R^\top $ by observing $\L = \U \P_1, \R = \V \P_2$.

Next, we provide the spectrum bounds for $\cL^{r}_{\L,\R}$. For any $\theta_{(\L,\R)} = [\theta_L^\top \quad \theta_R^\top]^\top \in \cH_{(\L,\R)} \widebar{\cM}_{r}^{q_1}$, where $\theta_L = (\U \S \P_2\W_{\L, \R}^{-1} \P_2^\top + \U_\perp \D_1 ) \P_2^{-\top}$, $\theta_R =(\V \S^\top \P_1 \V_{\L,\R}^{-1} \P_1^\top + \V_\perp \D_2  ) \P_1^{-\top} $ we have
\begin{equation} \label{ineq: upper-bound-theta-norm-general}
	\begin{split}
		&\quad \bar{g}_{(\L,\R)}(\theta_{(\L,\R)}, \theta_{(\L,\R)}) \\
		&= \tr(\W_{\L,\R} \theta_L^\top \theta_L) + \tr(\V_{\L,\R} \theta_R^\top \theta_R) = \|\theta_L \W_{\L,\R}^{1/2}\|_\F^2 + \|\theta_R \V_{\L,\R}^{1/2}\|_\F^2 \\
		& \leq (\|\S \P_2\W_{\L, \R}^{-1} \P_2^\top\|_\F^2 + \|\D_1\|_\F^2 ) \sigma_1^2(\P_2^{-\top} \W_{\L,\R}^{1/2})  + (\|\S^\top \P_1 \V_{\L,\R}^{-1} \P_1^\top\|_\F^2 + \|\D_2\|_\F^2  ) \sigma_1^2(\P_1^{-\top} \V_{\L,\R}^{1/2}) \\
		& \leq \frac{1}{ \sigma_r(\P_2 \W_{\L,\R}^{-1} \P_2^\top ) \wedge \sigma_r(\P_1 \V_{\L,\R}^{-1} \P_1^\top ) } (\|\S \P_2\W_{\L, \R}^{-1} \P_2^\top\|_\F^2 + \|\D_1\|_\F^2 + \|\S^\top \P_1 \V_{\L,\R}^{-1} \P_1^\top\|_\F^2 + \|\D_2\|_\F^2 ),
	\end{split}
\end{equation} and
\begin{equation}\label{ineq: quadratic-S'-bound-general1}
	\begin{split}
		&\langle  \P_1 \V_{\L,\R}^{-1} \P_1^\top \S, \S \P_2\W_{\L, \R}^{-1} \P_2^\top \rangle\\
		 =& \langle  (\P_1 \V_{\L,\R}^{-1} \P_1^\top)^{1/2} \S(\P_2\W_{\L, \R}^{-1} \P_2^\top)^{1/2}, (\P_1 \V_{\L,\R}^{-1} \P_1^\top)^{1/2}\S (\P_2\W_{\L, \R}^{-1} \P_2^\top)^{1/2} \rangle \geq 0.
	\end{split}
\end{equation} Thus
\begin{equation} \label{ineq: L_(L,R)^r-bound}
\begin{split}
	\|\cL^{r}_{\L,\R}(\theta_{(\L,\R)})\|_\F^2 = \|\xi^{\theta_{(\L,\R)}}_{\X}\|_\F^2 & \overset{ \eqref{def: xi-theta-correspondence-general1}  }= \|\P_1 \theta_R^\top \V + \U^\top \theta_L \P_2^\top \|_\F^2 +  \|\P_1 \theta_R^\top \V_\perp \|_\F^2 + \|\U_\perp^\top \theta_L \P_2^\top\|_\F^2 \\
	& = \| \P_1 \V_{\L,\R}^{-1} \P_1^\top \S +  \S \P_2\W_{\L, \R}^{-1} \P_2^\top\|_\F^2 + \|\D_1\|_\F^2 + \|\D_2\|_\F^2\\
	& \overset{\eqref{ineq: quadratic-S'-bound-general1} } \geq \|\S \P_2\W_{\L, \R}^{-1} \P_2^\top\|_\F^2 + \|\S^\top \P_1 \V_{\L,\R}^{-1} \P_1^\top\|_\F^2+ \|\D_1\|_\F^2 + \|\D_2\|_\F^2 \\
	& \overset{\eqref{ineq: upper-bound-theta-norm-general} } \geq (\sigma_r(\P_2 \W_{\L,\R}^{-1} \P_2^\top ) \wedge \sigma_r(\P_1 \V_{\L,\R}^{-1} \P_1^\top )) \bar{g}^{r}_{(\L,\R)}(\theta_{(\L,\R)},\theta_{(\L,\R)}),
\end{split}
\end{equation} and
\begin{equation*}
\begin{split}
\|\cL^{r}_{\L,\R}(\theta_{(\L,\R)})\|_\F^2 = \|\xi^{\theta_{(\L,\R)}}_{\X}\|_\F^2 & \overset{ \eqref{def: xi-theta-correspondence-general1}  }= \|\P_1 \theta_R^\top \V + \U^\top \theta_L \P_2^\top \|_\F^2 +  \|\P_1 \theta_R^\top \V_\perp \|_\F^2 + \|\U_\perp^\top \theta_L \P_2^\top\|_\F^2\\
& \leq 2( \|\P_1 \theta_R^\top \V  \|_\F^2 + \| \U^\top \theta_L \P_2^\top \|_\F^2) +  \|\P_1 \theta_R^\top \V_\perp \|_\F^2 + \|\U_\perp^\top \theta_L \P_2^\top\|_\F^2\\
& \leq 2(  \|\P_1 \theta_R^\top  \|_\F^2 + \| \theta_L \P_2^\top \|_\F^2 ) \\
& \overset{\eqref{ineq: upper-bound-theta-norm-general}} \leq  2( \sigma^2_1(\W_{\L,\R}^{-1/2} \P_2^\top ) \vee \sigma^2_1(\V_{\L,\R}^{-1/2} \P_1^\top ) )\bar{g}^{r}_{(\L,\R)}(\theta_{(\L,\R)},\theta_{(\L,\R)})\\
& =  2( \sigma_1(\P_2\W_{\L,\R}^{-1} \P_2^\top ) \vee \sigma_1(\P_1\V_{\L,\R}^{-1} \P_1^\top ) )\bar{g}^{r}_{(\L,\R)}(\theta_{(\L,\R)},\theta_{(\L,\R)})
\end{split}
\end{equation*}
This finishes the proof of this proposition.
\quad $\blacksquare$

\begin{proposition}[Bijection Between $\cH_{(\U,\B,\V)} \widebar{\cM}_{r}^{q_2}$ and $T_\X \cM_{r}^e$] \label{prop: general-bijection2} Suppose $\U \in \st(r,p_1)$, $\B \in \bbS_+(r)$, $\V \in \st(p_2,r)$ and $\X = \U \B \V^\top$. For any $\theta_{(\U,\B,\V)} = [\theta_U^\top \quad \theta_B^\top \quad \theta_V^\top]^\top \in \cH_{(\U,\B,\V)} \widebar{\cM}_{r}^{q_2}$ and $\xi_\X =  [\U \quad \U_\perp] \begin{bmatrix}
			\S & \D_2^\top\\
			\D_1 & \0
		\end{bmatrix} [\V \quad \V_\perp]^\top \in T_{\X}\cM^e_{r}$, define
\begin{equation} \label{def: xi-theta-correspondence-general2}
	\begin{split}
		\xi^{\theta_{(\U,\B,\V)}}_{\X}&:= [\U \quad \U_\perp] \begin{bmatrix}
		\U^\top \theta_U \B + \theta_B + \B \theta_V^\top \V & \B \theta_V^\top \V_\perp \\
		\U_\perp^\top \theta_U \B & \0
	\end{bmatrix}[\V \quad \V_\perp]^\top  \in T_{\X}\cM^e_{r}, \\
	\theta_{(\U,\B,\V)}^{\xi_\X} &:= [(\U_\perp \D_1 \B^{-1} + \U \bOmega')^\top \quad \S' \quad (\V_\perp \D_2 \B^{-1} - \V \bOmega')^\top]^\top \in \cH_{(\U,\B,\V)} \widebar{\cM}_{r}^{q_2},
	\end{split}
\end{equation} where $\S', \bOmega'$ are uniquely determined by the linear equation system: $\bOmega' \B + \S' - \B \bOmega^{'\top} = \S$, $\bOmega' = - \bOmega^{'\top}$, $\S' = \S^{'\top}$. Then we can find a linear bijective mapping $\cL_{\U,\B,\V}^{r}$ between $\cH_{(\U,\B,\V)} \widebar{\cM}_{r}^{q_2}$ and $T_\X \cM_{r}^e$,
\begin{equation*}
\begin{split}
	&\cL_{\U,\B,\V}^{r}: \theta_{(\U,\B,\V)} \in  \cH_{(\U,\B,\V)} \widebar{\cM}_{r}^{q_2} \longrightarrow \xi_{\X}^{\theta_{(\U,\B,\V)}} \in T_\X \cM_{r}^e, \\
	&(\cL_{\U,\B,\V}^{r})^{-1}: \xi_\X \in T_\X \cM_{r}^e \to \theta_{(\U,\B,\V)}^{\xi_\X} \in \cH_{(\U,\B,\V)} \widebar{\cM}_{r}^{q_2}
\end{split}
\end{equation*}such that $\cL_{\U,\B,\V}^{r}(\theta_{(\U,\B,\V)}) = \theta_U \B \V^\top + \U \theta_B \V^\top + \U \B \theta_V^\top$ holds for any $\theta_{(\U,\B,\V)} \in \cH_{(\U,\B,\V)} \widebar{\cM}_{r}^{q_2}$.	

Finally, we have the following spectrum bounds for $\cL_{\U,\B,\V}^{r}$: for all $\theta_{(\U,\B,\V)} \in \cH_{(\U,\B,\V)} \widebar{\cM}_{r}^{q_2}$,
		\begin{equation} \label{ineq: bijection-spectrum-general2}
		\sigma^2_r(\X)\bar{g}^{r}_{(\U,\B,\V)}(\theta_{(\U,\B,\V)}, \theta_{(\U,\B,\V)}) \leq \|\cL^{r}_{\U,\B,\V}(\theta_{(\U,\B,\V)})\|_\F^2 \leq 2 \sigma^2_1(\X)\bar{g}^{r}_{(\U,\B,\V)}(\theta_{(\U,\B,\V)}, \theta_{(\U,\B,\V)}),
	\end{equation}	
\end{proposition}
{\noindent \bf Proof of Proposition \ref{prop: general-bijection2}.} The proof is divided into two steps: in Step 1, we show $\xi_\X^{\theta_{(\U,\B,\V)}}$ and $\theta_{(\U,\B,\V)}^{\xi_\X}$ are well defined for any $\theta_{(\U,\B,\V)}$ and $\xi_\X$; in Step 2, we show $\cL_{\U\,\B,\V}^{r}$ is a bijection and prove its spectrum bounds.

{\bf Step 1.} First, it is clear for any $\theta_{(\U,\B,\V)} \in \cH_{(\U,\B,\V)} \widebar{\cM}_{r}^{q_2}$, $\xi_\X^{\theta_{(\U,\B,\V)}}$ is well defined. To show $\theta_{(\U,\B,\V)}^{\xi_\X}$ is well defined given any $\xi_\X \in T_{\X}\cM^e_{r}$, we need to show the equation system: $(i) \bOmega_1 \B + \S_1 - \B \bOmega_1^\top = \S$, (ii) $\bOmega_1 = - \bOmega_1^\top$, (iii) $\S_1 = \S_1^{\top}$ with respect to $\S_1, \bOmega_1$ has a unique solution. By (i) we have $\S_1 = \S - \bOmega_1 \B + \B \bOmega_1^\top$. By plugging it into (iii), we have $\B \bOmega_1^\top - \bOmega_1 \B = (\S^\top - \S)/2$. Combining it with (ii) states $\B \bOmega_1 + \bOmega_1 \B = (\S - \S^\top)/2$. So we conclude
\begin{equation*}
	\begin{split}
		&\{ (\S_1, \bOmega_1): \bOmega_1 \B + \S_1 - \B \bOmega_1^\top = \S, \bOmega_1 = - \bOmega_1^\top, \S_1 = \S_1^{\top} \} \\
		\subseteq & \{(\S_1, \bOmega_1): \B \bOmega_1 + \bOmega_1 \B = (\S - \S^\top)/2,  \S_1 = \S - \bOmega_1 \B + \B \bOmega_1^\top \}.
	\end{split}
\end{equation*} 
Note that $\B \bOmega_1 + \bOmega_1 \B = (\S - \S^\top)/2$ is a Sylvester equation with respect to $\bOmega_1$. Since $\B$ and $-\B$ have distinct spectra, the system has a unique solution \cite[Theorem VII.2.1]{bhatia2013matrix} and we denote it by $\bOmega'$. Let $\S' =  \S - \bOmega' \B + \B \bOmega^{'\top}$. If we can show $\bOmega' = - \bOmega^{'\top}$ and $\S' = \S^{'\top}$, then we can conclude $(\S',\bOmega')$ is the unique solution of the linear equation system $\bOmega_1 \B + \S_1 - \B \bOmega_1^\top = \S, \bOmega_1 = - \bOmega_1^\top, \S_1 = \S_1^{\top}$.

Let us first show $\bOmega' = - \bOmega^{'\top}$. We know $\bOmega'$ satisfies $\B \bOmega' + \bOmega' \B = (\S - \S^\top)/2$ and $\bOmega^{'\top}\B + \B\bOmega^{'\top} = (\S^\top - \S)/2$. By summing these two equations we have $\B (\bOmega' + \bOmega^{'\top}) + (\bOmega' + \bOmega^{'\top}) \B = \0$. This is a new Sylvester equation $\B \widebar{\bOmega} + \widebar{\bOmega} \B = \0$ with respect to $\widebar{\bOmega}$ and we know again by  \cite[Theorem VII.2.1]{bhatia2013matrix} that $\0$ is the unique solution to this system. So we have $\bOmega' + \bOmega^{'\top} = \0$, i.e., $\bOmega' = - \bOmega^{'\top}$. Then
\begin{equation*}
\begin{split}
	\S' =  \S - \bOmega' \B + \B \bOmega^{'\top} &= \S + \bOmega^{'\top} \B + \B \bOmega^{'\top} = \S + (\S^\top - \S)/2 \\
	&= (\S^\top + \S)/2 = \S^\top + \B \bOmega^{'} + \bOmega'\B = \S^\top - \B \bOmega^{'\top} + \bOmega'\B = \S^{'\top}.
\end{split}
\end{equation*} So we have shown  $\xi_\X^{\theta_{(\U,\B,\V)}}$ is well defined for any $\theta_{(\U,\B,\V)} \in \cH_{(\U,\B,\V)} \widebar{\cM}_{r}^{q_2}$.

{\bf Step 2.} Notice $\cH_{(\U,\B,\V)} \widebar{\cM}_{r}^{q_2}$ is of dimension $(p_1 + p_2 -r)r$ and it is the same with $\dim(T_\X \cM_{r}^e)$. Suppose $\cL_{\U,\B,\V}^{r'} : \xi_\X \in T_\X\cM^e_{r} \longrightarrow \theta_{(\U,\B,\V)}^{\xi_\X} \in \cH_{(\U,\B,\V)} \widebar{\cM}_{r}^{q_2}$. For any $\xi_\X = [\U \quad \U_\perp] \begin{bmatrix}
			\S & \D_2^\top\\
			\D_1 & \0
		\end{bmatrix} [\V \quad \V_\perp]^\top \in T_\X\cM^e_{r}$, we have
\begin{equation} \label{eq: bijection-general2}
\begin{split}
	\cL^{r}_{\U,\B,\V} ( \cL_{\U,\B,\V}^{r'}(\xi_\X) ) &= \cL^{r}_{\U,\B,\V} (\theta_{(\U,\B,\V)}^{\xi_\X}) \\
	&=  [\U \quad \U_\perp] \begin{bmatrix}
		\U^\top \theta_U^{\xi_\X} \B + \theta_B^{\xi_\X} + \B \theta_V^{\xi_\X \top} \V &  \B \theta^{\xi_\X \top}_V \V_\perp \\
		\U_\perp^\top \theta_U^{\xi_\X} \B & \0
	\end{bmatrix}[\V \quad \V_\perp]^\top = \xi_\X.
\end{split}
\end{equation} Since $\cL^{r}_{\U,\B,\V}$ and $\cL_{\U,\B,\V}^{r'}$ are linear maps, \eqref{eq: bijection-general2} implies $\cL^{r}_{\U,\B,\V}$ is a bijection and $\cL_{\U,\B,\V}^{r'} = (\cL^{r}_{\U,\B,\V})^{-1}$. At the same time, it is easy to check $\cL_{\U,\B,\V}^{r}(\theta_{(\U,\B,\V)}) = \theta_U \B \V^\top + \U \theta_B \V^\top + \U \B \theta_V^\top$ holds for any $\theta_{(\U,\B,\V)} \in \cH_{(\U,\B,\V)} \widebar{\cM}_{r}^{q_2}$.

Next, we provide the spectrum bounds for $\cL_{\U,\B,\V}^r$. For any $\theta_{(\U,\B,\V)} = [\theta_U^\top \quad \theta_B^\top \quad \theta_V^\top]^\top \in \cH_{(\U,\B,\V)} \widebar{\cM}_{r}^{q_2}$ with $\theta_U = \U_\perp \D_1 + \U \bOmega, \theta_V = \V_\perp \D_2 - \V \bOmega$, $\bOmega = - \bOmega^\top$, $\theta_B \in \bbS^{r \times r}$. We have
\begin{equation} \label{eq: upper-bound-theta-norm-general2}
	\begin{split}
		\bar{g}^{r}_{(\U,\B,\V)}(\theta_{(\U,\B,\V)}, \theta_{(\U,\B,\V)}) &= \|\theta_U\|_\F^2 + \|\theta_V\|_\F^2 + \tr(\B^{-1} \theta_B \B^{-1} \theta_B) \\
		& = \|\theta_U\|_\F^2 + \|\theta_V\|_\F^2 + \|\B^{-1/2} \theta_B \B^{-1/2}\|_\F^2\\
		& = \|\D_1\|_\F^2 + \|\D_2\|_\F^2 + 2 \|\bOmega\|_\F^2 + \|\B^{-1/2} \theta_B \B^{-1/2}\|_\F^2,
	\end{split}
\end{equation} and
\begin{equation} \label{ineq: quadratic-S'-bound-general2}
	\begin{split}
		\langle \bOmega \B - \B \bOmega^\top, \theta_B \rangle = \langle \B \bOmega^\top - \bOmega \B, \theta_B^\top \rangle \overset{\theta_B = \theta_B^\top}= - \langle \bOmega \B - \B \bOmega^\top, \theta_B \rangle \Longrightarrow\langle \bOmega \B - \B \bOmega^\top, \theta_B \rangle = 0.
			\end{split}
\end{equation}
Thus,
\begin{equation} \label{ineq: bijection-general2-ineq1}
	\begin{split}
		\|\cL^r_{\U,\B,\V}(\theta_{(\U,\B,\V)})\|_\F^2 = \|\xi_\X^{\theta_{(\U,\B,\V)}}\|_\F^2 &\overset{ \eqref{def: xi-theta-correspondence-general2}} = \|\U^\top \theta_U \B + \theta_B + \B \theta_V^\top \V\|_\F^2 + \|\B \theta_V^\top \V_\perp\|_\F^2 + \|\U_\perp^\top \theta_U \B\|_\F^2  \\
		&= \|\bOmega \B + \theta_B - \B \bOmega^\top\|_\F^2 + \|\B \D_2^\top\|_\F^2 + \|\D_1 \B\|_\F^2\\
		& \overset{ \eqref{ineq: quadratic-S'-bound-general2} } = \|\bOmega \B  - \B \bOmega^\top\|_\F^2 + \|\theta_B\|_\F^2 + \|\B \D_2^\top\|_\F^2 + \|\D_1 \B\|_\F^2 \\
		& \overset{ \bOmega = -\bOmega^\top } = 2 \|\bOmega \B \|_\F^2 + 2\|\B^{1/2} \bOmega \B^{1/2} \|_\F^2 + \|\theta_B\|_\F^2 + \|\B \D_2^\top\|_\F^2 + \|\D_1 \B\|_\F^2 \\
		& \overset{(a)} \geq \sigma^2_r(\X) \bar{g}^{r}_{(\U,\B,\V)}(\theta_{(\U,\B,\V)}, \theta_{(\U,\B,\V)}),
	\end{split}
\end{equation} where in (a) we use the fact $\sigma_r(\B) = \sigma_r(\X)$, and
\begin{equation*}
	\begin{split}
		\|\cL^r_{\U,\B,\V}(\theta_{(\U,\B,\V)})\|_\F^2 = \|\xi_\X^{\theta_{(\U,\B,\V)}}\|_\F^2 
		&\overset{ \eqref{ineq: bijection-general2-ineq1} } =   2 \|\bOmega \B \|_\F^2 + 2\|\B^{1/2} \bOmega \B^{1/2} \|_\F^2  + \|\theta_B\|_\F^2 + \|\B \D_2^\top\|_\F^2 + \|\D_1 \B\|_\F^2 \\
		& \leq 2 \sigma^2_1(\X)  \bar{g}^{r}_{(\U,\B,\V)}(\theta_{(\U,\B,\V)}, \theta_{(\U,\B,\V)}).
	\end{split}
\end{equation*}
 This finishes the proof of this proposition.
\quad $\blacksquare$

\begin{proposition}[Bijection Between $\cH_{(\U,\Y)} \widebar{\cM}_{r}^{q_3}$ and $T_\X \cM_{r}^e$] \label{prop: general-bijection3} Suppose $\U \in \st(r,p_1)$, $\Y \in \bbR^{p_2 \times r}_*$ and $\X = \U \Y^\top$ with top $r$ right singular subspace spanned by $\V$. For any $\theta_{(\U,\Y)} = [\theta_U^\top \quad \theta_Y^\top]^\top \in \cH_{(\U,\Y)} \widebar{\cM}_{r}^{q_3}$ and $\xi_\X =  [\U \quad \U_\perp] \begin{bmatrix}
			\S & \D_2^\top\\
			\D_1 & \0
		\end{bmatrix} [\V \quad \V_\perp]^\top \in T_{\X}\cM^e_{r}$, define
\begin{equation} \label{def: xi-theta-correspondence-general3}
	\begin{split}
		\xi^{\theta_{(\U,\Y)}}_{\X}&:= [\U \quad \U_\perp] \begin{bmatrix}
		\theta_Y^\top \V & \theta_Y^\top \V_\perp \\
		\U_\perp^\top \theta_U \Y^\top \V & \0
	\end{bmatrix}[\V \quad \V_\perp]^\top  \in T_{\X}\cM^e_{r}, \\
	\theta_{(\U,\Y)}^{\xi_\X} &:= [(\U_\perp \D_1 (\Y^\top \V)^{-1})^\top \quad (\V \S^\top + \V_\perp \D_2)^\top]^\top \in \cH_{(\U,\Y)} \widebar{\cM}_{r}^{q_3}.
	\end{split}
\end{equation}Then we can find a linear bijective mapping $\cL_{\U,\Y}^{r}$ between $\cH_{(\U,\Y)} \widebar{\cM}_{r}^{q_3}$ and $T_\X \cM_{r}^e$,
\begin{equation*}
\begin{split}
	&\cL_{\U,\Y}^{r}: \theta_{(\U,\Y)} \in  \cH_{(\U,\Y)} \widebar{\cM}_{r}^{q_3} \longrightarrow \xi_{\X}^{\theta_{(\U,\Y)}} \in T_\X \cM_{r}^e, \quad (\cL_{\U,\Y}^{r})^{-1}: \xi_\X \in T_\X \cM_{r}^e \to \theta_{(\U,\Y)}^{\xi_\X} \in \cH_{(\U,\Y)} \widebar{\cM}_{r}^{q_3}
\end{split}
\end{equation*}such that $\cL_{\U,\Y}^{r}(\theta_{(\U,\Y)}) = \U \theta_Y^\top + \theta_U \Y^\top$ holds for any $\theta_{(\U,\Y)} \in \cH_{(\U,\Y)} \widebar{\cM}_{r}^{q_3}$.	

Finally, we have the following spectrum bounds for $\cL_{\U,\Y}^{r}$: for all $\theta_{(\U,\Y)} \in \cH_{(\U,\Y)} \widebar{\cM}_{r}^{q_3}$,
		\begin{equation} \label{ineq: bijection-spectrum-general3}
		\begin{split}
			(\sigma_r(\W^{-1}_\Y) \wedge \sigma_r^2(\Y \V_\Y^{-1/2} )) \bar{g}^{r}_{(\U,\Y)}(\theta_{(\U,\Y)}, \theta_{(\U,\Y)}) \leq & \|\cL^{r}_{\U,\Y}(\theta_{(\U,\Y)})\|_\F^2\\
			 & \quad \leq ( \sigma_1(\W_\Y^{-1}) \vee \sigma^2_1(\Y\V_\Y^{-1/2})   )\bar{g}^{r}_{(\U,\Y)}(\theta_{(\U,\Y)}, \theta_{(\U,\Y)}).
		\end{split}
	\end{equation}	
\end{proposition}
{\noindent \bf Proof of Proposition \ref{prop: general-bijection3}.} First, it is easy to see $\xi^{\theta_{(\U,\Y)}}_{\X}$ and $\theta_{(\U,\Y)}^{\xi_\X}$ are well defined given any $\theta_{(\U,\Y)} \in \cH_{(\U,\Y)} \widebar{\cM}_{r}^{q_3}$ and $\xi_\X \in T_{\X}\cM^e_{r}$. 

Next, we show $\cL_{\U,\Y}^{r}$ is a bijection. Notice $\cH_{(\U,\Y)} \widebar{\cM}_{r}^{q_3}$ is of dimension $(p_1+p_2-r)r$, which is the same with $T_\X \cM_{r}^e$. Suppose $\cL_{\U,\Y}^{r'} : \xi_\X \in T_\X\cM^e_{r} \longrightarrow \theta_{(\U,\Y)}^{\xi_\X} \in \cH_{(\U,\Y)} \widebar{\cM}_{r}^{q_3}$. For any $\xi_\X = [\U \quad \U_\perp] \begin{bmatrix}
			\S & \D_2^\top\\
			\D_1 & \0
		\end{bmatrix} [\V \quad \V_\perp]^\top \in T_\X\cM^e_{r}$, we have
\begin{equation} \label{eq: bijection-general3}
	\cL^{r}_{\U,\Y} ( \cL_{\U,\Y}^{r'}(\xi_\X) ) = \cL^{r}_{\U,\Y} (\theta_{(\U,\Y)}^{\xi_\X}) =  [\U \quad \U_\perp] \begin{bmatrix}
		\theta_Y^{\xi_\X \top}\V & \theta^{\xi_\X \top}_Y \V_\perp \\
		\U_\perp^\top \theta_U^{\xi_\X} \Y^\top \V & \0
	\end{bmatrix}[\V \quad \V_\perp]^\top = \xi_\X.
\end{equation} Since $\cL^{r}_{\U,\Y}$ and $\cL_{\U,\Y}^{r'}$ are linear maps, \eqref{eq: bijection-general3} implies $\cL^{r}_{\U,\Y}$ is a bijection and $\cL_{\U,\Y}^{r'} = (\cL^{r}_{\U,\Y})^{-1}$. At the same time, it is easy to check $\cL_{\U,\Y}^{r}(\theta_{(\U,\Y)}) = \U \theta_Y^\top + \theta_U \Y^\top$ holds for any $\theta_{(\U,\Y)} \in \cH_{(\U,\Y)} \widebar{\cM}_{r}^{q_3}$ by observing $P_{\U_\perp} \theta_U = \theta_U$ and $\Y$ lies in the column space of $\V$.

Next, we provide the spectrum bounds for $\cL^{r}_{\U,\Y}$. For any $\theta_{(\U,\Y)} = [(\U_\perp \D)^\top \quad \theta_Y^\top]^\top \in \cH_{(\U,\Y)} \widebar{\cM}_{r}^{q_3}$, we have $\bar{g}^{r}_{(\U,\Y)}(\theta_{(\U,\Y)},\theta_{(\U,\Y)})  = \|\D\V_\Y^{1/2}\|_\F^2 + \|\theta_Y \W_\Y^{1/2} \|_\F^2$. Thus, we have
\begin{equation*}
	\begin{split}
		\|\cL_{\U,\Y}^r(\theta_{(\U,\Y)})\|_\F^2 &= \|\theta_Y^\top \V\|_\F^2 + \|\theta_Y^\top \V_\perp\|_\F^2 + \|\U_\perp^\top \theta_U \Y^\top \V\|_\F^2 \\
		& = \|\theta_Y\|_\F^2 +  \|\D \Y^\top \V\|_\F^2\\
		& =   \|\theta_Y \W_\Y^{1/2} \W_\Y^{-1/2}  \|_\F^2 +  \|\D \V_\Y^{1/2}  \V_\Y^{-1/2} \Y^\top \V\|_\F^2 \\
		&   \overset{(a)}\geq \|\D\V_\Y^{1/2}\|_\F^2 \sigma^2_r(\Y\V_\Y^{-1/2}) +  \|\theta_Y \W_\Y^{1/2} \|_\F^2 \sigma_r(\W_\Y^{-1}) \\
		&\geq (\sigma_r(\W^{-1}_\Y) \wedge \sigma_r^2(\Y \V_\Y^{-1/2} ))\bar{g}^{r}_{(\U,\Y)}(\theta_{(\U,\Y)}, \theta_{(\U,\Y)}),
	\end{split}
\end{equation*} where (a) is because $\Y$ lies in the column space of $\V$, and
\begin{equation*}
	\begin{split}
		\|\cL_{\U,\Y}^r(\theta_{(\U,\Y)})\|_\F^2 & = \|\theta_Y\|_\F^2 +  \|\D \Y^\top \V\|_\F^2 \leq ( \sigma_1(\W_\Y^{-1}) \vee \sigma^2_1(\Y\V_\Y^{-1/2})   )\bar{g}^{r}_{(\U,\Y)}(\theta_{(\U,\Y)}, \theta_{(\U,\Y)}). \quad \blacksquare
	\end{split}
\end{equation*}

Now, we are ready to present our main results on the geometric landscape connection of Riemannian fixed-rank matrix optimization \eqref{eq: general prob} under the embedded and quotient geometries.

\begin{theorem}[Geometric Landscape Connections of \eqref{eq: general prob} on $\cM_{r}^e$ and $\cM_{r}^{q_1}$] \label{th: embedded-quotient-connection-general1}
	  Suppose the conditions in Proposition \ref{prop: general-bijection1} hold and the $\W_{\L,\R}$ and $\V_{\L,\R}$ in $\bar{g}^r_{(\L,\R)}$ satisfy $\W_{\L,\R} = \M \W_{\L\M,\R \M^{-\top}} \M^\top$, $\V_{\L,\R} = \M^{-\top } \V_{\L\M,\R \M^{-\top}} \M^{-1}$ for any $\M \in \GL(r)$. Then
	\begin{equation} \label{eq: gradient-connect-general1}
	\begin{split}
		\grad f(\X) &= \overline{\grad_{\L}\, h_{r}([\L,\R])} \W_{\L,\R} \R^\dagger  + ( \overline{\grad_{\R}\, h_{r}([\L,\R])} \V_{\L,\R} \L^\dagger )^\top (\I_{p_2} - \R \R^\dagger) \\
		 \overline{\grad\, h_{r}([\L,\R])} &= \begin{bmatrix}
		 	\grad f(\X) \R \W_{\L,\R}^{-1} \\
		 	(\grad f(\X))^\top \L \V_{\L,\R}^{-1}
		 \end{bmatrix}.
	\end{split}
	\end{equation}
Furthermore, if $[\L,\R]$ is a Riemannian FOSP of \eqref{eq: general-opt-problem-quotient-sub1}, we have: 
\begin{equation} \label{eq: Hessian-connection-general1}
	\begin{split}
		\overline{\Hess \, h_{r}([\L,\R])}[\theta_{(\L,\R)}, \theta_{(\L,\R)}]  &= \Hess f(\X)[\cL_{\L,\R}^{r}(\theta_{(\L,\R)}),\cL_{\L,\R}^{r}(\theta_{(\L,\R)})], \quad \forall \theta_{(\L,\R)} \in \cH_{(\L,\R)} \widebar{\cM}_{r}^{q_1}.
	\end{split}
	\end{equation}
	Finally, $\overline{\Hess \, h_{r}([\L,\R])}$ and $\Hess f(\X)$ have $(p_1 + p_2-r)r$ eigenvalues and for $i = 1,\ldots, (p_1 + p_2-r)r$, we have
	 $\lambda_i(\overline{\Hess \, h_{r}([\L,\R])})$ is sandwiched between $\overline{\gamma} \lambda_i(\Hess f(\X)) $ and $2\underline{\Gamma} \lambda_i(\Hess f(\X)) $, where $\overline{\gamma} = \sup_{\L, \R: \L \R^\top = \X} \gamma_{\L,\R}$ and $\underline{\Gamma}  = \inf_{\L, \R: \L \R^\top = \X} \Gamma_{\L,\R}$ where $\gamma_{\L,\R}$ and $\Gamma_{\L,\R}$ are given in \eqref{ineq: bijection-spectrum-general1}.
\end{theorem}
{\noindent \bf Proof of Theorem \ref{th: embedded-quotient-connection-general1}.} First, recall $\U$ and $\V$ span the top $r$ left and right singular subspaces of $\L\R^\top$, respectively and we have $\L \L^\dagger = P_\U$, $\R \R^\dagger = P_\V$. So \eqref{eq: gradient-connect-general1} is by direct calculation from the gradient expressions in Proposition \ref{prop: gradient-hessian-exp-general}. 

Next, we prove \eqref{eq: Hessian-connection-general1}. Since $[\L,\R]$ is a Riemannian FOSP of \eqref{eq: general-opt-problem-quotient-sub1}, we have
\begin{equation} \label{eq: FOSP-condition-general1}
	\overline{ \grad \, h_{r}([\L,\R])} = \0,\quad \nabla f(\L \R^\top) \R = \0 \quad \text{ and }  \left(\nabla f(\L \R^\top) \right)^\top \L = \0
\end{equation} and have $\grad f(\X) = \0$ by \eqref{eq: gradient-connect-general1}. So $ \nabla f(\X) = P_{\U_\perp} \nabla f(\X) P_{\V_\perp}$. For any $\L, \R$ such that $\L \R^\top = \X$. Recall $\P_1 = \U^\top \L$, $\P_2 = \V^\top \R$, and let $\bSigma = \U^\top \X \V$. Given any $\theta_{(\L,\R)} =  [\theta_L^\top \quad \theta_R^\top]^\top \in \cH_{(\L,\R)} \widebar{\cM}_{r}^{q_1}$, we have
\begin{equation} \label{eq: Hessian-con-gradient-general1}
	\langle \nabla f(\X), P_{\U_\perp} \theta_L \P_2^{\top} \bSigma^{-1} \P_1 \theta_R^\top  P_{\V_\perp}    \rangle = \langle \nabla f(\X), \theta_L \theta_R^\top \rangle,
\end{equation} where the equality is because $\P_1,\P_2$ are  nonsingular, $\P_1 \P_2^\top = \bSigma$ and $ \nabla f(\X) = P_{\U_\perp} \nabla f(\X) P_{\V_\perp} $.

Then by Proposition \ref{prop: gradient-hessian-exp-general}:
\begin{equation*}
	\begin{split}
		 &\quad \overline{\Hess \, h_{r}([\L,\R])}[\theta_{(\L,\R)}, \theta_{(\L,\R)}] \\
		 &=  \nabla^2 f(\L \R^\top)[\L \theta_R^\top + \theta_L \R^\top,\L \theta_R^\top + \theta_L \R^\top] + 2 \langle \nabla f(\L \R^\top), \theta_L \theta_R^\top \rangle \\
		 & \quad + \langle \nabla f(\L \R^\top) \R \rmD \W_{\L,\R}^{-1}[\theta_{(\L,\R)}] , \theta_L \W_{\L,\R} \rangle + \langle  \left(\nabla f(\L \R^\top) \right)^\top \L \rmD \V_{\L,\R}^{-1}[\theta_{(\L,\R)}], \theta_R \V_{\L,\R} \rangle \\
		 & \quad +  \langle \rmD \W_{\L,\R}[ \overline{ \grad \, h_{r}([\L,\R])} ], \theta_L^\top \theta_L \rangle /2 +  \langle \rmD \V_{\L,\R}[ \overline{ \grad \, h_{r}([\L,\R])} ], \theta_R^\top \theta_R \rangle /2\\
		 & \overset{ \eqref{eq: FOSP-condition-general1} } =  \nabla^2 f(\L \R^\top)[\L \theta_R^\top + \theta_L \R^\top,\L \theta_R^\top + \theta_L \R^\top] + 2 \langle \nabla f(\L \R^\top), \theta_L \theta_R^\top \rangle \\
		 & \overset{ \text{Proposition } \ref{prop: general-bijection1}, \eqref{eq: Hessian-con-gradient-general1} } =   \nabla^2 f(\X)[\cL^r_{\L,\R}(\theta_{(\L,\R)}),\cL^r_{\L,\R}(\theta_{(\L,\R)})] + 2 \langle \nabla f(\X), P_{\U_\perp} \theta_L \P_2^{\top} \bSigma^{-1} \P_1 \theta_R^\top  P_{\V_\perp}    \rangle\\
		 & = \Hess f(\X)[\cL^r_{\L,\R}(\theta_{(\L,\R)}),\cL^r_{\L,\R}(\theta_{(\L,\R)})],
	\end{split}
\end{equation*} where the last equality follows from the expression of $\Hess f(\X)$ in \eqref{eq: embedded-gd-hessian-general} and the definition of $\cL^r_{\L,\R}$. 

Then, by \eqref{ineq: bijection-spectrum-general1}, \eqref{eq: Hessian-connection-general1} and Theorem \ref{th: hessian-sandwich}, we have $\overline{\Hess \, h_{r}([\L,\R])}$ and $\Hess f(\X)$ have $(p_1+p_2- r)r$ eigenvalues and $\widebar{\lambda}_i(\Hess \, h_{r}([\L,\R]))$ is sandwiched between $\gamma_{\L,\R} \lambda_i(\Hess f(\X)) $ and $2\Gamma_{\L,\R} \lambda_i(\Hess f(\X)) $ for $i = 1,\ldots,(p_1+p_2- r)r$, where $\gamma_{\L,\R}$ and $\Gamma_{\L,\R}$ are given in \eqref{ineq: bijection-spectrum-general1}. Since the same sandwich relationship holds for all $(\L, \R)$ such that $\L \R^\top = \X$, we have
	 $\lambda_i(\overline{\Hess \, h_{r}([\L,\R])})$ is sandwiched between $\overline{\gamma} \lambda_i(\Hess f(\X)) $ and $2\underline{\Gamma} \lambda_i(\Hess f(\X)) $. \quad $\blacksquare$

\begin{theorem}[Geometric Landscape Connections of \eqref{eq: general prob} on $\cM_{r}^e$ and $\cM_{r}^{q_2}$] \label{th: embedded-quotient-connection-general2}
	 Suppose $\U \in \st(r,p_1)$, $\B \in \bbS_+(r)$, $\V \in \st(p_2,r)$ and $\X = \U \B \V^\top$. Then
	\begin{equation} \label{eq: gradient-connect-general2}
	\begin{split}
		\grad f(\X) &= P_{\U_\perp} \overline{\grad_\U\, h_{r}([\U,\Y])} \B^{-1} \V^\top + \U \bDelta_1 \V^\top + (P_{\V_\perp} \overline{\grad_\V\, h_{r}([\U,\B,\V])} \B^{-1} \U^\top)^\top 
			\end{split}
	\end{equation} where $\bDelta_1$ is uniquely determined by the equation system:
	\begin{equation*}
		\begin{split}
			 &\B\skew^\top(\bDelta_1 ) + \skew^\top(\bDelta_1 ) \B = 2 \overline{\grad_\U\, h_{r}([\U,\Y])}^\top \U,\\
			  &\sym(\bDelta_1) = \B^{-1} \overline{\grad_\B\, h_{r}([\U,\Y])} \B^{-1},
		\end{split}
	\end{equation*} and
	\begin{equation}  \label{eq: gradient-connect-general2-secondone}
		\begin{split}
			 \overline{\grad\, h_{r}([\U,\B,\V])} &= \begin{bmatrix}
		 	P_{\U_\perp}\grad f(\X) \V \B + \U ( \skew(\bDelta_2) \B + \B \skew(\bDelta_2)  )/2 \\
		 	\B \sym(\bDelta_2) \\
		 	(\B \U^\top \grad f(\X) P_{\V_\perp} )^\top - \V ( \skew(\bDelta_2) \B + \B \skew(\bDelta_2)  )/2
		 \end{bmatrix},
		\end{split}
	\end{equation} where $\bDelta_2 = \U^\top \grad f(\X) \V$.
	
Furthermore, if $[\U,\B,\V]$ is a Riemannian FOSP of \eqref{eq: general-opt-problem-quotient-sub2}, we have: for any $ \theta_{(\U,\B,\V)} \in \cH_{(\U,\B,\V)} \widebar{\cM}_{r}^{q_2}$,
\begin{equation} \label{eq: Hessian-connection-general2}
	\begin{split}
		\overline{\Hess \, h_{r}([\U,\B,\V])}[\theta_{(\U,\B,\V)}, \theta_{(\U,\B,\V)}]  &= \Hess f(\X)[\cL_{\U,\B,\V}^{r}(\theta_{(\U,\B,\V)}),\cL_{\U,\B,\V}^{r}(\theta_{(\U,\B,\V)})].
	\end{split}
	\end{equation}
	Finally, $\overline{\Hess \, h_{r}([\U,\B,\V])}$ has $(p_1 + p_2-r)r$ eigenvalues and for $i = 1,\ldots, (p_1 + p_2-r)r$, we have
	 $\lambda_i(\overline{\Hess \, h_{r}([\U,\B,\V])})$ is sandwiched between $\sigma^2_r(\X)\lambda_i(\Hess f(\X)) $ and $2\sigma^2_1(\X) \lambda_i(\Hess f(\X)) $.
\end{theorem}
{\noindent \bf Proof of Theorem \ref{th: embedded-quotient-connection-general2}.} First, recall $\bDelta = \U^\top \nabla f(\U \B \V^\top) \V$ and notice $\B\skew^\top(\bDelta_1 ) + \skew^\top(\bDelta_1 ) \B = 2 \overline{\grad_\U\, h_{r}([\U,\Y])}^\top \U = \B\skew^\top(\bDelta ) + \skew^\top(\bDelta ) \B$ is a Sylvester equation with respect $\skew(\bDelta_1)^\top$, which has a unique solution as $\B,-\B$ have disjoint spectra and \cite[Theorem VII.2.1]{bhatia2013matrix}. Since $\bDelta_1 = \skew(\bDelta_1) + \sym(\bDelta_1)$, $\bDelta_1$ is uniquely determined by the equation system $ \B\skew^\top(\bDelta_1 ) + \skew^\top(\bDelta_1 ) \B = 2 \overline{\grad_\U\, h_{r}([\U,\Y])}^\top \U, \sym(\bDelta_1) = \B^{-1} \overline{\grad_\B\, h_{r}([\U,\Y])} \B^{-1}$. Finally, because $\bDelta$ is a solution to this equation system, we have $\bDelta_1 = \bDelta = \U^\top \nabla f(\U \B \V^\top) \V$. The rest proofs of \eqref{eq: gradient-connect-general2} and \eqref{eq: gradient-connect-general2-secondone} are by direct calculation from the gradient expressions in Proposition \ref{prop: gradient-hessian-exp-general}. 

Next, we prove \eqref{eq: Hessian-connection-general2}. Since $[\U,\B,\V]$ is a Riemannian FOSP of \eqref{eq: general-opt-problem-quotient-sub2}, we have
\begin{equation} \label{eq: FOSP-condition-general2}
	\overline{ \grad \, h_{r}([\U,\B,\V])} = \0,\quad \grad f(\X) \overset{\eqref{eq: gradient-connect-general2}}= \0, \quad \nabla f(\U \B \V^\top) \V = \0 \,\text{ and }\,   \U^\top\nabla f(\U \B \V^\top)  = \0.
\end{equation} So $ \nabla f(\X) = P_{\U_\perp} \nabla f(\X) P_{\V_\perp}$. Given any $\theta_{(\U,\B,\V)} =  [\theta_U^\top \quad \theta_B^\top \quad \theta_V^\top]^\top \in \cH_{(\U,\B,\V)} \widebar{\cM}_{r}^{q_2}$, we have
\begin{equation} \label{eq: Hessian-con-gradient-general2}
	\langle \nabla f(\X), P_{\U_\perp} \theta_U \B \B^{-1} \B \theta_V^\top  P_{\V_\perp}    \rangle = \langle \nabla f(\X), \theta_U \B \theta_V^\top \rangle,
\end{equation} where the equality is because $ \nabla f(\X) = P_{\U_\perp} \nabla f(\X) P_{\V_\perp} $.

Then by Proposition \ref{prop: gradient-hessian-exp-general} and recall $\bDelta' =  \theta_U^\top\nabla f(\U \B \V^\top) \V$ and $\bDelta'' =  \U^\top\nabla f(\U \B \V^\top) \theta_V$, we have
\begin{equation*}
	\begin{split}
		 &\quad \overline{\Hess \, h_{r}([\U,\B,\V])}[\theta_{(\U,\B,\V)}, \theta_{(\U,\B,\V)}] \\
		 &=  \nabla^2 f(\U \B \V^\top)[\theta_U \B \V^\top + \U \theta_B \V^\top + \U \B \theta_V^\top, \theta_U \B \V^\top + \U \theta_B \V^\top + \U \B \theta_V^\top] + 2\langle \nabla f(\U\B\V^\top), \theta_U \B \theta_V^\top \rangle\\
		& +   \left\langle \bDelta,  \sym(\U^\top \theta_U \U^\top\theta_U) \B + \B \sym(\V^\top \theta_V \U^\top \theta_U) -2\theta_U^\top \theta_U \B \right\rangle/2\\
		& + \left\langle \bDelta,  \B\sym(\V^\top \theta_V \V^\top\theta_V) + \sym(\U^\top \theta_U \V^\top \theta_V) \B -2\B\theta_V^\top \theta_V +  2\theta_B \B^{-1} \theta_B \right\rangle/2\\
		& + \langle \bDelta', 2\theta_B -  \U^\top \theta_U \B - \theta_U^\top \U \B/2 -\V^\top \theta_V \B/2 \rangle + \langle \bDelta'', 2\theta_B -  \B\theta_V^\top \V - \B\V^\top \theta_V/2 -\B \theta_U^\top \U /2 \rangle,\\
		 & \overset{ \eqref{eq: FOSP-condition-general2} } =  \nabla^2 f(\U \B \V^\top)[\theta_U \B \V^\top + \U \theta_B \V^\top + \U \B \theta_V^\top, \theta_U \B \V^\top + \U \theta_B \V^\top + \U \B \theta_V^\top] + 2\langle \nabla f(\U\B\V^\top), \theta_U \B \theta_V^\top \rangle \\
		 & \overset{ \text{Proposition } \ref{prop: general-bijection2}, \eqref{eq: Hessian-con-gradient-general2} } =   \nabla^2 f(\X)[\cL^r_{\U,\B,\V}(\theta_{(\U,\B,\V)}),\cL^r_{\U,\B,\V}(\theta_{(\U,\B,\V)})] + 2\langle \nabla f(\X), P_{\U_\perp} \theta_U \B \B^{-1} \B \theta_V^\top  P_{\V_\perp}    \rangle\\
		 & = \Hess f(\X)[\cL^r_{\U,\B,\V}(\theta_{(\U,\B,\V)}),\cL^r_{\U,\B,\V}(\theta_{(\U,\B,\V)})],
	\end{split}
\end{equation*} where the last equality follows from the expression of $\Hess f(\X)$ in \eqref{eq: embedded-gd-hessian-general} and the definition of $\cL^r_{\U,\B,\V}$. 

Then, by \eqref{ineq: bijection-spectrum-general2}, \eqref{eq: Hessian-connection-general2} and Theorem \ref{th: hessian-sandwich}, we have $\overline{\Hess \, h_{r}([\U,\B,\V])}$ has $(p_1+p_2- r)r$ eigenvalues and  $\widebar{\lambda}_i(\Hess \, h_{r}([\U,\B,\V]))$ is sandwiched between $\sigma^2_r(\X)\lambda_i(\Hess f(\X)) $ and $2\sigma^2_1(\X) \lambda_i(\Hess f(\X)) $ for $i = 1,\ldots,(p_1+p_2- r)r$. \quad $\blacksquare$

\begin{theorem}[Geometric Landscape Connections of \eqref{eq: general prob} on $\cM_{r}^e$ and $\cM_{r}^{q_3}$] \label{th: embedded-quotient-connection-general3}
	  Suppose the conditions in Proposition \ref{prop: general-bijection3} hold and the $\V_\Y,\W_\Y$ in $\bar{g}^r_{(\U,\Y)}$ satisfies $\V_\Y = \O \V_{\Y\O} \O^\top$ and $\W_\Y = \O \W_{\Y\O} \O^\top$ for any $\O \in \bbO_r$. Then
	\begin{equation} \label{eq: gradient-connect-general3}
	\begin{split}
		\grad f(\X) &= \overline{\grad_\U\, h_{r}([\U,\Y])}\V_\Y \Y^\dagger  + \left(\overline{\grad_\Y\, h_{r}([\U,\Y])} \W_\Y \U^\top \right)^\top  \\
		 \overline{\grad\, h_{r}([\U,\Y])} &= \begin{bmatrix}
		 	P_{\U_\perp}\grad f(\X) \Y \V_\Y^{-1} \\
		 	(\grad f(\X))^\top \U \W_\Y^{-1}
		 \end{bmatrix}.
	\end{split}
	\end{equation}
Furthermore, if $[\U,\Y]$ is a Riemannian FOSP of \eqref{eq: general-opt-problem-quotient-sub3}, we have: 
\begin{equation} \label{eq: Hessian-connection-general3}
	\begin{split}
		\overline{\Hess \, h_{r}([\U,\Y])}[\theta_{(\U,\Y)}, \theta_{(\U,\Y)}]  &= \Hess f(\X)[\cL_{\U,\Y}^{r}(\theta_{(\U,\Y)}),\cL_{\U,\Y}^{r}(\theta_{(\U,\Y)})], \quad \forall \theta_{(\U,\Y)} \in \cH_{(\U,\Y)} \widebar{\cM}_{r}^{q_3}.
	\end{split}
	\end{equation}
	Finally, $\overline{\Hess \, h_{r}([\U,\Y])}$ has $(p_1 + p_2-r)r$ eigenvalues and for $i = 1,\ldots, (p_1 + p_2-r)r$, we have
	 $\lambda_i(\overline{\Hess \, h_{r}([\U,\Y])})$ is sandwiched between $(\sigma_r(\W^{-1}_\Y) \wedge \sigma_r^2(\Y \V_\Y^{-1/2} )) \lambda_i(\Hess f(\X)) $ and $( \sigma_1(\W_\Y^{-1}) \vee \sigma^2_1(\Y\V_\Y^{-1/2})   ) \lambda_i(\Hess f(\X)) $.
\end{theorem}
{\noindent \bf Proof of Theorem \ref{th: embedded-quotient-connection-general3}.} The proof is similar to the proof of Theorem \ref{th: embedded-quotient-connection-general1} and is postponed to Appendix \ref{sec: additional-proofs-general}. \quad $\blacksquare$

In Table \ref{tab: illustraion-gap-coefficient-general}, we provide explicit values of gap coefficients of the sandwich inequalities in Theorems \ref{th: embedded-quotient-connection-general1}, \ref{th: embedded-quotient-connection-general2} and \ref{th: embedded-quotient-connection-general3} under several different metrics specified by $(\W_{\L,\R}, \V_{\L,\R})$, $(\V_\B,\W_\B)$ and $(\V_\Y, \W_\Y)$ in $\bar{g}^{r}_{(\L,\R)}$, $\bar{g}^{r}_{(\U,\B,\V)}$, and $\bar{g}^{r}_{(\U,\Y)}$, respectively.
\begin{table}[ht]
	\centering
	\begin{tabular}{c | c | c | c}
		\hline
		 & \multirow{2}{15em}{Choices of $(\W_{\L,\R}, \V_{\L,\R})$, $(\V_\B, \W_\B)$, and $(\V_\Y, \W_\Y)$ in $\bar{g}^{r}$} & \multirow{2}{8em}{Gap Coefficient Lower Bound} & \multirow{2}{8em}{Gap Coefficient Upper Bound} \\
		 & & & \\
		 \hline
		 \multirow{2}{7em}{$\cM_{r}^e$ v.s. $\cM^{q_1}_{r1}$} & $\W_{\L,\R} = (\L^\top \L)^{-1}$,  $\V_{\L,\R} = (\R^\top \R)^{-1}$  & $\sigma^2_r(\X)$  & $2 \sigma^2_1(\X)$ \\
		 \cline{2-4} 
		 & $\W_{\L,\R} = \R^\top \R$,  $\V_{\L,\R} = \L^\top \L$ & $1$ & $2$ \\
		 \hline
		 \multirow{1}{7em}{$\cM_{r}^e$ v.s. $\cM^{q_2}_{r}$} & $\V_\B = \I_r$, $\W_\B = \B^{-1}$  & $\sigma^2_r(\X)$  & $ 2\sigma^2_1(\X)$ \\
		  \hline
		 \multirow{3}{7em}{$\cM_{r}^e$ v.s. $\cM^{q_3}_{r}$} & $\V_\Y = \I_r$, $\W_\Y = \I_r$  & $\sigma^2_r(\X) \wedge 1 $  & $\sigma^2_1(\X) \vee  1$\\
		 \cline{2-4}
		 & $\V_\Y = \I_r, \W_\Y = (\Y^\top \Y)^{-1}$  & $ \sigma^2_r(\X) $  & $\sigma^2_1(\X)$\\
		 \cline{2-4}
		 & $\V_\Y = \Y^\top \Y, \W_\Y = \I_r$  & $1 $  & $1$\\
		 \hline
	\end{tabular}
	\caption{Gap coefficients in the sandwich inequalities in Theorems \ref{th: embedded-quotient-connection-general1}, \ref{th: embedded-quotient-connection-general2}, and \ref{th: embedded-quotient-connection-general3} under different Riemannian metrics. 
	} \label{tab: illustraion-gap-coefficient-general}
\end{table}

By Theorems \ref{th: embedded-quotient-connection-general1}, \ref{th: embedded-quotient-connection-general2} and \ref{th: embedded-quotient-connection-general3}, we have the following Corollary \ref{coro: landscape connection general case} on the equivalence of Riemannian FOSPs, SOSPs and strict saddles of optimization \eqref{eq: general prob} under the embedded and the quotient geometries. The proof is given in Appendix \ref{sec: additional-proofs-general}.
\begin{corollary}({\bf Equivalence on Riemannian FOSPs, SOSPs and strict saddles of \eqref{eq: general prob} Under Embedded and Quotient Geometries}) \label{coro: landscape connection general case} Suppose $\W_{\L,\R} = \M \W_{\L\M,\R \M^{-\top}} \M^\top$, $\V_{\L,\R} = \M^{-\top } \V_{\L\M,\R \M^{-\top}} \M^{-1}$ hold for any $\M \in \GL(r)$ and $\V_\Y = \O \V_{\Y\O} \O^\top$, $\W_\Y = \O \W_{\Y\O} \O^\top$ holds for any $\O \in \bbO_r$. Then we have
\begin{itemize}
\item[(a)] given $\L \in \bbR^{p_1 \times r}_*, \R \in \bbR^{p_2 \times r}_*, \U \in \st(r,p_1)$, $\B \in \bbS_{+}(r)$, $\V \in \st(r,p_2)$ and $\Y \in \bbR^{p_2 \times r}_*$, if $[\L,\R]$ ($[\U,\B,\V]$ or $[\U,\Y]$) is a Riemannian FOSP or SOSP or strict saddle of \eqref{eq: general-opt-problem-quotient-sub1} (\eqref{eq: general-opt-problem-quotient-sub2} or \eqref{eq: general-opt-problem-quotient-sub3}), then $\X = \L \R^\top$ ($\X = \U \B \V^\top$ or $\X = \U \Y^\top$) is a Riemannian FOSP or SOSP or strict saddle of \eqref{eq: general prob} under the embedded geometry;
\item[(b)]if $\X$ is a Riemannian FOSP or SOSP or strict saddle of \eqref{eq: general prob} under the embedded geometry, then there is a unique $[\L,\R]$ ($[\U,\B,\V]$ or $[\U,\Y]$) such that $\L \R^\top = \X$ ($\U \B \V^\top = \X$ or $\U \Y^\top =\X$) and it is a Riemannian FOSP or SOSP or strict saddle of \eqref{eq: general-opt-problem-quotient-sub1} (\eqref{eq: general-opt-problem-quotient-sub2} or \eqref{eq: general-opt-problem-quotient-sub3}).
\end{itemize}
\end{corollary}

\begin{remark}({\bf Geometric Connection of Non-convex Factorization and Quotient Manifold Formulations for \eqref{eq: PSD-manifold-formulation} and \eqref{eq: general prob}})
As we have discussed in the Introduction, another popular approach for handling the rank constraint in \eqref{eq: PSD-manifold-formulation} or \eqref{eq: general prob} is via factorizing $\X$ into $\Y \Y^\top$ or $\L \R^\top$ and then treating the new problem as unconstrained optimization in the Euclidean space. In the recent work \cite{luo2021nonconvex}, they showed for both \eqref{eq: PSD-manifold-formulation} and \eqref{eq: general prob} the geometric landscapes under the factorization and embedded submanifold formulations are almost equivalent. By combining their results and the results in this paper, we also have a geometric landscape connection of \eqref{eq: PSD-manifold-formulation} and \eqref{eq: general prob} under the factorization and quotient manifold formulations.
\end{remark}

\begin{remark}({\bf Algorithmic Connection of Embedded and Quotient Geometries})\label{rem: algorithmic-connection}
	In contrast to the geometric connections of an optimization problem under two different geometries, the relationship of two sequences of iterates generated by algorithms in two geometries is more subtle. First, the algorithms under the quotient geometry are performed in the horizontal space and they depend on the quotient structure and the Riemannian metric we pick. Thus, it is hard to expect a universal algorithmic connection under two geometries. On the other hand, with some careful calculation, we find that by taking some specific metrics under the quotient geometry, there are indeed some connections of gradient flows, i.e., solutions of dynamical systems $d \X /d t = -\grad f(\X)$, $d [\Z] /d t = -\grad h([\Z])$ \cite[Chapter 9]{lee2013smooth}, produced under the embedded and quotient geometries. This is particularly true when the metrics are chosen in a way such that the sandwich gap coefficients in the geometric connection are some universal constants (see Remark \ref{rem: effiect-of-metric-on-landscape} and Tables \ref{tab: illustraion-gap-coefficient-PSD} and \ref{tab: illustraion-gap-coefficient-general}). Next, we illustrate a few such examples of gradient flows under various embedded and quotient geometries. Specifically, we compare the dynamical system $d \X /d t = -\grad f(\X)$ under the embedded geometry and the induced dynamical system $d \ell([\Z]) /dt = - \rmD \ell([\Z])[ \grad h([\Z]) ] = - \rmD  \bar{\ell}(\Z) [ \overline{\grad h([\Z])}  ]$ from the quotient geometry with the link $\X = \ell([\Z])$: 
	\begin{itemize}[leftmargin=*]
		\item (PSD case) Take $\W_\Y = 2\Y^\top \Y$ in $\bar{g}_\Y^{r+}$ and $\V_\B = 2\B^2, \W_\B = \I_r$ in $\bar{g}_{(\U,\B)}^{r+}$, we have the following gradient flows of $\{ \X_t \}$ in \eqref{eq: PSD-manifold-formulation} under $\cM_{r+}^e$, $\cM_{r+}^{q_1}$ and $\cM_{r+}^{q_2}$: 	
		 \begin{equation*}
		\begin{split}
			\text{under } \cM_{r+}^e: &\quad d \X /d t = -\grad f(\X) = -P_{\U} \nabla f(\X) - \nabla f(\X)P_{\U} + P_{\U} \nabla f(\X)P_{\U},\\
			\text{under } \cM_{r+}^{q_1}:& \quad d \ell([\Z]) /d t = -\overline{\grad\, h_{r+}([\Y])}\Y^\top - \Y \left(\overline{\grad\, h_{r+}([\Y])}\right)^\top = -P_{\U} \nabla f(\X) - \nabla f(\X)P_{\U};\\
			\text{under } \cM_{r+}^{q_2}:& \quad d \ell([\Z]) / dt = - \overline{\grad_\U\, h_{r+}([\U,\B])} \B \U^\top - \U \overline{\grad_\B\, h_{r+}([\U,\B])} \U^\top \\
			& \quad \quad \quad \quad \quad \quad  - \U \B (\overline{\grad_\U\, h_{r+}([\U,\B])})^\top \\
			& \quad \quad \quad \quad \quad \,  = -P_{\U} \nabla f(\X) - \nabla f(\X)P_{\U} + P_{\U} \nabla f(\X)P_{\U}.
		\end{split}
	\end{equation*}
	\item (general case) Take $\W_{\L,\R} = \R^\top \R, \V_{\L,\R} = \L^\top \L$ in $\bar{g}^r_{(\L,\R)}$, $\V_\Y = \Y^\top \Y, \W_\Y = \I_r$ in $\bar{g}^r_{(\U,\Y)}$, we have the following gradient flows of $\{\X_t \}$ in \eqref{eq: general prob} under $\cM_{r}^e$, $\cM_{r}^{q_1}$ and $\cM_{r}^{q_3}$:
\begin{equation*}		
\begin{split}
			\text{under } \cM_{r}^e: &\quad d \X /d t = -\grad f(\X) = -P_{\U} \nabla f(\X) - \nabla f(\X)P_{\V} + P_{\U} \nabla f(\X)P_{\V},\\
			\text{under } \cM_{r}^{q_1}: & \quad d \ell([\Z]) /d t = -\overline{\grad_{\L}\, h_{r}([\L,\R])}\R^\top - \L \left(\overline{\grad_{\R}\, h_{r}([\L,\R])}\right)^\top = -P_{\U} \nabla f(\X) - \nabla f(\X)P_{\V},\\
			\text{under } \cM_{r}^{q_3}: & \quad d \ell([\Z]) /d t = - \overline{\grad_\U\, h_{r}([\U,\Y])} \Y^\top - \U (\overline{\grad_\Y\, h_{r}([\U,\Y])})^\top \\
			& \quad \quad \quad \quad \quad  = -P_{\U} \nabla f(\X) - \nabla f(\X)P_{\V} + P_{\U} \nabla f(\X)P_{\V}.
		\end{split}
\end{equation*}
	\end{itemize}
		
Consider $\cM_{r+}^e$ v.s. $\cM_{r+}^{q_1}$ or $\cM_{r}^e$ v.s. $\cM_{r}^{q_1}$ with the specific metrics chosen above, we can see that the gradient flows of \eqref{eq: PSD-manifold-formulation} or \eqref{eq: general prob} under these embedded and the quotient geometries only differ one term, which has a magnitude smaller than the other terms. In addition, we know from Tables \ref{tab: illustraion-gap-coefficient-PSD} and \ref{tab: illustraion-gap-coefficient-general} that the sandwich coefficients lower and upper bounds are universal constants $1$ and $2$, respectively, in these settings. Some empirical evidence which shows the remarkably similar algorithmic performance under these embedded and quotient geometries was provided in \cite{mishra2012riemannian} and here our geometric connection results provide more theoretical insights for this empirical observation. More interestingly, we discover the gradient flows under embedded and quotient geometries are {\it identical} when the spectrum of Riemannian Hessians under two geometries coincide at Riemannian FOSPs, e.g., see $\cM_{r+}^e$ v.s. $\cM_{r+}^{q_2}$ or $\cM_{r}^e$ v.s. $\cM_{r}^{q_3}$ with the specific metrics chosen above.

\end{remark}

\section{Conclusion and Discussions} \label{sec: conclusion}
In this paper, we propose a general procedure for establishing geometric connections of a Riemannian optimization problem under embedded and quotient geometries. By applying the general procedure to problems \eqref{eq: PSD-manifold-formulation} and \eqref{eq: general prob} under the embedded and quotient geometries, we establish an exact Riemannian gradient connection under two geometries at every point on the manifold and sandwich inequalities between the spectra of the Riemannian Hessians at Riemannian FOSPs. These results immediately imply an equivalence on the sets of Riemannian FOSPs, SOSPs, and strict saddles of \eqref{eq: PSD-manifold-formulation} and \eqref{eq: general prob} under embedded and quotient geometries. Moreover, we observe an intriguing algorithmic connection for fixed-rank matrix optimization under two geometries with some specific Riemannian metrics. After the first release of this manuscript, there are a couple of new related works that are worth mentioning. In \cite{levin2022effect}, the authors studied the landscape connection of a pair of optimization problems, including a related bounded-rank matrix optimization problem, connected by a smooth parameterization. \cite{zheng2022riemannian} studied the algorithmic connection of Riemannian gradient descent type algorithms under the embedded geometry and a few quotient geometries in the PSD fixed-rank optimization problem. These works provide a good complement to the results in this paper. 

There are many interesting extensions to the results in this paper to be explored in the future. First, as we have mentioned in the Example \ref{ex: example-1} in Section \ref{sec: general-strategy-outline}, our results on the connection of Riemannian Hessians under the embedded and the quotient geometries are established at Riemannian FOSPs. It is interesting to explore whether it is possible to connect the landscapes under two geometries at non-stationary points. Second, other than the geometries covered in this paper, there are many other embedded and quotient geometries for fixed-rank matrices, such as the ones in \cite{absil2009geometric,grubivsic2007efficient,vandereycken2013riemannian} for $\cM_{r+}$, it will be interesting to study the geometric landscape connection of \eqref{eq: PSD-manifold-formulation} and \eqref{eq: general prob} under these geometries. Third, another common manifold that has both embedded and quotient representations is the Stiefel manifold \cite{edelman1998geometry}. We believe our general procedure in Section \ref{sec: general-strategy-for-connection} can also be used to establish the geometric connection in that setting as well. Finally, our ultimate goal is to better understand the connections and comparisons of different Riemannian geometries and give some guidelines on how to choose them given a Riemannian optimization problem. Some progress and discussions on how to choose Riemannian metrics in quotient geometries can be found in \cite{mishra2016riemannian,vandereycken2013riemannian}. While there is still not too much study on how to choose different quotient structures and manifold classes. It is important to explore these directions in the future.

\section*{Acknowledgments.} The authors would like to thank the editors and two reviewers for their great suggestions and comments. Y. Luo would like to thank Nicol{\'a}s Garc{\'i}a Trillos for helpful discussions during the project.

\appendix

\section{Additional Preliminaries on Riemannian Connection} \label{sec: additional-preliminaries}
In this section, we provide a more detailed discussion of the Riemannian connection as it is critical in deriving the Riemannian Hessians. First, a {\it vector field} $\xi$ on a manifold $\cM$ is a mapping that assigns each point $\X \in \cM$ a tangent vector $\xi_\X \in T_\X \cM$. The Riemannian gradient is a typical example of the vector field. 

 The {\it coordinate-free} definition of the Riemannian connection follows from the fundamental theorem of differential geometry \cite[Chapter 2.1]{petersen2006riemannian} and can be characterized by the so-called {\it Koszul formula} given as follows:
\begin{equation} \label{eq: general-koszul-formula}
	2g_\X (\widebar{\nabla}_\xi \eta, \theta ) = {\rm D} g_\X (\eta,\theta)[\xi] + \rmD g_\X(\xi,\theta)[\eta] - \rmD g_\X(\eta, \xi)[\theta] + g_\X( \theta, [\xi, \eta] ) + g_\X (\eta,[\theta,\xi]) - g_\X (\xi,[\eta, \theta]),
\end{equation} where $[\xi,\eta]$ denotes the {\it Lie bracket} between vector fields $\xi$ and $\eta$. If $\cM$ is an open subset of a vector space, the Lie bracket of vector fields is given by
\begin{equation} \label{eq: lie-bracket}
	[\xi, \eta] = \rmD \eta [\xi] - \rmD \xi [\eta],
\end{equation}
where $\rmD \eta [\xi]$ is the classical direction derivative of a vector field in a vector space and its evaluation at $\X$ is denoted by $(\rmD \eta [\xi])_\X$ and given as follows: $(\rmD \eta [\xi])_\X = \lim_{t \to 0} ( \eta_{\X + t \xi_\X} - \eta_\X )/t $ \cite[Section 3.5]{meyer2011geometric}.

In the following Proposition \ref{prop: koszul-formula}, we present the Koszul formulas associated with each simple manifold in Table \ref{tab: basic-prop-simple-manifold} and their combinations. These formulas will be used in deriving Riemannian Hessians under the quotient geometries given in Sections \ref{sec: quotient-PSD} and \ref{sec: quotient-general}.  
\begin{proposition}({\bf Koszul Formulas Associated with $\bbR^{p \times r}_*$, $\st(r,p)$, $\bbS_{+}(r)$ and their combinations}) \label{prop: koszul-formula}

(i) ({\bf Koszul Formulas Associated with $\bbR^{p \times r}_*$, $\st(r,p)$ and $\bbS_{+}(r)$.}) Consider manifolds $\bbR^{p \times r}_*$, $\st(r,p)$ and $\bbS_{+}(r)$ endowed metrics $g$ given in Table \ref{tab: basic-prop-simple-manifold} with $\V_\lozenge = \I$ and $\W_\B = \B^{-1}$ . For vector fields $\xi, \eta, \theta$ on $\bbR^{p \times r}_*$ or $\st(r,p)$ or $\bbS_{+}(r)$, we have the following Koszul formulas:
\begin{equation} \label{eq: koszul-formula-simple-manifolds}
	\begin{split}
	 \bbR^{p \times r}_*: \, 2 g_\Y ( \widebar{\nabla}_\xi \eta, \theta  ) &= 2 g_\Y ( \rmD \eta[\xi], \theta ) + \tr( \rmD \W_\Y [\xi] \eta^\top \theta ) + \tr( \rmD \W_\Y [\eta] \xi^\top \theta ) - \tr( \rmD \W_\Y [\theta] \eta^\top \xi ),\\
	 \st(r,p): \,	2 g_\U (\widebar{\nabla}_\xi \eta, \theta  ) & = 2 g_\U (\rmD \eta[\xi], \theta ),\\
	 \bbS_{+}(r): \,	2 g_\B (\widebar{\nabla}_\xi \eta, \theta ) & = 2 g_\B (\rmD \eta[\xi], \theta) - \tr(\B^{-1} \xi \B^{-1} \eta \B^{-1} \theta ) - \tr(\B^{-1} \eta \B^{-1}  \xi \B^{-1} \theta ).
	\end{split}
\end{equation}	Here, the operation on vector fields is performed elementwise on the associated tangent vectors. For example, $\eta^\top \theta$ is mapping that assign each point $\Y \in \bbR^{p \times r}_*$ to the matrix $\eta_\Y^\top \theta_\Y$, i.e., $(\eta^\top \theta)_\Y = \eta_\Y^\top \theta_\Y$.

(ii) ({\bf Koszul Formula Associated with $\st(r,p)\times \bbS_{+}(r)$.})
 Suppose we have metric $g_{(\U, \B)}$ on $\st(r,p)\times \bbS_{+}(r)$ given by $g_{(\U,\B)}(\eta_{(\U, \B)}, \theta_{(\U, \B)}) = \tr(\V_\B \eta_U^\top \theta_U) + \tr(\W_\B \eta_B \W_\B \theta_B)$ for $\eta_{(\U, \B)} = [\eta_U^\top \quad \eta_B^\top]^\top, \theta_{(\U, \B)} = [\theta_U^\top \quad \theta_B^\top]^\top \in T_\U \st(r,p) \times T_\B \bbS_+(r)$. Then for vector fields $\xi, \eta, \theta$ on $\st(r,p)\times \bbS_{+}(r)$, 
 \begin{equation}\label{eq: koszul-formula-st-s+}
 \begin{split}
 	2g_{(\U, \B)} ( \widebar{\nabla}_{\xi_{(\U', \B')}} \eta, \theta_{(\U', \B')} )&= 2 g_{(\U, \B)}( \rmD \eta_{(\U', \B')}[\xi_{(\U', \B')}], \theta_{(\U', \B')} ) \\
 	& + \tr( \rmD \V_\B[\xi_{B'}] \eta_{U'}^\top \theta_{U'} ) + 2\tr( \sym( \rmD \W_\B[\xi_{B'}] \eta_{B'} \W_\B) \theta_{B'}   ) \\
 	& + \tr( \rmD \V_\B[\eta_{B'}] \xi_{U'}^\top \theta_{U'} ) + 2\tr( \sym( \rmD \W_\B[\eta_{B'}] \xi_{B'} \W_\B   ) \theta_{B'}   )\\
 	& - \tr( \rmD \V_\B[\theta_{B'}] \eta_{U'}^\top \xi_{U'} ) - 2\tr( \sym( \rmD \W_\B[\theta_{B'}] \eta_{B'} \W_\B ) \xi_{B'}   ).
 \end{split}
 \end{equation}

(iii) ({\bf Koszul Formula Associated with $\st(r,p_1) \times \bbR^{p_2 \times r}_*$.})  Suppose we have metric $g_{(\U, \Y)}$ on $\st(r,p_1) \times \bbR^{p_2 \times r}_*$ given by $g_{(\U,\Y)}(\eta_{(\U, \Y)}, \theta_{(\U, \Y)}) = \tr(\V_\Y \eta_U^\top \theta_U) + \tr(\W_\Y \eta_Y^\top \theta_Y)$ for $\eta_{(\U, \Y)} = [\eta_U^\top \quad \eta_Y^\top]^\top, \theta_{(\U, \Y)} = [\theta_U^\top \quad \theta_Y^\top]^\top \in T_\U \st(r,p_1) \times T_\Y \bbR^{p_2 \times r}_*$. Then for vector fields $\xi, \eta, \theta$ on $\st(r,p_1) \times \bbR^{p_2 \times r}_*$, 
 \begin{equation} \label{eq: koszul-formula-st-R}
 \begin{split}
 	2g_{(\U, \Y)} ( \widebar{\nabla}_{\xi_{(\U', \Y')}} \eta, \theta_{(\U', \Y')} )&= 2 g_{(\U, \Y)}( \rmD \eta_{(\U', \Y')}[\xi_{(\U', \Y')}], \theta_{(\U', \Y')} ) \\
 	& + \tr( \rmD \V_\Y[\xi_{Y'}] \eta_{U'}^\top \theta_{U'} ) + \tr( \rmD \W_\Y [\xi_{Y'}] \eta_{Y'}^\top \theta_{Y'}  ) \\
 	& + \tr( \rmD \V_\Y[\eta_{Y'}] \xi_{U'}^\top \theta_{U'} ) + \tr( \rmD \W_\Y [\eta_{Y'}] \xi_{Y'}^\top \theta_{Y'}  )\\
 	& - \tr( \rmD \V_\Y[\theta_{Y'}] \eta_{U'}^\top \xi_{U'} ) - \tr( \rmD \W_\Y [\theta_{Y'}] \eta_{Y'}^\top \xi_{Y'}  ).
 \end{split}
 \end{equation}

\end{proposition}
{\noindent \bf Proof of Proposition \ref{prop: koszul-formula}.} {\bf (Part (i))} The Koszul formulas associated with $\st(r,p)$ and $\bbS_{+}(r)$ can be found in \cite[Section 3.5.1]{meyer2011geometric} and \cite[Appendix B]{meyer2011geometric}, respectively. Next, we derive the Koszul formula associated with the $\bbR^{p \times r}_*$ by using the general formula \eqref{eq: general-koszul-formula}. First, for any vector fields $\eta, \theta, \xi$ on $\bbR^{p \times r}_*$, we have
\begin{equation} \label{eq: directional-derivative-gY-instance}
	\begin{split}
		(\rmD g_\Y (\eta, \theta )[\xi])_{\Y'} &= \lim_{t \to 0} \frac{\tr( \W_{\Y + t\xi_{\Y'} } \eta_{\Y' + t \xi_{\Y'} }^\top \theta_{\Y' + t\xi_{\Y'} }  )  - \tr( \W_\Y \eta_{\Y'}^\top \theta_{\Y'} ) }{t} \\
		&=  \lim_{t \to 0} \frac{\tr( \W_{\Y + t\xi_{\Y'} } \eta_{\Y' + t \xi_{\Y'} }^\top \theta_{\Y' + t\xi_{\Y'} }  )  - \tr( \W_\Y \eta_{\Y' + t \xi_{\Y'} }^\top \theta_{\Y' + t\xi_{\Y'} }  ) }{t} \\
		& \quad + g_\Y( \rmD \eta_{\Y'}[\xi_{\Y'}], \theta_{\Y'} ) + g_\Y( \eta_{\Y'}, \rmD \theta_{\Y'}[\xi_{\Y'}] ) \\
		& = \tr(\rmD \W_\Y [\xi_{\Y'}] \eta_{\Y'}^\top \theta_{\Y'} ) + g_\Y( \rmD \eta_{\Y'}[\xi_{\Y'}], \theta_{\Y'} ) + g_\Y( \eta_{\Y'}, \rmD \theta_{\Y'}[\xi_{\Y'}] ).
	\end{split}
\end{equation} This proves
\begin{equation} \label{eq: directional-derivative-gY}
	\rmD g_\Y (\eta, \theta)[\xi] = \tr(\rmD \W_\Y [\xi] \eta^\top \theta ) + g_\Y( \rmD \eta[\xi], \theta ) + g_\Y( \eta, \rmD \theta[\xi] ).
\end{equation}

Then by \eqref{eq: general-koszul-formula} and the fact $\bbR^{p \times r}_*$ is an open subset of $\bbR^{p \times r}$, we have
\begin{equation} \label{eq: koszul-formula-derive-gY}
	\begin{split}
		&\quad 2 g_\Y (\widebar{\nabla}_\xi \eta, \theta) \\
		&= {\rm D} g_\Y (\eta,\theta)[\xi] + \rmD g_\Y(\xi,\theta)[\eta] - \rmD g_\Y(\eta, \xi)[\theta] + g_\Y( \theta, [\xi, \eta] ) + g_\Y (\eta,[\theta,\xi]) - g_\Y (\xi,[\eta, \theta])\\
		& \overset{ \eqref{eq: lie-bracket},\eqref{eq: directional-derivative-gY}} = \tr(\rmD \W_\Y [\xi] \eta^\top \theta ) + g_\Y( \rmD \eta[\xi], \theta ) + g_\Y( \eta, \rmD \theta[\xi] ) + \tr(\rmD \W_\Y [\eta] \xi^\top \theta ) + g_\Y( \rmD \xi[\eta], \theta ) \\
		& \quad   + g_\Y( \xi, \rmD \theta[\eta] ) - \tr(\rmD \W_\Y [\theta] \eta^\top \xi ) - g_\Y( \rmD \eta[\theta], \xi ) - g_\Y( \eta, \rmD \xi[\theta] ) + g_\Y(\theta, \rmD \eta[\xi] - \rmD \xi[\eta]) \\
		& \quad + g_\Y(\eta, \rmD \xi[\theta] - \rmD \theta[\xi]) - g_\Y(\xi, \rmD \theta[\eta] - \rmD \eta[\theta])\\
		& = 2 g_\Y ( \rmD \eta[\xi], \theta ) + \tr( \rmD \W_\Y [\xi] \eta^\top \theta ) + \tr( \rmD \W_\Y [\eta] \xi^\top \theta ) - \tr( \rmD \W_\Y [\theta] \eta^\top \xi ). 
	\end{split}
\end{equation}

{\bf (Part (ii))}. For vector fields $\xi, \eta, \theta$ on $\st(r,p)\times \bbS_{+}(r)$,
\begin{equation}\label{eq: directional-derivative-gUB-instance}
	\begin{split}
		&\rmD g_{(\U, \B)}(\eta_{(\U', \B')}, \theta_{(\U', \B')})[\xi_{(\U', \B')}] \\
		=& \lim_{t \to 0} \frac{ \substack{\tr(\V_{\B + t \xi_{B'}} \eta_{U' + t \xi_{U'} }^\top \theta_{U' + t\xi_{U'} }) + \tr(\W_{\B + t \xi_{B'}} \eta_{B' + t \xi_{B'} } \W_{\B + t \xi_{B'}} \theta_{B'+ t \xi_{B'}}) - \tr(\V_\B \eta_{U'}^\top \theta_{U'}) - \tr(\W_\B \eta_{B'} \W_\B \theta_{B'}) }}{t} \\
		= & \tr( \rmD \V_\B[ \xi_{B'} ] \eta_{U'}^\top \theta_{U'} ) + \tr( \rmD \W_{\B}[\xi_{B'}] \eta_{B'} \W_\B \theta_{B'} )   + \tr( \W_\B \eta_{B'} \rmD \W_{\B}[\xi_{B'}] \theta_{B'} )\\
		& + \tr( \V_\B \rmD \eta_{U'} [ \xi_{U'} ]^\top \theta_{U'} )  + \tr( \V_\B \eta_{U'}\rmD \theta_{U'}[ \xi_{U'} ] ) + \tr( \W_\B \rmD \eta_{B'} [\xi_{B'} ]  \W_\B \theta_{B'} ) + \tr( \W_\B \eta_{B'}  \W_\B \rmD \theta_{B'} [\xi_{B'} ]  ) \\
		= & g_{(\U, \B)}(  \rmD \eta_{(\U', \B')}[\xi_{(\U', \B')}], \theta_{(\U', \B')}  ) + g_{(\U, \B)}(  \eta_{(\U', \B')}, \rmD \theta_{(\U', \B')}[\xi_{(\U', \B')}] ) \\
		& +  \tr( \rmD \V_\B[ \xi_{B'} ] \eta_{U'}^\top \theta_{U'} ) + \tr( \rmD \W_{\B}[\xi_{B'}] \eta_{B'} \W_\B \theta_{B'} )   + \tr( \W_\B \eta_{B'} \rmD \W_{\B}[\xi_{B'}] \theta_{B'} ).
	\end{split}
\end{equation} Combining \eqref{eq: general-koszul-formula} with \eqref{eq: directional-derivative-gUB-instance}, we conclude \eqref{eq: koszul-formula-st-s+}.

{\bf (Part (iii))}.  For vector fields $\xi, \eta, \theta$ on $\st(r,p_1) \times \bbR^{p_2 \times r}_*$, we have
\begin{equation}\label{eq: directional-derivative-gUY-instance}
	\begin{split}
		& \rmD g_{(\U, \Y)}(\eta_{(\U', \Y')}, \theta_{(\U', \Y')})[\xi_{(\U', \Y')}] \\
		= &  \lim_{t \to 0} \frac{\tr(\V_{\Y + t \xi_{Y'}} \eta_{U' + t \xi_{U'} }^\top \theta_{U' + t\xi_{U'} }) + \tr(\W_{\Y + t \xi_{Y'}} \eta_{Y' + t \xi_{Y'} }^\top \theta_{Y'+ t \xi_{Y'}})  - \tr(\V_\Y \eta_{U'}^\top \theta_{U'}) - \tr(\W_\Y \eta_{Y'}^\top \theta_{Y'}) }{t}\\
		= & \tr( \rmD \V_\Y[\xi_{Y'}] \eta_{U'}^\top \theta_{U'} )   + \tr(\rmD \W_\Y[\xi_{Y'}] \eta_{Y'}^\top \theta_{Y'}) \\
		& + g_{(\U,\Y)}( \rmD\eta_{(\U', \Y')} [\xi_{(\U', \Y')} ], \theta_{(\U', \Y')}  )  + g_{(\U,\Y)}(\eta_{(\U', \Y')}, \rmD\theta_{(\U', \Y')}[\xi_{(\U', \Y')} ]  ).
	\end{split}
\end{equation} Combining \eqref{eq: general-koszul-formula} with \eqref{eq: directional-derivative-gUY-instance}, we conclude \eqref{eq: koszul-formula-st-R}. This finishes the proof. \quad $\blacksquare$.


\section{Additional Proofs in Section \ref{sec: general-strategy-for-connection}} \label{sec: additional-proof-general-strategy-in-connection}

\subsection{Proof of Lemma \ref{lm: local-minimizer-connection}.}
Since $\ell$ is a diffeomorphism between $\cM^q$ and $\cM^e$, $\ell$ and $\ell^{-1}$ are both open maps \cite[Proposition 2.15]{lee2013smooth}. If $\X$ is a minimizer of \eqref{eq: general-strategy-example-embedded} in a neighborhood of $\X$ say $\cN_1$, $\ell^{-1}$ maps $\cN_1$ to an open neighborhood of $\ell^{-1}(\X)$ say $\cN_2$, and $\ell^{-1}(\X)$ is a minimizer of \eqref{eq: general-strategy-example-quotient} for all $[\Z]\in \cN_2$. This implies $\ell^{-1}(\X)$ is a local minimizer of \eqref{eq: general-strategy-example-quotient}. This argument also applies in another direction. \quad $\blacksquare$

\subsection{Proof of Theorem \ref{th: hessian-sandwich}}
First, because $\overline{\Hess\, h([\Z])}$ and $\Hess f(\X)$ are by definition self-adjoint linear maps from $\cH_\Z \widebar{\cM}^q$ and $T_\X \cM^e$ to $\cH_\Z \widebar{\cM}^q$ and $T_\X \cM^e$, respectively, both $\Hess f(\X)$ and $\overline{\Hess\, h([\Z])}$ have $p$ eigenvalues as $\dim(T_\X \cM^e) = \dim(\cH_\Z \widebar{\cM}^q) = d$. Suppose $\u_1,\ldots,\u_d$ are eigenvectors corresponding to $\lambda_1(\Hess f(\X))$,$\ldots$, $\lambda_p(\Hess f(\X))$ and $\mathbf{v}_1,\ldots,\mathbf{v}_d$ are eigenvectors corresponding to $\lambda_1(\overline{\Hess\, h([\Z])}),\ldots,\lambda_p(\overline{\Hess\, h([\Z])})$. For $k = 1,\ldots, d$, define
\begin{equation*}
	\begin{split}
		&\mathcal{U}_k = \rspan\{\mathbf{u}_1,\ldots,\mathbf{u}_k\}, \quad \mathcal{U}'_k = \rspan\{\cL^{-1}(\mathbf{u}_1),\ldots,\cL^{-1}(\mathbf{u}_k)\},\\
		&\mathcal{V}_k = \rspan\{\mathbf{v}_1,\ldots,\mathbf{v}_k\}, \quad \mathcal{V}'_k = \rspan\{\cL(\mathbf{v}_1),\ldots,\cL(\mathbf{v}_k)\}.
	\end{split}
\end{equation*}

Let us first consider the case that $\lambda_k(\Hess f(\X)) \geq 0$. The Max-min theorem for eigenvalues (Lemma \ref{lm: max-min-theorem}) yields
\begin{equation} \label{ineq: spectrum-ineq1}
	\begin{split}
		\lambda_k(\overline{\Hess\, h([\Z])}) &\geq \min_{\mathbf{u}' \in \cU_k', \mathbf{u}' \neq 0  } \frac{\overline{\Hess\, h([\Z])}[\mathbf{u}',\mathbf{u}']}{\bar{g}_\Z(\mathbf{u}',\mathbf{u}')} \overset{ \eqref{eq: general-setting-hessian-connection} }= \min_{\mathbf{u}' \in \cU_k', \mathbf{u}' \neq 0} \frac{\Hess f(\X)[\cL(\mathbf{u}'),\cL(\mathbf{u}')]}{\bar{g}_\Z(\mathbf{u}',\mathbf{u}')} \\
		& =  \min_{\mathbf{u} \in \cU_k, \mathbf{u} \neq 0} \frac{\Hess f(\X)[\mathbf{u},\mathbf{u}]}{\bar{g}_\Z(\cL^{-1}(\mathbf{u}),\cL^{-1}(\mathbf{u}))} \geq  \min_{\mathbf{u} \in \cU_k,\mathbf{u} \neq \0 }\frac{\lambda_k(\Hess f(\X))g_\X(\mathbf{u},\mathbf{u})}{\bar{g}_\Z(\cL^{-1}(\mathbf{u}),\cL^{-1}(\mathbf{u}))} \\
		&\overset{\eqref{eq: spectrum-bound-general-setting-L} } \geq \alpha \lambda_k(\Hess f(\X)) \geq 0.
	\end{split}
\end{equation} On the other hand, we have
\begin{equation}\label{ineq: spectrum-ineq2}
	\begin{split}
		& \quad \lambda_k(\Hess f(\X)) \\
		&\overset{\text{Lemma } \ref{lm: max-min-theorem}} \geq \min_{\mathbf{v}' \in \cV_k',\mathbf{v}' \neq \0 } \frac{\Hess f(\X)[\cL \cL^{-1}(\mathbf{v}'),\cL \cL^{-1}(\mathbf{v}')]}{g_\X(\mathbf{v}',\mathbf{v}')} \overset{ \eqref{eq: general-setting-hessian-connection} } = \min_{\mathbf{v}' \in \cV_k',\mathbf{v}' \neq \0 } \frac{\overline{\Hess\, h([\Z])}[\cL^{-1}(\mathbf{v}'),\cL^{-1}(\mathbf{v}')]}{g_\X(\mathbf{v}',\mathbf{v}')}\\
		 &= \min_{\mathbf{v} \in \cV_k,\mathbf{v} \neq \0 } \frac{\overline{\Hess\, h([\Z])}[\mathbf{v},\mathbf{v}] }{g_\X(\cL(\mathbf{v}),\cL(\mathbf{v}))} \geq \min_{\mathbf{v} \in \cV_k,\mathbf{v} \neq \0 } \frac{\lambda_k(\overline{\Hess\, h([\Z])}) \bar{g}_\Z(\mathbf{v},\mathbf{v}) }{g_\X(\cL(\mathbf{v}),\cL(\mathbf{v}))} \\
		 & \overset{\eqref{ineq: spectrum-ineq1}, \eqref{eq: spectrum-bound-general-setting-L} }\geq \lambda_k(\overline{\Hess\, h([\Z])} )/\beta.
	\end{split}
\end{equation} So we have proved the result for the case that $\lambda_k(\Hess f(\X)) \geq 0$. When $\lambda_k(\Hess f(\X)) < 0$, we have $\lambda_{p+1-k}(-\Hess f(\X)) = -\lambda_k(\Hess f(\X)) > 0$. Following the same proof of \eqref{ineq: spectrum-ineq1} and \eqref{ineq: spectrum-ineq2}, we have
\begin{equation*}
	\begin{split}
		-\lambda_k(\overline{\Hess\, h([\Z])}) = \lambda_{p+1-k}(-\overline{\Hess\, h([\Z])}) \geq \alpha \lambda_{p+1-k}(-\Hess f(\X)) = -\alpha\lambda_k(\Hess f(\X)) > 0,\\
		-\lambda_k(\Hess f(\X)) =  \lambda_{p+1-k}(-\Hess f(\X)) \geq  \lambda_{p+1-k}(-\overline{\Hess\, h([\Z])})/\beta = -\lambda_k(\overline{\Hess\, h([\Z])})/\beta.
	\end{split}
\end{equation*} This finishes the proof of this theorem. \quad $\blacksquare$

\section{Additional Proofs in Section \ref{sec: embedded-quotient-fixed-rank-matrix} } \label{proof-sec: embedded-quotient-fixed-rank-matrix}

\subsection{Proof of Lemma \ref{lm: psd-quotient-manifold1-prop}.}

It has been shown in \cite[Section 4]{journee2010low} that the vertical space of $\widebar{\cM}_{r+}^{q_1}$ is $\{ \theta_\Y: \theta_\Y = \Y \bOmega, \bOmega = - \bOmega^\top \}$ with dimension $(r^2-r)/2$. By observing $\Y = \U \P$, it is easy to check
\begin{equation*}
\begin{split}
	\{ \theta_\Y: \theta_\Y = \U \bOmega' \P^{-\top}, \bOmega' = - \bOmega^{'\top} \} &=\{ \theta_\Y: \theta_\Y = \Y \P^{-1} \bOmega' \P^{-\top}, \bOmega' = - \bOmega^{'\top}  \} \\
	& = \{ \theta_\Y: \theta_\Y = \Y \bOmega, \bOmega = - \bOmega^\top  \} = \cV_\Y \widebar{\cM}_{r+}^{q_1}.
\end{split}
\end{equation*} In addition, notice $\dim(\cH_\Y \widebar{\cM}_{r+}^{q_1}) = pr-(r^2-r)/2$ and $(r^2-r)/2 + (pr-(r^2-r)/2) = pr$, then $\cH_\Y \widebar{\cM}_{r+}^{q_1}$ is the horizontal space if we can show it is orthogonal to $\cV_\Y \widebar{\cM}_{r+}^{q_1}$ with respect to $\bar{g}^{r+}_\Y$. For $\eta_\Y = \U \bOmega \P^{-\top} \in \cV_\Y \widebar{\cM}_{r+}^{q_1}$ and $\theta_\Y = (\U \S + \U_\perp \D) \P^{-\top} \in \cH_\Y \widebar{\cM}_{r+}^{q_1}$, 
\begin{equation} \label{eq: lm-psd-quotient-horizontal-eq1}
\begin{split}
	\bar{g}^{r+}_\Y(\eta_\Y, \theta_\Y) = \langle  \bOmega, \S \P^{-\top}\W_\Y \P^{-1} \rangle \overset{(a)}= - \langle  \bOmega^\top , \S \P^{-\top}\W_\Y \P^{-1} \rangle &\overset{(b)}= - \langle  \bOmega, \S \P^{-\top}\W_\Y \P^{-1} \rangle \\
	&= - \bar{g}^{r+}_\Y(\eta_\Y, \theta_\Y),
\end{split}
\end{equation} where (a) is because $\bOmega = - \bOmega^\top$ and (b) is because $\S \P^{-\top}\W_\Y \P^{-1}$ is symmetric by the construction of $\cH_\Y \widebar{\cM}_{r+}^{q_1}$. \eqref{eq: lm-psd-quotient-horizontal-eq1} implies $\bar{g}^{r+}_\Y(\eta_\Y, \theta_\Y) = 0$. This finishes the proof for the first part.

From Section \ref{sec: Riemannian-opt-quotient}, we know to show $\cM_{r+}^{q_1}$ is a Riemannian quotient manifold endowed with the Riemannian metric $g^{r+}_{[\Y]}$ induced from $\bar{g}^{r+}_\Y$, it is enough to show $\bar{g}^{r+}_\Y (\eta_\Y, \theta_{\Y})= \bar{g}^{r+}_{\Y \O}(\eta_{\Y \O}, \theta_{\Y\O})$ for any $\O \in \bbO_r$, where $\eta_\Y, \theta_\Y$ and $\eta_{\Y \O}, \theta_{\Y\O}$ are the horizontal lifts of $\eta_{[\Y]}$ and $\theta_{[\Y]}$ at $\Y$ and $\Y\O$, respectively. By Lemma \ref{lm: horizontal-lifts-connection}(i), we have the the horizontal lifts of $\eta_{[\Y]}$ at $\Y$ and $\Y \O$ are related as $\eta_{\Y\O} = \eta_{\Y} \O$. So
\begin{equation*}
	\bar{g}^{r+}_{\Y \O}(\eta_{\Y \O}, \theta_{\Y\O}) = \tr(\W_{\Y\O} \eta_{\Y \O}^\top \theta_{\Y\O} ) = \tr(\W_{\Y\O} \O^\top \eta_\Y^\top \theta_\Y \O ),
\end{equation*} and $\bar{g}^{r+}_{\Y \O}(\eta_{\Y \O}, \theta_{\Y\O}) = \bar{g}^{r+}_\Y (\eta_\Y, \theta_{\Y})$ holds for any $\eta_\Y, \theta_{\Y}$ if and only if $\W_\Y = \O \W_{\Y\O} \O^\top$ holds for any $\O \in \bbO_r$. This shows $\cM_{r+}^{q_1}$ is a Riemannian quotient manifold. \quad $\blacksquare$

\subsection{Proof of Lemma \ref{lm: psd-quotient-manifold2-prop} }
The expressions for vertical and horizontal spaces and their dimensions are provided in \cite[Theorem 1]{bonnabel2010riemannian}. Next, we prove statement (ii). From Section \ref{sec: Riemannian-opt-quotient}, we know to show $\cM_{r+}^{q_2}$ is a Riemannian quotient manifold endowed with the Riemannian metric $g^{r+}_{[\U, \B]}$ induced from $\bar{g}^{r+}_{(\U,\B)}$, it is enough to show $\bar{g}^{r+}_{(\U,\B)} (\eta_{(\U,\B)}, \theta_{{(\U,\B)}})= \bar{g}^{r+}_{(\U\O, \O^\top \B \O )}(\eta_{(\U\O, \O^\top \B \O )}, \theta_{(\U\O, \O^\top \B \O )})$ for any $\O \in \bbO_r$, where $\eta_{(\U,\B)}, \theta_{(\U,\B)}$ and $\eta_{(\U\O, \O^\top \B \O )}, \theta_{(\U\O, \O^\top \B \O )}$ are the horizontal lifts of $\eta_{[\U, \B]}$ and $\theta_{[\U, \B]}$ at $(\U,\B)$ and $(\U\O, \O^\top \B \O )$, respectively. By Lemma \ref{lm: horizontal-lifts-connection}(ii), we have the the horizontal lifts of $\theta_{[\U, \B]}$ at $(\U, \B)$ and $(\U \O, \O^\top \B \O)$ are related as $\theta_{(\U \O, \O^\top \B \O)} = [(\theta_U \O)^\top \quad (\O^\top \theta_{B} \O )^\top ]^\top$ given $\theta_{(\U,\B)} = [\theta_U^\top \quad \theta_B^\top]^\top$. So
\begin{equation*}
	\bar{g}^{r+}_{(\U\O, \O^\top \B \O) }(\eta_{(\U\O, \O^\top \B \O) }, \theta_{(\U\O, \O^\top \B \O) }) = \tr(\V_{\O^\top \B \O} \O^\top \eta_U^\top \theta_U \O  ) +  \tr(\W_{\O^\top \B\O} \O^\top \eta_{B} \O \W_{\O^\top \B\O} \O^\top  \theta_{B} \O ), 
\end{equation*} and it is equal to $\bar{g}^{r+}_{(\U, \B)} (\eta_{(\U, \B)}, \theta_{(\U, \B)})$ if and only if $\V_\B = \O \V_{\O^\top \B \O} \O^\top$ and $\W_B = \O \W_{\O^\top \B \O} \O^\top$ hold for any $\O \in \bbO_r$. This finishes the proof. \quad $\blacksquare$

\subsection{Proof of Lemma \ref{lm: completeness-quotient-PSD}.}
The results under the full-rank factorization are provided in \cite[Proposition A.1]{massart2020quotient}. Next, we prove the results under the polar factorization. If $\U_2 = \U_1 \O $ and $\B_2 = \O^\top \B_1 \O$ for some $\O \in \bbO_r$, it is easy to check $\U_1 \B_1 \U_1^\top = \U_2 \B_2 \U_2^\top$. 

Now, if $\U_1 \B_1 \U_1^\top = \U_2 \B_2 \U_2^\top$ for some $\U_1, \U_2 \in \st(r,p)$, $\B_1, \B_2 \in \bbS_+(r)$. Then the column spans of $\U_1 \B_1 \U_1^\top$ and $\U_2 \B_2 \U_2^\top$ are $\U_1$ and $\U_2$, respectively and they are the same. Hence there is $\O \in \bbR^{r\times r}$ such that $\U_2 = \U_1 \O$. Moreover, since $\U_1, \U_2 \in \st(r,p)$, $\O^\top \O = \I$, i.e., $\O \in \bbO_r$. So $\U_1 \B_1 \U_1^\top = \U_2 \B_2 \U_2^\top$ and $\U_2 = \U_1 \O$ yield $\B_2 = \O^\top \B_1 \O$.

Finally, $\cM_{r+} = \{\U \B \U^\top: \U \in \st(r,p), \B \in \bbS_+(r) \}$ follows from the facts: (1) $\U \B \U^\top \in \cM_{r+}$ for $\U \in \st(r,p), \B \in \bbS_+(r)$ and (2) any matrix $\X \in \cM_{r+}$ admits eigendecomposition which has form $\U \B \U^\top$. \quad $\blacksquare$

\subsection{Proof of Lemma \ref{lm: general-quotient-manifold1-prop}.}

The vertical space of $\widebar{\cM}_r^{q_1}$ was given in \cite[Eq. (7)]{absil2014two} and it has dimension $r^2$. Moreover, we can check $\dim(\cH_{(\L,\R)} \widebar{\cM}_{r}^{q_1}) = (p_1 + p_2-r)r$. If we can show $\cV_{(\L,\R)} \widebar{\cM}_{r}^{q_1} \perp \cH_{(\L,\R)} \widebar{\cM}_{r}^{q_1}$ with respect to $\bar{g}^{r}_{(\L,\R)}$, then it implies $\cH_{(\L,\R)} \widebar{\cM}_{r}^{q_1}$ is a valid horizontal space choice. Suppose $\theta_{(\L,\R)} = \begin{bmatrix}
	\U \S \P_2^{-\top}\\
	-\V \S^\top \P_1^{-\top}
\end{bmatrix} \in \cV_{(\L,\R)} \widebar{\cM}_{r}^{q_1}$ and $\eta_{(\L,\R)} = \begin{bmatrix}
	(\U \S' \P_2\W_{\L, \R}^{-1} \P_2^\top + \U_\perp \D_1 ) \P_2^{-\top}\\
	(\V \S^{'\top} \P_1 \V_{\L,\R}^{-1} \P_1^\top + \V_\perp \D_2  ) \P_1^{-\top}
\end{bmatrix} \in \cH_{(\L,\R)} \widebar{\cM}_{r}^{q_1}$, then by some simple calculations we have $\bar{g}_{(\L,\R)}^r ( \eta_{(\L,\R)}, \theta_{(\L,\R)} ) = \tr(\S^{'\top} \S ) - \tr(\S^{'\top} \S ) = 0$.

To show $\cM_{r}^{q_1}$ is a Riemannian quotient manifold endowed with the Riemannian metric $g^{r}_{[\L,\R]}$, it is enough to show $\bar{g}^{r}_{(\L,\R)} (\eta_{(\L,\R)}, \theta_{(\L,\R)})= \bar{g}^{r}_{(\L\M,\R\M^{-\top})}(\eta_{(\L\M,\R\M^{-\top})}, \theta_{(\L\M,\R\M^{-\top})})$ for any $\M \in \GL(r)$, where $\eta_{(\L,\R)}, \theta_{(\L,\R)}$ and $\eta_{(\L\M,\R\M^{-\top})}, \theta_{(\L\M,\R\M^{-\top})}$ are the horizontal lifts of $\eta_{[\L,\R]}$ and $\theta_{[\L,\R]}$ at $(\L,\R)$ and $(\L\M,\R\M^{-\top})$, respectively. By Lemma \ref{lm: horizontal-lifts-connection}(iii), we have the the horizontal lifts of $\theta_{[\L,\R]}$ at $(\L,\R)$ and $(\L\M,\R\M^{-\top})$ are related as $\theta_{(\L\M,\R\M^{-\top})} = [(\theta_L \M)^\top \quad (\theta_{R} \M^{-\top})^\top ]^\top$ given $\theta_{(\L,\R)} = [\theta_L^\top \quad \theta_R^\top]^\top$. So
\begin{equation*}
	\begin{split}
		&\bar{g}^{r}_{(\L\M,\R\M^{-\top})}(\eta_{(\L\M,\R\M^{-\top})}, \theta_{(\L\M,\R\M^{-\top})}) \\
		=& \tr\left(\W_{\L\M,\R\M^{-\top}} (\eta_L \M)^\top \theta_L \M \right) + \tr\left(\V_{\L\M,\R\M^{-\top}} (\eta_R \M^{-\top})^\top \theta_R \M^{-\top } \right),
	\end{split}
\end{equation*} and it is equal to $\bar{g}_{(\L,\R)}^r ( \eta_{(\L,\R)}, \theta_{(\L,\R)} )$ for any $\eta_{(\L,\R)}, \theta_{(\L,\R)}$ if and only if $\W_{\L,\R}$ and $\V_{\L,\R}$ satisfy the assumption given in the lemma. This shows $\cM_{r}^{q_1}$ is a Riemannian quotient manifold endowed with the Riemannian metric $g^{r}_{[\L,\R]}$. This finishes the proof. \quad $\blacksquare$

\subsection{Proof of Lemma \ref{lm: general-quotient-manifold2-prop}.}
The vertical space of $\cM_{r}^{q_2}$ was given in \cite[Eq. (10)]{mishra2013low}. Distinct from \cite{mishra2013low}, here the horizontal space in \eqref{eq: vertical-horizontal-quotient-general-manifold2} is chosen in a non-canonical way, i.e., it is not orthogonal to $\cV_{(\U,\B,\V)} \widebar{\cM}_{r}^{q_2}$ with respect to $\bar{g}_{(\U,\B,\V)}^r$, so we need to verify $\cH_{(\U,\B,\V)} \widebar{\cM}_{r}^{q_2}$ in \eqref{eq: vertical-horizontal-quotient-general-manifold2} is a valid horizontal space choice. To show this, we need to show $T_{(\U,\B,\V)} \widebar{\cM}_r^{q_2} = \cV_{(\U,\B,\V)} \widebar{\cM}_{r}^{q_2} \oplus \cH_{(\U,\B,\V)} \widebar{\cM}_{r}^{q_2}$, where $T_{(\U,\B,\V)} \widebar{\cM}_r^{q_2} = T_\U\st(r,p_1) \times T_\B \bbS_+(r) \times T_\V \st(r,p_2)$.

First, it is easy to check $\dim(\cV_{(\U,\B,\V)} \widebar{\cM}_{r}^{q_2}) = (r^2-r)/2$, $\dim(\cH_{(\U,\B,\V)} \widebar{\cM}_{r}^{q_2}) = (p_1 + p_2 - r)r$ and $\dim( T_{(\U,\B,\V)} \widebar{\cM}_r^{q_2} ) = (p_1 + p_2 - r)r + (r^2-r)/2$. So $\dim(\cV_{(\U,\B,\V)} \widebar{\cM}_{r}^{q_2}) + \dim(\cH_{(\U,\B,\V)} \widebar{\cM}_{r}^{q_2}) = \dim( T_{(\U,\B,\V)} \widebar{\cM}_r^{q_2} )$. Moreover, for $\eta_{(\U,\B,\V)} \in \cV_{(\U,\B,\V)} \widebar{\cM}_{r}^{q_2}$, $\theta_{(\U,\B,\V)} \in \cH_{(\U,\B,\V)} \widebar{\cM}_{r}^{q_2}$, we have $\eta_{(\U,\B,\V)} + \theta_{(\U,\B,\V)} \in T_{(\U,\B,\V)} \widebar{\cM}_r^{q_2}$. Finally, given $\xi_{(\U,\B, \V)} = [(\U \bOmega + \U_\perp \D)^\top \quad \S \quad (\V \bOmega' + \V_\perp \D')^\top ]^\top$, we can find $\eta_{(\U,\B,\V)} = [(\U \bOmega_1)^\top \quad (\B\bOmega_1 - \bOmega_1\B)^\top \quad (\V \bOmega_1)^\top  ]^\top \in \cV_{(\U,\B,\V)} \widebar{\cM}_{r}^{q_2}$, $\theta_{(\U,\B,\V)} = [ (\U_\perp \D_1 + \U \bOmega_2 )^\top \quad \theta_B \quad (\V_\perp \D_2 - \V \bOmega_2 )^\top ]^\top \in \cH_{(\U,\B,\V)} \widebar{\cM}_{r}^{q_2}$ such that $\eta_{(\U,\B,\V)} + \theta_{(\U,\B,\V)} = \xi_{(\U,\B, \V)}$, where $\bOmega_1, \bOmega_2,\theta_B,\D_1,\D_2$ are uniquely determined by the equation system $\bOmega_1 + \bOmega_2 = \bOmega, \bOmega_1 - \bOmega_2 = \bOmega',\theta_B = \S - (\B \bOmega_1 - \bOmega_1 \B), \D_1 = \D, \D_2 = \D'.$ This proves $T_{(\U,\B,\V)} \widebar{\cM}_r^{q_2} = \cV_{(\U,\B,\V)} \widebar{\cM}_{r}^{q_2} \oplus \cH_{(\U,\B,\V)} \widebar{\cM}_{r}^{q_2}$

The claim $\cM_{r}^{q_2}$ endowed with metric $g^{r}_{[\U,\B,\V]}$ induced from $\bar{g}^{r}_{(\U,\B,\V)}$ is a Riemannian quotient manifold has appeared in \cite[Section 3]{mishra2013low} without a proof, here we provide the proof for completeness. The result follows if we can show $$\bar{g}^{r}_{(\U,\B,\V)} (\eta_{(\U,\B,\V)}, \theta_{(\U,\B,\V)})= \bar{g}^{r}_{(\U\O,\O^\top\B\O,\V\O)}(\eta_{(\U\O,\O^\top\B\O,\V\O)}, \theta_{(\U\O,\O^\top\B\O,\V\O)})$$ holds for any $\O \in \bbO_r$, where $\eta_{(\U,\B,\V)}, \theta_{(\U,\B,\V)}$ and $\eta_{(\U\O,\O^\top\B\O,\V\O)}, \theta_{(\U\O,\O^\top\B\O,\V\O)}$ are the horizontal lifts of $\eta_{[\U,\B,\V]}$ and $\theta_{[\U,\B,\V]}$ at $(\U,\B,\V)$ and $(\U\O,\O^\top\B\O,\V\O)$, respectively. By Lemma \ref{lm: horizontal-lifts-connection}(iv), we have the the horizontal lifts of $\theta_{[\U,\B,\V]}$ at $(\U,\B,\V)$ and $(\U\O,\O^\top\B\O,\V\O)$ are related as $\theta_{(\U\O,\O^\top\B\O,\V\O)} = [(\theta_U \O)^\top \quad (\O^\top \theta_{B} \O)^\top \quad (\theta_V \O)^\top ]^\top$ given $\theta_{(\U,\B,\V)} = [\theta_U^\top \quad \theta_B^\top \quad \theta_V^\top]^\top$. So
\begin{equation*}
	\begin{split}
		&\bar{g}^{r}_{(\U\O,\O^\top\B\O,\V\O)}(\eta_{(\U\O,\O^\top\B\O,\V\O)}, \theta_{(\U\O,\O^\top\B\O,\V\O)})\\
		 = & \tr(\O^\top \eta_U^\top \theta_U \O) + \tr\left((\O^\top\B\O)^{-1} \O^\top \eta_B \O (\O^\top\B\O)^{-1} \O^\top \theta_B \O   \right) + \tr(\O^\top \eta_V^\top \theta_V \O)\\
		 = & \tr(\eta_U^\top \theta_U ) + \tr(\B^{-1} \eta_B \B^{-1} \theta_B ) + \tr(\eta_V^\top \theta_V) = \bar{g}_{(\U,\B,\V)}^r ( \eta_{(\U,\B,\V)}, \theta_{(\U,\B,\V)} ).
	\end{split}
\end{equation*} This shows $\cM_{r}^{q_2}$ endowed with the metric $g^{r}_{[\U,\B,\V]}$ is a Riemannian quotient manifold. This finishes the proof of this lemma. \quad $\blacksquare$

\subsection{Proof of Lemma \ref{lm: general-quotient-manifold3-prop}.}
The vertical space of $\cM_r^{q_3}$ was given in \cite[Eq. (33)]{absil2014two}. Distinct from \cite{absil2014two}, here the horizontal space in \eqref{eq: vertical-horizontal-quotient-general-manifold3} is chosen in a non-canonical way, so we need to verify $\cH_{(\U,\Y)} \widebar{\cM}_{r}^{q_3}$ in \eqref{eq: vertical-horizontal-quotient-general-manifold3} is a valid horizontal space choice. To show this, we need to show $T_{(\U,\Y)} \widebar{\cM}_r^{q_3} = \cV_{(\U,\Y)} \widebar{\cM}_{r}^{q_3} \oplus \cH_{(\U,\Y)} \widebar{\cM}_{r}^{q_3}$, where $T_{(\U,\Y)} \widebar{\cM}_r^{q_3} = T_\U\st(r,p_1) \times \bbR^{p_2 \times r}_*$.

First, it is easy to check $\dim(\cV_{(\U,\Y)} \widebar{\cM}_{r}^{q_3}) = (r^2-r)/2$, $\dim(\cH_{(\U,\Y)} \widebar{\cM}_{r}^{q_3}) = (p_1 + p_2 - r)r$ and $\dim( T_{(\U,\Y)} \widebar{\cM}_r^{q_3} ) = (p_1 + p_2 - r)r + (r^2-r)/2$. So $\dim(\cV_{(\U,\Y)} \widebar{\cM}_{r}^{q_3}) + \dim(\cH_{(\U,\Y)} \widebar{\cM}_{r}^{q_3}) = \dim( T_{(\U,\Y)} \widebar{\cM}_r^{q_3} )$. Moreover, for $\eta_{(\U,\Y)} \in \cV_{(\U,\Y)} \widebar{\cM}_{r}^{q_3}$, $\theta_{(\U,\Y)} \in \cH_{(\U,\Y)} \widebar{\cM}_{r}^{q_3}$, we have $\eta_{(\U,\Y)} + \theta_{(\U,\Y)} \in T_{(\U,\Y)} \widebar{\cM}_r^{q_3}$. Finally, given $\xi_{(\U,\Y)} = [(\U \bOmega + \U_\perp \D)^\top \quad \xi_Y^\top ]^\top$, we can find $\eta_{(\U,\Y)} = [(\U \bOmega_1)^\top \quad (\Y \bOmega_1)^\top  ]^\top \in \cV_{(\U,\Y)} \widebar{\cM}_{r}^{q_3}$, $\theta_{(\U,\Y)} = [ (\U_\perp \D_1)^\top \quad \theta_Y^\top ]^\top \in \cH_{(\U,\Y)} \widebar{\cM}_{r}^{q_3}$ such that $\eta_{(\U,\Y)} + \theta_{(\U,\Y)} = \xi_{(\U,\Y)}$, where $\bOmega_1, \D_1, \theta_Y$ are uniquely determined as by the equation system $\bOmega_1 = \bOmega, \D_1 = \D, \theta_Y = \xi_Y - \Y \bOmega_1.$ This proves $T_{(\U,\Y)} \widebar{\cM}_r^{q_3} = \cV_{(\U,\Y)} \widebar{\cM}_{r}^{q_3} \oplus \cH_{(\U,\Y)} \widebar{\cM}_{r}^{q_3}$.

To show $\cM_{r}^{q_3}$ is a Riemannian quotient manifold, we only need to show $\bar{g}^{r}_{(\U, \Y)} (\eta_{(\U, \Y)}, \theta_{(\U,\Y)})= \bar{g}^{r}_{(\U\O,\Y\O)}(\eta_{(\U\O,\Y\O)}, \theta_{(\U\O,\Y\O)})$ for any $\O \in \bbO_r$, where $\eta_{(\U, \Y)}, \theta_{(\U, \Y)}$ and $\eta_{(\U\O,\Y\O)}, \theta_{(\U\O,\Y\O)}$ are the horizontal lifts of $\eta_{[\U,\Y]}$ and $\theta_{[\U,\Y]}$ at $(\U,\Y)$ and $(\U\O,\Y\O)$, respectively. By Lemma \ref{lm: horizontal-lifts-connection}(v), we have the the horizontal lifts of $\eta_{[\U,\Y]}$ at $(\U,\Y)$ and $(\U\O,\Y\O)$ are related as $\eta_{(\U\O,\Y\O)} = [(\eta_U \O)^\top \quad (\eta_Y \O)^\top]^\top$ given $\eta_{(\U, \Y)} = [\eta_U^\top \quad \eta_Y^\top]^\top$. So
\begin{equation*}
	\begin{split}
		\bar{g}^{r}_{(\U\O,\Y\O)}(\eta_{(\U\O,\Y\O)}, \theta_{(\U\O,\Y\O)}) = \tr( \V_{\Y \O} \O^\top \eta_U^\top \theta_U \O) + \tr( \W_{\Y \O} \O^\top \eta_Y^\top \theta_Y \O ),
	\end{split}
\end{equation*} and it is equal to $\bar{g}^{r}_{(\U, \Y)} (\eta_{(\U, \Y)}, \theta_{(\U,\Y)})$ for any $\eta_{(\U, \Y)}, \theta_{(\U,\Y)}$ if and only if $\V_\Y = \O \V_{\Y\O} \O^\top$ and $\W_\Y = \O \W_{\Y\O} \O^\top$ holds for any $\O \in \bbO_r$. This shows $\cM_{r}^{q_3}$ is a Riemannian quotient manifold endowed with metric $g_{[\U,\Y]}^r$ induced from $\bar{g}_{(\U,\Y)}^r$. This finishes the proof of this lemma. \quad $\blacksquare$

\subsection{Proof of Lemma \ref{lm: completeness-quotient-general}.}
{\bf (Full-rank factorization)} If $\L_2 = \L_1 \M$ and $\R_2 = \R_1 \M^{-\top} $ for some $\M \in \GL(r)$, it is easy to check $\L_1 \R_1^\top = \L_2 \R_2^\top$. 

Now, if $\L_1 \R_1^\top = \L_2 \R_2^\top$ for some $\L_1, \L_2 \in \bbR_*^{p_1 \times r}$ and $\R_1, \R_2 \in \bbR_*^{p_2 \times r}$. Then the column spans of $\L_1 \R_1^\top$ and $\L_2 \R_2^\top$ are $\rspan(\L_1)$ and $\rspan(\L_2)$, respectively and they are the same. Hence there exists $\M \in \GL(r)$ such that $\L_2 = \L_1 \M$. Thus, $\L_1 \R_1^\top = \L_2 \R_2^\top$ yields $\R_2 = \R_1 \M^{-\top}$.

Moreover, $\cM_r = \{\L \R^\top: \L \in \bbR^{p_1 \times r}_*, \R \in \bbR_*^{p_2 \times r} \}$ follows from the facts: (1) $\L \R^\top \in \cM_r$ for any $\L \in \bbR^{p_1 \times r}_*, \R \in \bbR_*^{p_2 \times r}$ and (2) any matrix $\X \in \cM_r$ admits SVD $\U \bSigma \V^\top$ and by taking $\L = \U \bSigma^{1/2}$ and $\R = \V \bSigma^{1/2} $, we have $\X = \L \R^\top$.

{\bf (Polar factorization)} If $\U_2 = \U_1 \O $, $\B_2 = \O^\top \B_1 \O$ and $\V_2 = \V_1 \O$ for some $\O \in \bbO_r$, it is easy to check $\U_1 \B_1 \V_1^\top = \U_2 \B_2 \V_2^\top$. 

Now, if $\U_1 \B_1 \V_1^\top = \U_2 \B_2 \V_2^\top$ for some $\U_1, \U_2 \in \st(r,p_1)$, $\B_1, \B_2 \in \bbS_+(r)$ and $\V_1, \V_2 \in \st(r,p_2)$. Then the column spans of $\U_1 \B_1 \V_1^\top$ and $\U_2 \B_2 \V_2^\top$ are $\rspan(\U_1)$ and $\rspan(\U_2)$, respectively and they are the same. Hence there exists $\O_1 \in \bbR^{r\times r}$ such that $\U_2 = \U_1 \O_1$. Moreover, since $\U_1, \U_2 \in \st(r,p)$, $\O_1^\top \O_1 = \I$, i.e., $\O_1 \in \bbO_r$. Similarly, we can conclude $\V_2 = \V_1 \O_2$ for some $\O_2 \in \bbO_r$. Thus, $\U_1 \B_1 \V_1^\top = \U_2 \B_2 \V_2^\top$ yields $\B_2 = \O_1^\top \B_1 \O_2$. Since $\B_2 \in \bbS_+(r)$, by the characterization of PSD factorization in Lemma \ref{lm: completeness-quotient-PSD}, we know $\O_1 = \O_2$.

Moreover, $\cM_r =\{\U \B \V^\top: \U \in \st(r,p_1), \B \in \bbS_+(r), \V \in \st(r, p_2) \}$ follows from the facts: (1) $\U \B \V^\top \in \cM_r$ for any $\U \in \st(r,p_1), \B \in \bbS_+(r), \V \in \st(r, p_2)$ and (2) any matrix $\X \in \cM_r$ admits SVD which has form $\U \B \V^\top$.

{\bf (Subspace-projection factorization)} If $\U_2 = \U_1 \O $ and $\Y_2 = \Y_1 \O$ for some $\O \in \bbO_r$, it is easy to check $\U_1 \Y_1^\top = \U_2 \Y_2^\top$. 

Now, if $\U_1 \Y_1^\top = \U_2 \Y_2^\top$ for some $\U_1, \U_2 \in \st(r,p_1)$ and $\Y_1, \Y_2 \in \bbR_*^{p_2 \times r}$. Then the column spans of $\U_1 \Y_1^\top$ and $\U_2 \Y_2^\top$ are $\rspan(\U_1)$ and $\rspan(\U_2)$, respectively and they are the same. Following the same argument as in the Polar factorization, we know $\U_2 = \U_1 \O$ for some $\O \in \bbO_r$. Thus, $\U_1 \Y_1^\top = \U_2 \Y_2^\top$ yields $\Y_2 = \Y_1 \O$.

Moreover, $\cM_r = \{\U \Y^\top: \U \in \st(r,p_1), \Y \in \bbR_*^{p_2 \times r} \}$ follows from the facts: (1) $\U \Y^\top \in \cM_r$ for any $ \U \in \st(r,p_1), \Y \in \bbR_*^{p_2 \times r} $ and (2) any matrix $\X \in \cM_{r}$ admits SVD $\U \bSigma \V^\top$ and by taking $\Y = \V \bSigma$, we have $\X = \U \Y^\top$. \quad $\blacksquare$

\section{Additional Proofs in Section \ref{sec: connection-PSD} } \label{sec: addition-proof-psd}

\subsection{Proof of Proposition \ref{prop: gradient-hessian-exp-PSD}.}
The proof is divided into four steps: in Step 1, we derive the expressions for Riemannian gradients; in Step 2, we derive the Riemannian Hessian of \eqref{eq: PSD-manifold-formulation} under the embedded geometry; in Step 3, we derive the Riemannian Hessian for $h_{r+}([\Y])$; in Step 4, we derive the Riemannian Hessian for $h_{r+}([\U,\B])$.

{\bf Step 1}. First, the Riemannian gradient expression of \eqref{eq: PSD-manifold-formulation} under the embedded geometry can be found in \cite[Proposition 1]{luo2021nonconvex}. Next, we compute the Riemannian gradients under the quotient geometry. By definition
\begin{equation} \label{eq: gradient-comp-psd-1}
	\begin{split}
		 \langle \grad \, \bar{h}_{r+}(\Y) \W_\Y, \eta_\Y \rangle=\bar{g}_\Y^{r+}(\grad \, \bar{h}_{r+}(\Y), \eta_\Y)  &= \rmD \bar{h}_{r+}(\Y)[\eta_\Y] = \langle \nabla \bar{h}_{r+}(\Y), \eta_\Y \rangle \\
		& = 2\langle \nabla f(\Y \Y^\top) \Y, \eta_\Y \rangle, \quad \forall \eta_\Y \in T_\Y \widebar{\cM}_{r+}^{q_1},
	\end{split}
\end{equation} where $\nabla \bar{h}_{r+}(\Y)$ denotes the Euclidean gradient of $\bar{h}_{r+}$ at $\Y$, and
\begin{equation} \label{eq: gradient-comp-psd-2}
	\begin{split}
		&\langle  \grad_\U\, \bar{h}_{r+}(\U, \B)\V_\B, \eta_U \rangle + \langle \W_\B \grad_\B\, \bar{h}_{r+}(\U, \B) \W_\B, \eta_B \rangle \\
		=& \bar{g}^{r+}_{(\U, \B)} ( \grad\, \bar{h}_{r+}(\U, \B), \eta_{(\U, \B)} ) \\
		=& \rmD \bar{h}_{r+}(\U, \B)[\eta_{(\U, \B)}] = \langle \nabla_\U \bar{h}_{r+}(\U,\B) , \eta_U \rangle + \langle \nabla_\B \bar{h}_{r+}(\U,\B) , \eta_B \rangle \\
		=& 2\langle \nabla f(\U \B \U^\top)\U \B, \eta_U \rangle + \langle \U^\top \nabla f(\U \B \U^\top) \U, \eta_B \rangle, \quad \forall \eta_{(\U, \B)} \in T_{(\U, \B)} \widebar{\cM}_{r+}^{q_2}.
	\end{split}
\end{equation} From \eqref{eq: gradient-comp-psd-1}, we have the following expression for $\grad \, \bar{h}_{r+}(\Y)$:
\begin{equation} \label{eq: Rgradient-total-space-psd-quotient}
	\begin{split}
		\grad \, \bar{h}_{r+}(\Y) = \nabla \bar{h}_{r+}(\Y) \W_\Y^{-1} = 2\nabla f(\Y \Y^\top) \Y \W_\Y^{-1}. 
	\end{split}
\end{equation} Since $\cV_\Y \widebar{\cM}_{r+}^{q_1} \perp \cH_\Y \widebar{\cM}_{r+}^{q_1}$, by Lemma \ref{lm: quotient-ortho-hori-gradient-express}, we have $\overline{\grad \, h_{r+}([\Y])} = \grad \, \bar{h}_{r+}(\Y)$. 

At the same time, from \eqref{eq: gradient-comp-psd-2}, we have $\grad_\B \, \bar{h}_{r+}(\U, \B) =\W^{-1}_\B \U^\top \nabla f(\U \B \U^\top) \U \W^{-1}_\B $. In addition, suppose $\grad_\U \, \bar{h}_{r+}(\U, \B) = \U_\perp \D + \U \bOmega \in T_\U \st(r,p)$, then by \eqref{eq: gradient-comp-psd-2}, we have $P_{ T_\U \st(r,p) } ( \U_\perp \D \V_\B + \U \bOmega \V_\B - 2 \nabla f(\U \B \U^\top)\U \B ) = 0$. Thus $\D = 2 \U_\perp^\top \nabla f(\U \B \U^\top) \U \B \V_\B^{-1} $ and $\bOmega$ is uniquely determined by the following Sylvester equation $\V_\B \bOmega^\top + \bOmega^\top \V_\B = 2( \B \U^\top \nabla f(\U \B \U^\top) \U - \U^\top \nabla f( \U \B \U^\top) \U \B  ) $ . So
\begin{equation*}
	\overline{\grad \, h_{r+}([\U, \B])} = P_{(\U, \B)}^{\cH} (\grad \, \bar{h}_{r+}(\U, \B)) = \begin{bmatrix}
				2 P_{\U_\perp} \nabla f(\U \B \U^\top) \U \B \V_\B^{-1} \\
				\W_\B^{-1} \U^\top \nabla f(\U \B \U^\top) \U \W_\B^{-1}
			\end{bmatrix},
\end{equation*} here $P_{(\U, \B)}^{\cH}(\cdot)$ is the projection operator onto $\cH_{(\U,\B)} \widebar{\cM}^{q_1}_{r+}$ and it satisfies $P_{(\U, \B)}^{\cH}(\eta_{(\U,\B)}) = [(P_{\U_\perp} \eta_U)^\top \quad \eta_B^\top ]^\top$ for any $\eta_{(\U,\B)} \in T_{(\U,\B)} \widebar{\cM}_{r+}^{q_2}$.

{\bf Step 2.}  Suppose $\X$ has eigendecomposition $\U' \bSigma' \U^{'\top}$ with $\U = \U' \O$ and $\O \in \bbO_r$. Then for
\begin{equation*}
	\begin{split}
		\xi_\X &= [\U \quad \U_\perp] \begin{bmatrix}
			\S & \D^\top\\
			\D & \0
		\end{bmatrix} [\U \quad \U_\perp]^\top = [\U' \quad \U_\perp] \begin{bmatrix}
			\O\S \O^\top & \O\D^\top\\
			\D\O^\top & \0
		\end{bmatrix} [\U' \quad \U_\perp]^\top,
	\end{split}
\end{equation*} by \cite[Proposition 2]{luo2021nonconvex} we have in the PSD case
\begin{equation*}
	\begin{split}
		\Hess f(\X)[\xi_\X, \xi_\X] &= \nabla^2 f(\X)[\xi_\X, \xi_\X] + 2\langle \nabla f(\X), \U_\perp \D \O^\top \bSigma^{'-1} \O \D^\top \U_\perp^\top \rangle \\
		& = \nabla^2 f(\X)[\xi_\X, \xi_\X] + 2\langle \nabla f(\X), \U_\perp \D \bSigma^{-1} \D^\top \U_\perp^\top \rangle
	\end{split}
\end{equation*} where $\bSigma = \O^\top \bSigma^{'} \O = \U^\top \X \U$.

{\bf Step 3.} We first compute
\begin{equation} \label{eq: directional-derivative-quotient-gradient-psd-1}
	\begin{split}
		\rmD \overline{ \grad \, h_{r+}([\Y])}[\theta_\Y] = \rmD \grad \, \bar{h}_{r+}(\Y)[\theta_\Y] & \overset{ \eqref{eq: Rgradient-total-space-psd-quotient} } =  \lim_{t \to 0} \left( \nabla \bar{h}_{r+} (\Y + t \theta_\Y) \W_{\Y + t\theta_\Y}^{-1} - \nabla \bar{h}_{r+}(\Y) \W_\Y^{-1}  \right)/t\\
		& = \nabla \bar{h}_{r+}(\Y) \rmD \W_\Y^{-1} [\theta_\Y] + \nabla^2 \bar{h}_{r+}(\Y)[\theta_\Y] \W_\Y^{-1}.
	\end{split}
\end{equation}

 By the definition of the bilinear form of the Riemannian Hessian, we have
\begin{equation} \label{eq: quotient-hessian-psd-manifold1}
	\begin{split}
		\overline{\Hess \, h_{r+}([\Y])}[\theta_\Y, \theta_\Y] &= \bar{g}_\Y^{r+} (\overline{\Hess \, h_{r+}([\Y])[\theta_{[\Y]}]}, \theta_\Y )\\
		&\overset{ \eqref{eq: quotient-hessian-linear-from} } = \bar{g}_\Y^{r+} \left( P_{\Y}^{\cH} \left( \widebar{\nabla}_{\theta_\Y} \overline{ \grad \, h_{r+}([\Y]) } \right) , \theta_\Y \right) \\
		& \overset{ \text{Lemma } \ref{lm: quotient-hessian-riemannian-connection-equal} } = \bar{g}_\Y^{r+} \left( \widebar{\nabla}_{\theta_\Y} \overline{ \grad \, h_{r+}([\Y])} , \theta_\Y \right)\\
		& \overset{ \eqref{eq: koszul-formula-simple-manifolds} } = \bar{g}_\Y^{r+} \left( \rmD \overline{ \grad \, h_{r+}([\Y])}[\theta_\Y] , \theta_\Y \right) + \tr\left( \rmD \W_\Y\left[\overline{ \grad \, h_{r+}([\Y])}\right] \theta_\Y^\top \theta_\Y  \right)/2\\
		& \overset{ \eqref{eq: Rgradient-total-space-psd-quotient}, \eqref{eq: directional-derivative-quotient-gradient-psd-1} } = \nabla^2 \bar{h}_{r+}(\Y)[\theta_\Y, \theta_\Y] + 2\tr\left(\W_\Y \left( \nabla f(\Y \Y^\top) \Y  \rmD \W_\Y^{-1} [\theta_\Y]\right)^\top \theta_\Y \right) \\
		& \quad + \tr\left( \rmD \W_\Y\left[\overline{ \grad \, h_{r+}([\Y])}\right] \theta_\Y^\top \theta_\Y  \right)/2.
	\end{split}
\end{equation} Finally, we have $\nabla^2 \bar{h}_{r+}(\Y)[\theta_\Y, \theta_\Y] = \nabla^2 f(\Y \Y^\top)[\Y\theta_\Y^\top + \theta_\Y \Y^\top , \Y\theta_\Y^\top + \theta_\Y \Y^\top ] + 2\langle \nabla f(\Y \Y^\top ), \theta_\Y \theta_\Y^\top \rangle$ from \cite[Proposition 2]{luo2021nonconvex} and plugging it into \eqref{eq: quotient-hessian-psd-manifold1}, we get the expression for the Riemannian Hessian of $h_{r+}([\Y])$.

{\bf Step 4.} 
For convenience, we denote $\bar{g}_{\U}^{r+}(\eta_U, \theta_U) = \tr(\V_\B \eta_U^\top \theta_U) $ as the $\U$ part inner product of $\bar{g}^{r+}_{(\U, \B)}$ and $\bar{g}_{\B}^{r+}(\eta_B, \theta_B) = \tr(\W_\B \eta_B \W_\B \theta_B)$ as the $\B$ part inner product of $\bar{g}^{r+}_{(\U, \B)}$. We have
\begin{equation} \label{eq: quotient-hessian-manifold2}
	\begin{split}
		& \quad \overline{\Hess \, h_{r+}([\U,\B])}[\theta_{(\U,\B)}, \theta_{(\U,\B)}] \\
		&= \bar{g}_{(\U, \B)}^{r+}\left( \overline{\Hess \, h_{r+}([\U,\B])[\theta_{[\U,\B]}]} , \theta_{(\U,\B)} \right)\\
		&\overset{ \eqref{eq: quotient-hessian-linear-from} } = \bar{g}_{(\U, \B)}^{r+} \left( P_{(\U, \B)}^{\cH} \left( \widebar{\nabla}_{\theta_{(\U, \B)}} \overline{ \grad \, h_{r+}([\U, \B]) } \right) , \theta_{(\U, \B)} \right)\\
		& \overset{ \text{Lemma } \ref{lm: quotient-hessian-riemannian-connection-equal} } = \bar{g}_{(\U, \B)}^{r+} \left( \widebar{\nabla}_{\theta_{(\U, \B)}} \overline{ \grad \, h_{r+}([\U, \B])} , \theta_{(\U, \B)} \right)\\
		& \overset{ \eqref{eq: koszul-formula-st-s+} } = \bar{g}_{\U}^{r+} \left( \rmD \overline{ \grad_\U \, h_{r+}([\U, \B])}[\theta_{(\U, \B)}] , \theta_{U} \right) + \bar{g}_{\B}^{r+} \left( \rmD \overline{ \grad_\B \, h_{r+}([\U, \B])}[\theta_{(\U, \B)}], \theta_B \right)\\
		&\quad  + \tr\left( \rmD \V_\B [ \overline{ \grad_\B \, h_{r+}([\U, \B])} ] \theta_U^\top \theta_U \right)/2\\
		&\quad   + \tr( \sym( \W_\B \theta_B \rmD \W_\B [\overline{ \grad_\B \, h_{r+}([\U, \B])}  ] ) \theta_B   ).
	\end{split}
\end{equation} Next, we compute $\rmD \overline{ \grad_\U \, h_{r+}([\U, \B])}[\theta_{(\U, \B)}]$ and $\rmD \overline{ \grad_\B \, h_{r+}([\U, \B])}[\theta_{(\U, \B)}]$ separately in \eqref{eq: quotient-hessian-manifold2}.
\begin{equation} \label{eq: quotient-hessian-manifold2-comp1}
\begin{split}
	&\rmD \overline{ \grad_\U \, h_{r+}([\U, \B])}[\theta_{(\U, \B)}]\\
	= & 2 \lim_{t \to 0} \Big( (\I_p - P_{\U + t \theta_U}) \nabla f\left( (\U + t\theta_U ) (\B+t\theta_B) (\U + t\theta_U)^\top\right) (\U + t \theta_U) (\B + t\theta_B)\V^{-1}_{\B+ t \theta_B} \\
	& \quad \quad - P_{\U_\perp} \nabla f(\U \B \U^\top) \U \B \V^{-1}_{\B} \Big)/t\\
	=&  -4 \sym(\U \theta_U^\top) \nabla f(\U \B \U^\top) \U \B\V^{-1}_{\B} + 2 P_{\U_\perp}\left( \nabla f(\U \B \U^\top) \theta_U \B\V^{-1}_{\B} + \nabla f(\U \B \U^\top) \U \theta_B \V^{-1}_{\B} \right) \\
	& + 2 P_{\U_\perp} \nabla f(\U \B \U^\top) \U \B \rmD \V_{\B}^{-1}[\theta_B]  \\
	& +   2 \lim_{t \to 0} P_{\U_\perp}\left(  \nabla f\left( (\U + t\theta_U ) (\B+t\theta_B) (\U + t\theta_U)^\top\right) - \nabla f(\U \B \U^\top) \right) \U \B\V^{-1}_{\B}/t\\
	=& -4 \sym(\U \theta_U^\top) \nabla f(\U \B \U^\top) \U \B\V^{-1}_{\B} + 2 P_{\U_\perp}\left( \nabla f(\U \B \U^\top) \theta_U \B\V^{-1}_{\B} + \nabla f(\U \B \U^\top) \U \theta_B \V^{-1}_{\B} \right) \\
	& + 2 P_{\U_\perp} \nabla f(\U \B \U^\top) \U \B \rmD \V_{\B}^{-1}[\theta_B]  \\
	&+ 2 P_{\U_\perp} \nabla^2 f(\U\B\U^\top)[ \U \B \theta_U^\top+ \U \theta_B \U^\top  + \theta_U \B \U^\top] \U \B\V^{-1}_{\B},
\end{split}
\end{equation} and 
\begin{equation} \label{eq: quotient-hessian-manifold2-comp2}
	\begin{split}
		&\rmD \overline{ \grad_\B \, h_{r+}([\U, \B])}[\theta_{(\U, \B)}]\\
		=& \lim_{t \to 0} \big( \W^{-1}_{\B + t \theta_B} (\U + t\theta_U)^\top \nabla f\left((\U + t\theta_U) (\B + t\theta_B) (\U + t\theta_U)^\top\right) (\U + t\theta_U) \W^{-1}_{\B + t \theta_B} \\
		& - \W^{-1}_{\B} \U^\top \nabla f(\U \B \U^\top) \U \W^{-1}_{\B} \big)/t\\
		=& (\rmD \W_\B^{-1}[\theta_B] \U^\top + \W^{-1}_{\B}  \theta_U^\top) \nabla f(\U \B \U^\top) \U \W^{-1}_{\B}   + \W^{-1}_{\B} \U^\top \nabla f(\U \B \U^\top) (\U\rmD \W_\B^{-1}[\theta_B]  + \theta_U \W^{-1}_{\B} ) \\
		& + \W^{-1}_{\B}  \U^\top  \nabla^2 f(\U\B\U^\top)[ \U \B \theta_U^\top+ \U \theta_B \U^\top  + \theta_U \B \U^\top] \U \W^{-1}_{\B} .
	\end{split}
\end{equation}
By observing $f(\X)$ is symmetric in $\X$, $\theta_{(\U, \B)} \in \cH_{(\U, \B)} \widebar{\cM}_{r+}^{q_2}$ and plugging \eqref{eq: quotient-hessian-manifold2-comp1} and \eqref{eq: quotient-hessian-manifold2-comp2} into \eqref{eq: quotient-hessian-manifold2}, we have
\begin{equation*}
	\begin{split}
		 &\overline{\Hess \, h_{r+}([\U,\B])}[\theta_{(\U,\B)}, \theta_{(\U,\B)}] \\
		=& \nabla^2 f(\U\B\U^\top)[ \U \B \theta_U^\top+ \U \theta_B \U^\top  + \theta_U \B \U^\top,  \U \B \theta_U^\top + \theta_U \B \U^\top] + 2 \langle  \nabla f(\U \B \U^\top), \theta_U \B \theta_U^\top \rangle \\
		& - 4\langle \sym(\U \theta_U^\top) \nabla f(\U \B \U^\top) \U \B, \theta_U \rangle + 2 \langle \nabla f(\U \B \U^\top) \U \theta_B , \theta_U \rangle \\
		& + 2 \langle \nabla f(\U \B \U^\top) \U \B \rmD \V_\B^{-1} [\theta_B] \V_\B, \theta_U \rangle \\
		& + \nabla^2 f(\U\B\U^\top)[ \U \B \theta_U^\top+ \U \theta_B \U^\top  + \theta_U \B \U^\top,  \U \theta_B \U^\top]\\
		&  + \tr(\W_{\B} (\rmD \W_\B^{-1} [\theta_B] \U^\top + \W^{-1}_{\B} \theta_U^\top) \nabla f(\U \B \U^\top) \U  \theta_B ) \\
		& + \tr(\U^\top \nabla f(\U \B \U^\top) (\U \rmD \W_{\B}^{-1} [\theta_B]  + \theta_U \W^{-1}_{\B})  \W_{\B} \theta_B )\\
		& + \tr\left( \rmD \V_\B [ \overline{ \grad_\B \, h_{r+}([\U, \B])} ] \theta_U^\top \theta_U \right)/2 \\
		& + \tr( \sym( \W_\B \theta_B \rmD \W_\B [\overline{ \grad_\B \, h_{r+}([\U, \B])}  ] ) \theta_B   ) \\
		=& \nabla^2 f(\U \B \U^\top)[\U \B \theta_U^\top + \U \theta_B \U^\top + \theta_U \B \U^\top, \U \B \theta_U^\top + \U \theta_B \U^\top + \theta_U \B \U^\top]  + 2 \langle  \nabla f(\U \B \U^\top), \theta_U \B \theta_U^\top \rangle\\
			& \,+2\langle \nabla f(\U \B \U^\top) \U, 2 \theta_U \theta_B + \U \rmD \W_\B^{-1}[ \theta_B] \W_\B \theta_B + \theta_U \V_\B \rmD \V_\B^{-1}[\theta_B] \B - \theta_U \U^\top \theta_U \B - \U \theta_U^\top \theta_U \B \rangle \\
			& \, + \tr( \rmD \V_\B [ \overline{ \grad_\B \, h_{r+}([\U, \B])}  ] \theta_U^\top \theta_U   )/2 + \tr( \sym(  \W_\B \theta_B \rmD \W_\B [  \overline{ \grad_\B \, h_{r+}([\U, \B])}  ]  ) \theta_B ). \quad \blacksquare
	\end{split}
\end{equation*}

\subsection{Proof of Theorem \ref{th: embedded-quotient-connection-PSD2}.}
First, \eqref{eq: gradient-connect-PSD2} is by direct calculation from the gradient expressions in Proposition \ref{prop: gradient-hessian-exp-PSD}. Next, we prove \eqref{eq: Hessian-connection-PSD2}. Since $[\U,\B]$ is a Riemannian FOSP of \eqref{eq: PSD-opt-problem-quotient-sub2}, we have
\begin{equation} \label{eq: FOSP-condition-PSD2}
	\overline{ \grad \, h_{r+}([\U,\B])} = \0 \quad \text{ and } \quad \nabla f(\U \B \U^\top) \U = \0 
\end{equation} and have $\grad f(\X) = \0$ by \eqref{eq: gradient-connect-PSD2}. So $ \nabla f(\X) = P_{\U_\perp} \nabla f(\X) P_{\U_\perp}$. Given any $\theta_{(\U,\B)} = [\theta_U^\top \quad \theta_B^\top]^\top \in \cH_{(\U,\B)} \widebar{\cM}_{r+}^{q_2}$, we have
\begin{equation} \label{eq: Hessian-con-gradient-2}
	\langle \nabla f(\X), P_{\U_\perp} \theta_U \B \B^{-1} \B \theta_U^\top  P_{\U_\perp}    \rangle = \langle \nabla f(\X), \theta_U \B \theta_U^\top \rangle,
\end{equation} where the equality is because $ \nabla f(\X) = P_{\U_\perp} \nabla f(\X) P_{\U_\perp} $.

Then by Proposition \ref{prop: gradient-hessian-exp-PSD}:
\begin{equation*}
\begin{split}
	& \quad \overline{\Hess \, h_{r+}([\U,\B])}[\theta_{(\U,\B)}, \theta_{(\U,\B)}] \\
&= \nabla^2 f(\U \B \U^\top)[\U \B \theta_U^\top + \U \theta_B \U^\top + \theta_U \B \U^\top, \U \B \theta_U^\top + \U \theta_B \U^\top + \theta_U \B \U^\top]  + 2 \langle  \nabla f(\U \B \U^\top), \theta_U \B \theta_U^\top \rangle\\
			& \,+2\langle \nabla f(\U \B \U^\top) \U, 2 \theta_U \theta_B + \U \rmD \W_\B^{-1}[ \theta_B] \W_\B \theta_B + \theta_U \V_\B \rmD \V_\B^{-1}[\theta_B] \B - \theta_U \U^\top \theta_U \B - \U \theta_U^\top \theta_U \B \rangle \\
			& \, + \tr( \rmD \V_\B [ \overline{ \grad_\B \, h_{r+}([\U, \B])}  ] \theta_U^\top \theta_U   )/2 + \tr( \sym(  \W_\B \theta_B \rmD \W_\B [  \overline{ \grad_\B \, h_{r+}([\U, \B])}  ]  ) \theta_B )\\
	& \overset{ \eqref{eq: FOSP-condition-PSD2} }= \nabla^2 f(\U \B \U^\top)[\U \B \theta_U^\top + \U \theta_B \U^\top + \theta_U \B \U^\top, \U \B \theta_U^\top + \U \theta_B \U^\top + \theta_U \B \U^\top]  \\
	& \quad  + 2 \langle  \nabla f(\U \B \U^\top), \theta_U \B \theta_U^\top \rangle\\
	& \overset{ \text{Proposition } \ref{prop: psd-bijection2}, \eqref{eq: Hessian-con-gradient-2} }= \nabla^2 f(\X)[\cL_{\U,\B}^{r+}(\theta_{(\U,\B)}),\cL_{\U,\B}^{r+}(\theta_{(\U,\B)})] + 2\langle \nabla f(\X), P_{\U_\perp} \theta_U \B \B^{-1} \B \theta_U^\top  P_{\U_\perp}    \rangle\\
	& = \Hess f(\X)[\cL_{\U,\B}^{r+}(\theta_{(\U,\B)}),\cL_{\U,\B}^{r+}(\theta_{(\U,\B)})],
\end{split}
\end{equation*} 
where the last equality follows from the expression of  $\Hess f(\X)$ in \eqref{eq: embedded-gd-hessian-psd} and the definition of $\cL_{\U,\B}^{r+}$.

Then,  by \eqref{ineq: bijection-spectrum-psd2}, \eqref{eq: Hessian-connection-PSD2} and Theorem \ref{th: hessian-sandwich}, we have $\overline{\Hess \, h_{r+}([\U,\B])}$ has $(pr-(r^2-r)/2)$ eigenvalues and $\widebar{\lambda}_i(\Hess \, h_{r+}([\U,\B]))$ is sandwiched between $(\sigma^2_r(\W_\B^{-1}) \wedge 2 \sigma_r^2(\V_\B^{-1/2} \B) ) \lambda_i(\Hess f(\X))$ and $(\sigma^2_1(\W^{-1}_\B) \vee 2 \sigma^2_1( \V_\B^{-1/2} \B ) ) \lambda_i(\Hess f(\X))$ for $i = 1,\ldots,pr-(r^2-r)/2$. \quad $\blacksquare$

\subsection{Proof of Corollary \ref{coro: landscape connection PSD}.}
Here we prove the Riemannian FOSP, SOSP, and strict saddle equivalence on $\cM_{r+}^e$ and $\cM_{r+}^{q_1}$, similar proof applies to the equivalence of \eqref{eq: PSD-manifold-formulation} on $\cM_{r+}^e$ and $\cM_{r+}^{q_2}$. First, by the connection of Riemannian gradients in \eqref{eq: gradient-connect-PSD1}, the connection of Riemannian FOSPs under two geometries clearly holds. 

Suppose $[\Y]$ is a Riemannian SOSP of \eqref{eq: PSD-opt-problem-quotient-sub1} and let $\X = \Y \Y^\top$. Given any  $\xi_\X \in T_{\X}\cM^e_{r+}$, we have $\Hess f(\X)[\xi_\X, \xi_\X] \overset{ \eqref{eq: Hessian-connection-PSD1} }=\overline{\Hess \, h_{r+}([\Y])}[(\cL^{r+}_\Y)^{-1}(\xi_\X), (\cL^{r+}_\Y)^{-1}(\xi_\X)] \geq 0$, where the inequality is by the SOSP assumption on $[\Y]$. Combining the fact $\X$ is a Riemannian FOSP, this shows $\X = \Y \Y^\top$ is a Riemannian SOSP under the embedded geometry.

Next, let us show the other direction: suppose $\X$ is a Riemannian SOSP under the embedded geometry, there is a unique $[\Y]$ such that $\Y \Y^\top = \X$ and it is a Riemannian SOSP of \eqref{eq: PSD-opt-problem-quotient-sub1}. To see this, first the uniqueness of $[\Y]$ is guaranteed by the fact $\bar{\ell}: \Y \in \widebar{\cM}_{r+}^{q_1} \to  \Y \Y^\top \in \cM_{r+}^e $ induces a diffeomorphism between $\cM_{r+}^e$ and $\cM_{r+}^{q_1}$ \cite[Proposition A.7]{massart2020quotient}. In addition, we have shown $[\Y]$ is a Riemannian FOSP of \eqref{eq: PSD-opt-problem-quotient-sub1}. Then by \eqref{eq: Hessian-connection-PSD1}, we have for any $\theta_\Y \in \cH_\Y \widebar{\cM}_{r+}^{q_1}$, $\overline{\Hess \, h_{r+}([\Y])}[\theta_\Y, \theta_\Y]  = \Hess f(\X)[\cL_\Y^{r+}(\theta_\Y),\cL_\Y^{r+}(\theta_\Y)] \geq 0$.

Finally, the equivalence on strict saddles also follows easily from the sandwich inequality from Theorem \ref{th: embedded-quotient-connection-PSD1} and the definition of the strict saddle. \quad $\blacksquare$ 

\section{Additional Proofs in Section \ref{sec: connection-general} } \label{sec: additional-proofs-general}

\subsection{Proof of Proposition \ref{prop: gradient-hessian-exp-general}.}
The proof is divided into five steps: in Step 1, we derive the expressions for Riemannian gradients; in Step 2, we derive the Riemannian Hessian of \eqref{eq: general prob} under the embedded geometry; in Step 3, we derive the Riemannian Hessian for $h_{r}([\L,\R])$; in Step 4, we derive the Riemannian Hessian for $h_{r}([\U,\B,\V])$; in Step 5, we derive the Riemannian Hessian for $h_r([\U,\Y])$.

{\bf Step 1.}  First, the Riemannian gradient expression of \eqref{eq: general prob} under the embedded geometry can be found in \cite[Proposition 1]{luo2021nonconvex}. Next, we compute $\grad \, \bar{h}_{r}(\L,\R)$, $\grad \, \bar{h}_{r}(\U, \B,\V)$ and $\grad \, \bar{h}_{r}(\U, \Y)$ from their definitions:
\begin{equation} \label{eq: quotient-gradient-comp-general-eq1}
\begin{split}
	& \quad \langle \grad_{\L} \, \bar{h}_r(\L,\R) \W_{\L,\R}, \eta_L \rangle + \langle \grad_{\R} \, \bar{h}_r(\L,\R) \V_{\L,\R}, \eta_R \rangle \\
	&= \bar{g}_{(\L,\R)}^r(\grad \, \bar{h}_{r}(\L,\R), \eta_{(\L,\R)} )\\
	 &= \rmD \bar{h}_{r}(\L,\R)[\eta_{(\L,\R)}] = \langle \nabla_{\L} \bar{h}_r(\L,\R) , \eta_L \rangle +\langle \nabla_{\R} \bar{h}_r(\L,\R) , \eta_R \rangle \\
	 &= \langle \nabla f(\L \R^\top) \R , \eta_L \rangle +\langle (\nabla f(\L \R^\top))^\top \L , \eta_R \rangle, \quad \forall \eta_{(\L,\R)} \in T_{(\L,\R)} \widebar{\cM}_r^{q_1},
\end{split}
\end{equation} where $\nabla_{\L} \bar{h}_r(\L,\R)$ denotes the Euclidean gradient of $\bar{h}_r(\L,\R)$ with respect to $\L$,
\begin{equation}\label{eq: quotient-gradient-comp-general-eq2}
	\begin{split}
		&\langle \grad_\U \, \bar{h}_r(\U,\B,\V), \eta_U \rangle + \langle \B^{-1} \grad_{\B} \, \bar{h}_r(\U,\B,\V) \B^{-1}, \eta_B \rangle + \langle \grad_{\V} \, \bar{h}_r(\U,\B,\V), \eta_V \rangle\\
		= & \bar{g}_{(\U,\B,\V)}^r (\grad \, \bar{h}_r(\U,\B,\V), \eta_{(\U,\B,\V)} ) = \rmD \bar{h}_r(\U,\B,\V)[\eta_{(\U,\B,\V)}] \\
		= & \langle \nabla_{\U} \bar{h}_r(\U,\B,\V) , \eta_U \rangle +\langle \nabla_{\B} \bar{h}_r(\U,\B,\V) , \eta_B \rangle + \langle \nabla_{\V} \bar{h}_r(\U,\B,\V) , \eta_V \rangle \\
		= & \langle \nabla f(\U \B \V^\top) \V \B , \eta_U \rangle +\langle \U^\top \nabla f(\U \B \V^\top) \V, \eta_B \rangle \\
		&+ \langle (\nabla f(\U \B \V^\top))^\top \U \B , \eta_V \rangle, \quad \forall \eta_{(\U,\B,\V)} \in T_{(\U,\B,\V)} \widebar{\cM}_r^{q_2},
	\end{split}
\end{equation} and
\begin{equation}\label{eq: quotient-gradient-comp-general-eq3}
	\begin{split}
		& \quad \langle \grad_\U \, \bar{h}_r(\U,\Y) \V_\Y, \eta_U \rangle + \langle \grad_\Y \, \bar{h}_r(\U,\Y) \W_\Y, \eta_Y \rangle\\
		 &= \bar{g}_{(\U,\Y)}^r (\grad \, \bar{h}_r(\U,\Y), \eta_{(\U,\Y)} ) =\rmD \bar{h}_r(\U,\Y)[\eta_{(\U,\Y)}] \\
		& = \langle \nabla_\U \bar{h}_r(\U,\Y) , \eta_U \rangle +\langle \nabla_\Y \bar{h}_r(\U,\Y) , \eta_Y \rangle\\
		& = \langle \nabla f(\U \Y^\top) \Y , \eta_U \rangle +\langle (\nabla f(\U \Y^\top))^\top \U  , \eta_Y \rangle, \quad \forall \eta_{(\U,\Y)} \in T_{(\U,\Y)} \widebar{\cM}_r^{q_3}.
	\end{split}
\end{equation}
From \eqref{eq: quotient-gradient-comp-general-eq1} and \eqref{eq: quotient-gradient-comp-general-eq2}, we have
\begin{equation} \label{eq: RGradient-total-space-general}
	\begin{split}
		\grad \, \bar{h}_{r}(\L,\R) &=\begin{bmatrix}
			\nabla f(\L \R^\top) \R \W_{\L,\R}^{-1} \\
			(\nabla f(\L \R^\top))^\top \L \V_{\L,\R}^{-1}
		\end{bmatrix},   \\
		\grad \, \bar{h}_r(\U,\B,\V) &= \begin{bmatrix}
			P_{ T_\U \st(r,p_1) } \left(  \nabla f(\U \B \V^\top) \V \B \right)\\
			\B\sym( \U^\top \nabla f(\U \B \V^\top) \V )\B\\
			P_{ T_\V \st(r,p_2) } \left( (\nabla f(\U \B \V^\top))^\top \U \B  \right)
  		\end{bmatrix}.
	\end{split}
\end{equation}
Moreover, from \eqref{eq: quotient-gradient-comp-general-eq3}, we know $\grad_\Y \, \bar{h}_r(\U,\Y)= (\nabla f(\U \Y^\top))^\top \U \W_\Y^{-1}$. Suppose $\grad_\U \, \bar{h}_r(\U,\Y) = \U_\perp \D + \U \bOmega \in T_\U \st(r,p_1)$. Then from \eqref{eq: quotient-gradient-comp-general-eq3}, we have $P_{T_\U \st(r,p_1)}(  \grad_\U \, \bar{h}_r(\U,\Y) \V_\Y -  \nabla f(\U \Y^\top) \Y  ) = 0$. This yields $\D = \U_\perp^\top \nabla f(\U \Y^\top) \Y \V_\Y^{-1} $ and $\bOmega$ is uniquely determined by the Sylvester equation: $\V_\Y \bOmega^\top + \bOmega^\top \V_\Y = \Y^\top (\nabla f(\U \Y^\top))^\top \U - \U^\top \nabla f(\U \Y^\top) \Y  $. So 
\begin{equation*}
		\grad \, \bar{h}_r(\U,\Y) = \begin{bmatrix}
			P_{\U_\perp}\nabla f(\U \Y^\top) \Y \V_\Y^{-1}  + \U \bOmega  \\
			(\nabla f(\U \Y^\top))^\top \U \W_\Y^{-1}
		\end{bmatrix}.
\end{equation*}

By Lemma \ref{lm: general-quotient-manifold1-prop}, $\cV_{(\L,\R)} \widebar{\cM}_{r}^{q_1} $ is orthogonal to $ \cH_{(\L,\R)} \widebar{\cM}_{r}^{q_1}$ with respect to $\bar{g}^{r}_{(\L,\R)}$, we have $\overline{\grad\, h_{r}([\L,\R])} = \grad \, \bar{h}_{r}(\L,\R)$ by Lemma \ref{lm: quotient-ortho-hori-gradient-express}. 

Next, we compute $\overline{\grad\, h_{r}([\U,\B,\V])}$. Given $\eta_{(\U,\B, \V)} = [(\U \bOmega_1 + \U_\perp \D_1)^\top \quad \S \quad  ( \V \bOmega_2 + \V_\perp \D_2)^\top]^\top\in T_{(\U,\B,\V)} \widebar{\cM}_r^{q_2}$, and suppose $P_{(\U,\B,\V)}^\cH (\eta_{(\U,\B, \V)}) = [ (\U \widetilde{\bOmega} + \U_\perp \widetilde{\D}_1)^\top \quad \widetilde{\S} \quad (-\V \widetilde{\bOmega} + \V_\perp \widetilde{\D}_2)^\top ] \in \cH_{(\U,\B,\V)} \widebar{\cM}_r^{q_2}$. By definition 
 \begin{equation*} 
 	\begin{split}
 		&(\widetilde{\bOmega}, \widetilde{\S}, \widetilde{\D}_1, \widetilde{\D}_2) \\
 		= &\argmin_{ \substack{\bOmega' = -\bOmega^{'\top} \in \bbR^{r\times r}, \S' \in \bbS^{r \times r},\\ \D_1' \in \bbR^{(p_1 - r) \times r}, \D_2' \in \bbR^{(p_2 - r) \times r}} } \bar{g}_{(\U,\B,\V)}^r \left( \begin{bmatrix}
 \U (\bOmega_1 - \bOmega') + \U_\perp (\D_1 - \D_1') \\
 \S-\S' \\
 \V (\bOmega_2 + \bOmega') + \V_\perp (\D_2 - \D_2')	
 \end{bmatrix}, \begin{bmatrix}
 \U (\bOmega_1 - \bOmega') + \U_\perp (\D_1 - \D_1') \\
 \S-\S' \\
 \V (\bOmega_2 + \bOmega') + \V_\perp (\D_2 - \D_2')	
 \end{bmatrix}\right) \\
 = &\argmin_{ \substack{\bOmega' = -\bOmega^{'\top} \in \bbR^{r\times r}, \S' \in \bbS^{r \times r},\\ \D_1' \in \bbR^{(p_1 - r) \times r}, \D_2' \in \bbR^{(p_2 - r) \times r}} } \|\bOmega_1 - \bOmega'\|^2_\F + \|\bOmega_2 + \bOmega'\|_\F^2 + \sum_{i=1}^2\| \D_i - \D_i'\|_\F^2  + \|\B^{-1/2} (\S - \S') \B^{-1/2} \|_\F^2.
 	\end{split}
 \end{equation*} So we have $(\widetilde{\bOmega}, \widetilde{\S}, \widetilde{\D}_1, \widetilde{\D}_2) = ( (\bOmega_1 - \bOmega_2)/2, \S, \D_1, \D_2)$, i.e.,
 \begin{equation} \label{eq: horizontal-projection-UBV}
	P_{(\U,\B,\V)}^\cH(\eta_{(\U,\B,\V)}) = \begin{bmatrix}
		P_{\U_\perp} \eta_U + \U(\U^\top \eta_U - \V^\top \eta_V )/2\\
		\eta_B \\
		P_{\V_\perp} \eta_V - \V(\U^\top \eta_U - \V^\top \eta_V )/2
	\end{bmatrix}, \quad \eta_{(\U,\B,\V)} \in T_{(\U,\B,\V)} \widebar{\cM}_r^{q_2}.
\end{equation}
 Then recall the definition of the projection operator $P_{ T_\U \st(r,p_1) }$ from Table \ref{tab: basic-prop-simple-manifold}, we have
  \begin{equation*}
  	\begin{split}
  		\overline{\grad\, h_{r}([\U,\B,\V])} \overset{ \eqref{eq: quotient-gradient-general} } =& P_{(\U,\B,\V)}^\cH \left( \grad \, \bar{h}_r(\U,\B,\V) \right) 
  		= \begin{bmatrix}
				P_{\U_\perp} \nabla f(\U \B \V^\top) \V \B + \U \bOmega \\
				\B \sym(\U^\top \nabla f(\U \B \V^\top) \V) \B \\
				P_{\V_\perp} (\nabla f(\U \B \V^\top))^\top \U \B - \V \bOmega 
			\end{bmatrix},
  	\end{split}
  \end{equation*} where 
  \begin{equation*}
  	\begin{split}
  		\bOmega &= \left(\skew\left( \U^\top \nabla f(\U \B \V^\top) \V \B  \right) - \skew\left(\V^\top (\nabla f(\U \B \V^\top))^\top \U \B \right)\right)/2\\
  		& = \skew(\U^\top \nabla f(\U \B \V^\top) \V ) \B/2 + \B \skew(\U^\top \nabla f(\U \B \V^\top) \V )/2.
  	\end{split}
  \end{equation*}

Finally,
\begin{equation*}
	\begin{split}
		\overline{\grad\, h_{r}([\U,\Y])} \overset{ \eqref{eq: quotient-gradient-general} } = P_{(\U,\V)}^{\cH}( \grad \, \bar{h}_r(\U,\Y) )  &=\begin{bmatrix}
			P_{\U_\perp}\left( P_{\U_\perp}\nabla f(\U \Y^\top) \Y \V_\Y^{-1}  + \U \bOmega  \right)\\
			(\nabla f(\U \Y^\top))^\top \U \W_\Y^{-1}
		\end{bmatrix}\\
		& = \begin{bmatrix}
			P_{\U_\perp} \nabla f(\U \Y^\top) \Y \V^{-1}_\Y \\
			(\nabla f(\U \Y^\top))^\top \U \W_\Y^{-1}
		\end{bmatrix}.
	\end{split}
\end{equation*}

{\bf Step 2.} Suppose $\X$ has has SVD $\U' \bSigma' \V^{'\top}$ with $\U = \U' \O_1,\V = \V' \O_2$ ($\O_1,\O_2 \in \bbO_r$). Then for
\begin{equation*}
	\begin{split}
		\xi_\X &= [\U \quad \U_\perp] \begin{bmatrix}
			\S & \D_2^\top\\
			\D_1 & \0
		\end{bmatrix} [\V \quad \V_\perp]^\top = [\U' \quad \U_\perp] \begin{bmatrix}
			\O_1\S\O_2^\top & \O_1\D_2^\top\\
			\D_1 \O_2^\top & \0
		\end{bmatrix} [\V' \quad \V_\perp]^\top,
	\end{split}
\end{equation*} by \cite[Proposition 2]{luo2021nonconvex} we have in the general case
\begin{equation*}
	\begin{split}
		\Hess f(\X)[\xi_\X, \xi_\X] &= \nabla^2 f(\X)[\xi_\X, \xi_\X] + 2\langle \nabla f(\X), \U_\perp \D_1 \O_2^\top \bSigma^{'-1} \O_1 \D_2^\top \V_\perp^\top \rangle \\
		& =  \nabla^2 f(\X)[\xi_\X, \xi_\X] + 2\langle \nabla f(\X), \U_\perp \D_1 \bSigma^{-1} \D_2^\top \V_\perp^\top \rangle,
	\end{split}
\end{equation*} where $\bSigma = \O_1^\top \bSigma' \O_2 = \U^\top \X \V$.

{\bf Step 3.} For convenience, we denote $\bar{g}_{\L}^{r}(\eta_L, \theta_L) = \tr(\W_{\L,\R}\eta_L^\top \theta_L) $ as the $\L$ part inner product of $\bar{g}^r_{(\L, \R)}$ and $\bar{g}_{\R}^{r}(\eta_R, \theta_R) = \tr(\V_{\L,\R}\eta_R^\top \theta_R)$ as the $\R$ part inner product of $\bar{g}^r_{(\L, \R)}$. First following similar arguments as in \eqref{eq: directional-derivative-gY-instance} and \eqref{eq: koszul-formula-derive-gY}, for three vector fields $\eta, \theta, \xi$ on $\widebar{\cM}_r^{q_1}$, we have
\begin{equation*} \label{eq: directional-derivative-gLR}
\begin{split}
	\rmD \bar{g}_{(\L,\R)} (\eta, \theta)[\xi] = \tr(\rmD \W_{\L,\R} [\xi] \eta_L^\top \theta_L ) + \tr(\rmD \V_{\L,\R} [\xi] \eta_R^\top \theta_R ) +\bar{g}_{(\L,\R)}( \rmD \eta[\xi], \theta ) + \bar{g}_{(\L,\R)}( \eta, \rmD \theta[\xi] )
\end{split}
\end{equation*} and
\begin{equation} \label{eq: koszul-formula-gLR}
	\begin{split}
		2 \bar{g}_{(\L,\R)} (\widebar{\nabla}_\xi \eta, \theta) &= 2 \bar{g}_{(\L,\R)} ( \rmD \eta[\xi], \theta ) + \tr( \rmD \W_{\L,\R} [\xi] \eta_L^\top \theta_L ) + \tr( \rmD \V_{\L,\R} [\xi] \eta_R^\top \theta_R ) + \tr( \rmD \W_{\L,\R} [\eta] \xi_L^\top \theta_L )\\
		& \quad + \tr( \rmD \V_{\L,\R} [\eta] \xi_R^\top \theta_R )- \tr( \rmD \W_{\L,\R} [\theta] \eta_L^\top \xi_L ) - \tr( \rmD \V_{\L,\R} [\theta] \eta_R^\top \xi_R ).
	\end{split}
\end{equation}
 Then
\begin{equation} \label{eq: koszul-formula-LR}
	\begin{split}
		 &\quad \overline{\Hess \, h_{r}([\L,\R])}[\theta_{(\L,\R)}, \theta_{(\L,\R)}] \\
		 &= \bar{g}_{(\L,\R)}^{r}\left( \overline{\Hess \, h_{r}([\L,\R])[\theta_{[\L,\R]}]} , \theta_{(\L,\R)} \right)\\
		&\overset{ \eqref{eq: quotient-hessian-linear-from} } =\bar{g}_{(\L,\R)}^{r} \left( P_{(\L,\R)}^{\cH} \left( \widebar{\nabla}_{\theta_{(\L,\R)}} \overline{ \grad \, h_{r}([\L,\R]) } \right) , \theta_{(\L,\R)} \right)\\
		 &\overset{ \text{Lemma } \ref{lm: quotient-hessian-riemannian-connection-equal} } = \bar{g}_{(\L,\R)}^{r} \left( \widebar{\nabla}_{\theta_{(\L,\R)}} \overline{ \grad \, h_{r}([\L,\R])} , \theta_{(\L,\R)} \right) \\
		 & \overset{ \eqref{eq: koszul-formula-gLR} } = \bar{g}_{\L}^{r} \left( \rmD \overline{ \grad_{\L} \, h_{r}([\L,\R])}[\theta_{(\L,\R)}] , \theta_{L} \right) + \bar{g}_{\R}^{r} \left( \rmD \overline{ \grad_{\R} \, h_{r}([\L,\R])}[\theta_{(\L, \R)}], \theta_R \right)\\
		 & \quad  + \tr( \rmD \W_{\L,\R}[ \overline{ \grad \, h_{r}([\L,\R])} ] \theta_L^\top \theta_L )/2 +  \tr( \rmD \V_{\L,\R}[ \overline{ \grad \, h_{r}([\L,\R])} ] \theta_R^\top \theta_R )/2. 
	\end{split}
\end{equation}

Next, we compute $\rmD \overline{ \grad_{\L} \, h_{r}([\L,\R])}[\theta_{(\L,\R)}]$.
\begin{equation} \label{eq: directional-derivative-L}
	\begin{split}
		&\rmD \overline{ \grad_{\L} \, h_{r}([\L,\R])}[\theta_{(\L,\R)}] \\
		=& \lim_{t \to 0} \left( \nabla f\left((\L+ t\theta_L) (\R + t\theta_R )^\top\right) (\R + t\theta_R ) \W_{\L+ t\theta_L,\R + t\theta_R}^{-1} -\nabla f(\L \R^\top) \R \W_{\L,\R}^{-1}  \right)/t \\
		=& \nabla f(\L \R^\top) \theta_R \W_{\L,\R}^{-1} + \nabla f(\L \R^\top) \R \rmD \W_{\L,\R}^{-1}[\theta_{(\L,\R)}] + \nabla^2 f(\L \R^\top)[\L \theta_R^\top + \theta_L \R^\top] \R  \W_{\L,\R}^{-1}.
	\end{split}
\end{equation} Similarly we have
\begin{equation} \label{eq: directional-derivative-R}
	\begin{split} 
		&\rmD \overline{ \grad_{\R} \, h_{r}([\L,\R])}[\theta_{(\L, \R)}] \\
		=& \left(\nabla f(\L \R^\top) \right)^\top \theta_L \V_{\L,\R}^{-1} + \left(\nabla f(\L \R^\top) \right)^\top \L \rmD \V_{\L,\R}^{-1}[\theta_{(\L,\R)}] + \left(\nabla^2 f(\L \R^\top)[\L \theta_R^\top + \theta_L \R^\top]\right)^\top \L  \V_{\L,\R}^{-1}.
	\end{split}
\end{equation}
Plugging \eqref{eq: directional-derivative-L} and \eqref{eq: directional-derivative-R} into \eqref{eq: koszul-formula-LR}, we have
\begin{equation*}
	\begin{split}
		 &\quad \overline{\Hess \, h_{r}([\L,\R])}[\theta_{(\L,\R)}, \theta_{(\L,\R)}] \\
		 &=  \nabla^2 f(\L \R^\top)[\L \theta_R^\top + \theta_L \R^\top,\L \theta_R^\top + \theta_L \R^\top] + 2 \langle \nabla f(\L \R^\top), \theta_L \theta_R^\top \rangle \\
		 & \quad + \langle \nabla f(\L \R^\top) \R \rmD \W_{\L,\R}^{-1}[\theta_{(\L,\R)}] , \theta_L \W_{\L,\R} \rangle + \langle  \left(\nabla f(\L \R^\top) \right)^\top \L \rmD \V_{\L,\R}^{-1}[\theta_{(\L,\R)}], \theta_R \V_{\L,\R} \rangle \\
		 & \quad +  \tr( \rmD \W_{\L,\R}[ \overline{ \grad \, h_{r}([\L,\R])} ] \theta_L^\top \theta_L )/2 +  \tr( \rmD \V_{\L,\R}[ \overline{ \grad \, h_{r}([\L,\R])} ] \theta_R^\top \theta_R )/2. 
	\end{split}
\end{equation*}

{\bf Step 4.} For convenience, we denote $\bar{g}_{\U}^{r}(\eta_U, \theta_U) = \tr(\eta_U^\top \theta_U) $ as the $\U$ part inner product of $\bar{g}^r_{(\U, \B,\V)}$, $\bar{g}^r_{\B}(\eta_B, \theta_B) = \tr(\B^{-1} \eta_B \B^{-1} \theta_B)$ as the $\B$ part inner product of $\bar{g}^r_{(\U, \B,\V)}$ and $\bar{g}_{\V}^{r}(\eta_V, \theta_V) = \tr(\eta_V^\top \theta_V) $ as the $\V$ part inner product of $\bar{g}^r_{(\U, \B,\V)}$. Then
\begin{equation} \label{eq: quotient-hessian-UBV}
	\begin{split}
		& \quad \overline{\Hess \, h_{r}([\U,\B,\V])}[\theta_{(\U,\B,\V)}, \theta_{(\U,\B,\V)}] \\
		&= \bar{g}_{(\U,\B,\V)}^{r}\left( \overline{\Hess \, h_{r}([\U,\B,\V])[\theta_{[\U,\B,\V]}]} , \theta_{(\U,\B,\V)} \right)\\
		& \overset{\eqref{eq: quotient-hessian-linear-from}, \text{Lemma } \ref{lm: quotient-hessian-riemannian-connection-equal} } = \bar{g}_{(\U,\B,\V)}^{r} \left( \widebar{\nabla}_{\theta_{(\U,\B,\V)}} \overline{ \grad \, h_{r}([\U,\B,\V])} , \theta_{(\U, \B,\V)} \right)\\
		& \overset{ \eqref{eq: koszul-formula-simple-manifolds} } = \bar{g}_{\U}^{r} \left( \rmD \overline{ \grad_\U \, h_{r}([\U, \B,\V])}[\theta_{(\U, \B,\V)}] , \theta_{U} \right) + \bar{g}_{\B}^{r} \left( \rmD \overline{ \grad_\B \, h_{r}([\U, \B,\V])}[\theta_{(\U, \B,\V)}], \theta_B \right)\\
		& \quad - \tr\left(\B^{-1} \sym\left( \theta_B \B^{-1} \overline{ \grad_\B \, h_{r}([\U, \B,\V])} \right) \B^{-1} \theta_B \right) + \bar{g}_{\V}^{r} \left( \rmD \overline{ \grad_\V \, h_{r}([\U, \B,\V])}[\theta_{(\U, \B,\V)}] , \theta_{V} \right).
	\end{split}
\end{equation}
Recall the horizontal projection operator $P_{(\U,\B,\V)}^\cH$ given in \eqref{eq: horizontal-projection-UBV}, and let $P_{\U}^\cH$, $P_{\B}^\cH$ and $P_{\V}^\cH$ be the restriction of $P_{(\U,\B,\V)}^\cH$ to the $\U$, $\B$ and $\V$ components, respectively. Moreover, define $\U_t = \U + t \theta_U, \V_t = \V + t \theta_V, \B_t = \B + t\theta_B$. By \eqref{eq: RGradient-total-space-general}, we have
\begin{equation} \label{eq: directional-derivative-gradient-UBV}
	\begin{split}
		\bDelta_U:=&\lim_{t \to 0} \left( \grad_\U \, \bar{h}_r(\U_t,\B_t,\V_t) -  \grad_\U \, \bar{h}_r(\U,\B,\V))  \right)/t\\
		=& 	P_{ T_\U \st(r,p_1) } \left(  \nabla f(\U \B \V^\top) \theta_V \B + \nabla f(\U \B \V^\top) \V \theta_B + \nabla^2 f(\U \B \V^\top)[\theta_U \B \V^\top + \U \theta_B \V^\top + \U \B \theta_V^\top] \V \B  \right)\\
		& - 2 \sym(\U \theta_U^\top) \nabla f(\U \B \V^\top) \V \B + \theta_U \skew\left(\U^\top\nabla f(\U \B \V^\top) \V \B \right) + \U \skew\left(  \theta_U^\top\nabla f(\U \B \V^\top) \V \B \right);\\
		\bDelta_V:= &\lim_{t \to 0} \left( \grad_\V \, \bar{h}_r(\U_t,\B_t,\V_t) -  \grad_\V \, \bar{h}_r(\U,\B,\V))  \right)/t\\
		=& 	P_{ T_\V \st(r,p_2) } \Big(  \left(\nabla f(\U \B \V^\top)\right)^\top \theta_U \B + \left(\nabla f(\U \B \V^\top)\right)^\top \U \theta_B \\
		&+ \left(\nabla^2 f(\U \B \V^\top)[\theta_U \B \V^\top + \U \theta_B \V^\top + \U \B \theta_V^\top]\right)^\top \U \B  \Big)- 2 \sym(\V \theta_V^\top) \left(\nabla f(\U \B \V^\top)\right)^\top \U \B\\
		& + \theta_V \skew\left(\V^\top\left(\nabla f(\U \B \V^\top)\right)^\top \U \B \right) + \V \skew\left(  \theta_V^\top\left(\nabla f(\U \B \V^\top)\right)^\top \U \B \right). 
	\end{split}
\end{equation}
Then
\begin{equation} \label{eq: U-part-UBV-directional-derivative}
	\begin{split}
		&\rmD \overline{ \grad_\U \, h_{r}([\U, \B,\V])}[\theta_{(\U, \B,\V)}] \\
		=& \lim_{t \to 0} \Big(  P_{\U_t}^\cH (  \grad_\U \, \bar{h}_r(\U_t,\B_t,\V_t) ) - P_\U^\cH ( \grad_\U \, \bar{h}_r(\U,\B,\V)) \Big)/t \\
		=& P_\U^\cH \left( \lim_{t \to 0} \left( \grad_\U \, \bar{h}_r(\U_t,\B_t,\V_t) -  \grad_\U \, \bar{h}_r(\U,\B,\V)  \right)/t  \right) - \sym(\U \theta_U^\top) \grad_\U \, \bar{h}_r(\U,\B,\V) \\
		& - (\theta_U \V^\top + \U \theta_V^\top) \grad_\V \, \bar{h}_r(\U,\B,\V)/2 \\
		\overset{ \eqref{eq: directional-derivative-gradient-UBV} } = & P_\U^\cH (\bDelta_U)  - \sym(\U \theta_U^\top) \grad_\U \, \bar{h}_r(\U,\B,\V) - (\theta_U \V^\top + \U \theta_V^\top) \grad_\V \, \bar{h}_r(\U,\B,\V)/2,
	\end{split}
\end{equation}
and
\begin{equation} \label{eq: V-part-UBV-directional-derivative}
	\begin{split}
		&\rmD \overline{ \grad_\V \, h_{r}([\U, \B,\V])}[\theta_{(\U, \B,\V)}] \\
		=& \lim_{t \to 0} \Big(  P_{\V_t}^\cH (  \grad_\V \, \bar{h}_r(\U_t,\B_t,\V_t) ) - P_\V^\cH ( \grad_\V \, \bar{h}_r(\U,\B,\V)) \Big)/t \\
		=& P_\V^\cH \left( \lim_{t \to 0} \left( \grad_\V \, \bar{h}_r(\U_t,\B_t,\V_t) -  \grad_\V \, \bar{h}_r(\U,\B,\V)  \right)/t  \right) - \sym(\V \theta_V^\top) \grad_\V \, \bar{h}_r(\U,\B,\V) \\
		& - (\theta_V \U^\top + \V \theta_U^\top) \grad_\V \, \bar{h}_r(\U,\B,\V)/2 \\
		\overset{ \eqref{eq: directional-derivative-gradient-UBV} } = & P_\U^\cH (\bDelta_V)  - \sym(\V \theta_V^\top) \grad_\V \, \bar{h}_r(\U,\B,\V) - (\theta_V \U^\top + \V \theta_U^\top) \grad_\U   \, \bar{h}_r(\U,\B,\V)/2.
	\end{split}
\end{equation}

In addition
\begin{equation} \label{eq: B-part-UBV-directional-derivative}
	\begin{split}
		&\rmD \overline{ \grad_\B \, h_{r}([\U, \B,\V])}[\theta_{(\U, \B,\V)}] \\
		=& \lim_{t \to 0} \Big(   \grad_\B \, \bar{h}_r(\U_t,\B_t,\V_t)  - \grad_\B \, \bar{h}_r(\U,\B,\V) \Big)/t \\
		=& \B\Big( \sym\left(\theta_U^\top \nabla f(\U \B \V^\top) \V \right) + \sym\left(\U^\top \nabla f(\U \B \V^\top) \theta_V\right) \\
		&+ \U^\top  \nabla^2 f(\U \B \V^\top)[\theta_U \B \V^\top + \U \theta_B \V^\top + \U \B \theta_V^\top] \V   \Big) \B \\
		& + \B\sym(\U^\top \nabla f(\U \B \V^\top) \V ) \theta_B + \theta_B \sym(\U^\top \nabla f(\U \B \V^\top) \V ) \B.
	\end{split}
\end{equation}

Let $\bDelta =  \U^\top\nabla f(\U \B \V^\top) \V$, $\bDelta' =  \theta_U^\top\nabla f(\U \B \V^\top) \V$ and $\bDelta'' =  \U^\top\nabla f(\U \B \V^\top) \theta_V$. By \eqref{eq: U-part-UBV-directional-derivative}, \eqref{eq: V-part-UBV-directional-derivative} and \eqref{eq: B-part-UBV-directional-derivative} and observing the fact $\theta_{(\U,\B, \V)} \in \cH_{(\U,\B,\V)} \widebar{\cM}_r^{q_2}$, we have
\begin{equation} \label{eq: U-part-UBV-directional-derivative2}
	\begin{split}
		 &\bar{g}_{\U}^{r} \left( \rmD \overline{ \grad_\U \, h_{r}([\U, \B,\V])}[\theta_{(\U, \B,\V)}] , \theta_{U} \right) \\
		 =& \left\langle \bDelta_U  - \sym(\U \theta_U^\top) \grad_\U \, \bar{h}_r(\U,\B,\V) - (\theta_U \V^\top + \U \theta_V^\top) \grad_\V \, \bar{h}_r(\U,\B,\V)/2, \theta_U  \right\rangle \\
		=&  \left\langle \bDelta_U , \theta_U  \right\rangle+ \langle \bDelta, \sym(\U^\top \theta_U \U^\top\theta_U) \B + \B \sym(\V^\top \theta_V \U^\top \theta_U) \rangle/2 - \langle \bDelta', \U^\top \theta_U \B \rangle/2 - \langle \bDelta'', \B \theta_U^\top \U \rangle/2;\\
		\overset{\eqref{eq: directional-derivative-gradient-UBV} } =&  \nabla^2 f(\U \B \V^\top)[\theta_U \B \V^\top + \U \theta_B \V^\top + \U \B \theta_V^\top, \theta_U \B \V^\top] + \langle \nabla f(\U\B\V^\top), \theta_U \B \theta_V^\top \rangle  \\
		& - \langle \bDelta, \theta_U^\top \theta_U \B \rangle + \langle \bDelta', 2\theta_B -  \U^\top \theta_U \B - \theta_U^\top \U \B \rangle/2 +   \langle \bDelta, \sym(\U^\top \theta_U \U^\top\theta_U) \B + \B \sym(\V^\top \theta_V \U^\top \theta_U) \rangle/2  \\
		& - \langle \bDelta', \U^\top \theta_U \B \rangle/2 - \langle \bDelta'', \B \theta_U^\top \U \rangle/2\\
		= & \nabla^2 f(\U \B \V^\top)[\theta_U \B \V^\top + \U \theta_B \V^\top + \U \B \theta_V^\top, \theta_U \B \V^\top] + \langle \nabla f(\U\B\V^\top), \theta_U \B \theta_V^\top \rangle  \\
		& + \left\langle \bDelta,  \sym(\U^\top \theta_U \U^\top\theta_U) \B + \B \sym(\V^\top \theta_V \U^\top \theta_U) -2\theta_U^\top \theta_U \B \right\rangle/2 \\
		& + \langle \bDelta', \theta_B -  \U^\top \theta_U \B - \theta_U^\top \U \B/2 \rangle - \langle \bDelta'', \B \theta_U^\top \U \rangle/2,
	\end{split}
\end{equation} and similarly
\begin{equation} \label{eq: V-part-UBV-directional-derivative2}
	\begin{split}
		& \bar{g}_{\V}^{r} \left( \rmD \overline{ \grad_\V \, h_{r}([\U, \B,\V])}[\theta_{(\U, \B,\V)}] , \theta_{V} \right)  \\
		= & \nabla^2 f(\U \B \V^\top)[\theta_U \B \V^\top + \U \theta_B \V^\top + \U \B \theta_V^\top, \U \B \theta_V^\top] + \langle \nabla f(\U\B\V^\top), \theta_U \B \theta_V^\top \rangle\\
		& + \left\langle \bDelta,  \B\sym(\V^\top \theta_V \V^\top\theta_V) + \sym(\U^\top \theta_U \V^\top \theta_V) \B -2\B\theta_V^\top \theta_V \right\rangle/2 \\
		& + \langle \bDelta'', \theta_B -  \B\theta_V^\top \V - \B\V^\top \theta_V/2 \rangle - \langle \bDelta', \V^\top \theta_V \B \rangle/2,
	\end{split}
\end{equation}
 and 
 \begin{equation} \label{eq: B-part-UBV-directional-derivative2}
 	\begin{split}
 		 & \bar{g}_{\B}^{r} \left( \rmD \overline{ \grad_\B \, h_{r}([\U, \B,\V])}[\theta_{(\U, \B,\V)}], \theta_B \right)  - \tr\left(\B^{-1} \sym\left( \theta_B \B^{-1} \overline{ \grad_\B \, h_{r}([\U, \B,\V])} \right) \B^{-1} \theta_B \right) \\
 		 = & \nabla^2 f(\U \B \V^\top)[\theta_U \B \V^\top + \U \theta_B \V^\top + \U \B \theta_V^\top, \U \theta_B \V^\top] + 2 \langle \bDelta, \theta_B \B^{-1} \theta_B \rangle \\
 		 & + \langle \bDelta', \theta_B \rangle + \langle \bDelta'', \theta_B \rangle - \langle \bDelta, \theta_B \B^{-1} \theta_B \rangle\\
 		 = & \nabla^2 f(\U \B \V^\top)[\theta_U \B \V^\top + \U \theta_B \V^\top + \U \B \theta_V^\top, \U \theta_B \V^\top] + \langle \bDelta, \theta_B \B^{-1} \theta_B \rangle + \langle \bDelta', \theta_B \rangle + \langle \bDelta'', \theta_B \rangle. 
 	\end{split}
 \end{equation} By plugging \eqref{eq: U-part-UBV-directional-derivative2}, \eqref{eq: V-part-UBV-directional-derivative2} and \eqref{eq: B-part-UBV-directional-derivative2} into \eqref{eq: quotient-hessian-UBV}, we finally have
\begin{equation*}
	\begin{split}
		 &\overline{\Hess \, h_{r}([\U,\B,\V])}[\theta_{(\U,\B,\V)}, \theta_{(\U,\B,\V)}] \\
		 =& \bar{g}_{\U}^{r} \left( \rmD \overline{ \grad_\U \, h_{r}([\U, \B,\V])}[\theta_{(\U, \B,\V)}] , \theta_{U} \right) + \bar{g}_{\B}^{r} \left( \rmD \overline{ \grad_\B \, h_{r}([\U, \B,\V])}[\theta_{(\U, \B,\V)}], \theta_B \right)\\
		& - \tr\left(\B^{-1} \sym\left( \theta_B \B^{-1} \overline{ \grad_\B \, h_{r}([\U, \B,\V])} \right) \B^{-1} \theta_B \right) + \bar{g}_{\V}^{r} \left( \rmD \overline{ \grad_\V \, h_{r}([\U, \B,\V])}[\theta_{(\U, \B,\V)}] , \theta_{V} \right)\\
		=& \nabla^2 f(\U \B \V^\top)[\theta_U \B \V^\top + \U \theta_B \V^\top + \U \B \theta_V^\top, \theta_U \B \V^\top + \U \theta_B \V^\top + \U \B \theta_V^\top] + 2\langle \nabla f(\U\B\V^\top), \theta_U \B \theta_V^\top \rangle\\
		& +   \left\langle \bDelta,  \sym(\U^\top \theta_U \U^\top\theta_U) \B + \B \sym(\V^\top \theta_V \U^\top \theta_U) -2\theta_U^\top \theta_U \B \right\rangle/2\\
		& + \left\langle \bDelta,  \B\sym(\V^\top \theta_V \V^\top\theta_V) + \sym(\U^\top \theta_U \V^\top \theta_V) \B -2\B\theta_V^\top \theta_V +  2\theta_B \B^{-1} \theta_B \right\rangle/2\\
		& + \langle \bDelta', 2\theta_B -  \U^\top \theta_U \B - \theta_U^\top \U \B/2 -\V^\top \theta_V \B/2 \rangle + \langle \bDelta'', 2\theta_B -  \B\theta_V^\top \V - \B\V^\top \theta_V/2 -\B \theta_U^\top \U /2 \rangle.
	\end{split}
\end{equation*}

{\bf Step 5.}  For convenience, we denote $\bar{g}_{\U}^{r}(\eta_U, \theta_U) = \tr(\V_\Y \eta_U^\top \theta_U) $ as the $\U$ part inner product of $\bar{g}^r_{(\U, \Y)}$ and $\bar{g}_{\Y}^{r}(\eta_Y, \theta_Y) = \tr(\W_\Y\eta_Y^\top \theta_Y)$ as the $\Y$ part inner product of $\bar{g}^r_{(\U, \Y)}$.  Then
\begin{equation} \label{eq: koszul-formula-UY}
	\begin{split}
		 &\quad \overline{\Hess \, h_{r}([\U,\Y])}[\theta_{(\U,\Y)}, \theta_{(\U,\Y)}] \\
		 &= \bar{g}_{(\U,\Y)}^{r}\left( \overline{\Hess \, h_{r}([\U,\Y])[\theta_{[\U,\Y]}]} , \theta_{(\U,\Y)} \right)\\
		 &\overset{\eqref{eq: quotient-hessian-linear-from}, \text{Lemma } \ref{lm: quotient-hessian-riemannian-connection-equal} } = \bar{g}_{(\U,\Y)}^{r} \left( \widebar{\nabla}_{\theta_{(\U,\Y)}} \overline{ \grad \, h_{r}([\U,\Y])} , \theta_{(\U,\Y)} \right) \\
		 & \overset{ \eqref{eq: koszul-formula-simple-manifolds} } = \bar{g}_{\U}^{r} \left( \rmD \overline{ \grad_\U \, h_{r}([\U,\Y])}[\theta_{(\U,\Y)}] , \theta_{U} \right) + \bar{g}_{\Y}^{r} \left( \rmD \overline{ \grad_\Y \, h_{r}([\U,\Y])}[\theta_{(\U, \Y)}], \theta_Y \right)\\
		 & \quad +  \tr\left( \rmD \V_\Y\left[\overline{ \grad_\Y \, h_{r}([\U,\Y])}\right] \theta_U^\top \theta_U  \right)/2  +  \tr\left( \rmD \W_\Y\left[\overline{ \grad_\Y \, h_{r}([\U,\Y])}\right] \theta_Y^\top \theta_Y  \right)/2. 
	\end{split}
\end{equation}

Next, we compute $ \rmD \overline{ \grad_\U \, h_{r}([\U,\Y])}[\theta_{(\U,\Y)}]$ and $ \rmD \overline{ \grad_\Y \, h_{r}([\U,\Y])}[\theta_{(\U, \Y)}]$ separately. We have
\begin{equation} \label{eq: directional-derivative-UY-Upart}
	\begin{split}
		&\rmD \overline{ \grad_\U \, h_{r}([\U,\Y])}[\theta_{(\U,\Y)}] \\
		=& \lim_{t \to 0} \left( (\I_{p_1} - P_{\U + t\theta_U}) \nabla f\left( (\U + t\theta_U) (\Y + t \theta_Y )^\top \right) (\Y+ t\theta_Y) \V^{-1}_{\Y + t \theta_Y } - P_{\U_\perp} \nabla f(\U \Y^\top) \Y \V^{-1}_\Y \right)/t\\
		=&  P_{\U_\perp}\left( \nabla f(\U \Y^\top) \theta_Y \V^{-1}_\Y +  \nabla f(\U \Y^\top) \Y \rmD \V_\Y^{-1}[\theta_Y]  + \nabla^2 f(\U \Y^\top)[\U \theta_Y^\top + \theta_U \Y^\top] \Y \V^{-1}_\Y  \right)\\
		& - (\U \theta_U^\top  + \theta_U \U^\top) \nabla f(\U \Y^\top) \Y\V^{-1}_\Y,
	\end{split}
\end{equation} and
\begin{equation}\label{eq: directional-derivative-UY-Ypart}
	\begin{split}
		 &\rmD \overline{ \grad_\Y \, h_{r}([\U,\Y])}[\theta_{(\U, \Y)}] \\
		 = & \lim_{t \to 0} \left( (\nabla f\left((\U + t\theta_U) (\Y + t \theta_Y )^\top\right))^\top (\U + t\theta_U) \W_{(\Y + t \theta_Y )}^{-1}  - (\nabla f(\U \Y^\top))^\top \U \W_\Y^{-1} \right)/t\\
		 = & (\nabla f(\U \Y^\top))^\top \theta_U \W_\Y^{-1} + (\nabla f(\U \Y^\top))^\top \U \rmD \W_\Y^{-1} [\theta_Y] + \left( \nabla^2 f(\U \Y^\top)[\U \theta_Y^\top + \theta_U \Y^\top] \right)^\top \U \W_\Y^{-1}.
	\end{split}
\end{equation}
By plugging \eqref{eq: directional-derivative-UY-Upart} and \eqref{eq: directional-derivative-UY-Ypart} into \eqref{eq: koszul-formula-UY} and observing the fact $\theta_{(\U,\Y)} \in \cH_{(\U,\Y)} \widebar{\cM}_r^{q_3}$, we have
\begin{equation*}
	\begin{split}
		&\overline{\Hess \, h_{r}([\U,\Y])}[\theta_{(\U,\Y)}, \theta_{(\U,\Y)}] \\
		=& \nabla^2 f(\U \Y^\top)[\U \theta_Y^\top + \theta_U \Y^\top, \U \theta_Y^\top + \theta_U \Y^\top] + 2 \langle  \nabla f(\U \Y^\top), \theta_U \theta_Y^\top \rangle - \langle  \U^\top \nabla f(\U \Y^\top) \Y, \theta_U^\top \theta_U \rangle \\
		& + \langle (\nabla f(\U \Y^\top))^\top \U \rmD \W_\Y^{-1} [\theta_Y], \theta_Y \W_\Y \rangle + \langle \nabla f(\U \Y^\top) \Y \rmD \V_\Y^{-1} [\theta_Y], \theta_U \V_\Y \rangle\\
		& + \tr\left( \rmD \V_\Y\left[\overline{ \grad_\Y \, h_{r}([\U,\Y])}\right] \theta_U^\top \theta_U  \right)/2 + \langle \rmD \W_\Y\left[\overline{ \grad_\Y \, h_{r}([\U,\Y])}\right], \theta_Y^\top \theta_Y  \rangle /2.
	\end{split}
\end{equation*} This finishes the proof of this proposition. \quad $\blacksquare$

\subsection{Proof of Theorem \ref{th: embedded-quotient-connection-general3}.}
First, recall $\X = \U \Y^\top$ and $\V$ spans the right singular subspace of $\X$, so $\Y$ lies in the column space of $\V$ and $\Y \Y^\dagger = P_\V$. Thus, \eqref{eq: gradient-connect-general3} is by direct calculation from the gradient expressions in Proposition \ref{prop: gradient-hessian-exp-general}. 

Next, we prove \eqref{eq: Hessian-connection-general3}. Since $[\U,\Y]$ is a Riemannian FOSP of \eqref{eq: general-opt-problem-quotient-sub3}, we have
\begin{equation} \label{eq: FOSP-condition-general3}
	\overline{ \grad \, h_{r}([\U,\Y])} = \0 \quad, \quad    (\nabla f(\U \Y^\top))^\top \U = \0 \quad,  \text{ and } \quad  \nabla f(\U \Y^\top) \Y = \0
\end{equation} and have $\grad f(\X) = \0$ by \eqref{eq: gradient-connect-general3}. So $ \nabla f(\X) = P_{\U_\perp} \nabla f(\X) P_{\V_\perp}$. Let $\bSigma = \U^\top \X \V$. Given any $\theta_{(\U,\Y)} =  [\theta_U^\top \quad \theta_Y^\top]^\top \in \cH_{(\U,\Y)} \widebar{\cM}_{r}^{q_3}$, we have
\begin{equation} \label{eq: Hessian-con-gradient-general3}
	\langle \nabla f(\X), P_{\U_\perp} \theta_U \Y^\top \V \bSigma^{-1}  \theta_Y^\top  P_{\V_\perp}    \rangle = \langle \nabla f(\X), \theta_U \theta_Y^\top \rangle,
\end{equation} where the equality is because $\Y^\top \V = \U^\top \X \V = \bSigma$ and $ \nabla f(\X) = P_{\U_\perp} \nabla f(\X) P_{\V_\perp} $.

Then by Proposition \ref{prop: gradient-hessian-exp-general}:
\begin{equation*}
	\begin{split}
		 &\quad \overline{\Hess \, h_{r}([\U,\Y])}[\theta_{(\U,\Y)}, \theta_{(\U,\Y)}] \\
		=& \nabla^2 f(\U \Y^\top)[\U \theta_Y^\top + \theta_U \Y^\top, \U \theta_Y^\top + \theta_U \Y^\top] + 2 \langle  \nabla f(\U \Y^\top), \theta_U \theta_Y^\top \rangle - \langle  \U^\top \nabla f(\U \Y^\top) \Y, \theta_U^\top \theta_U \rangle \\
		& + \langle (\nabla f(\U \Y^\top))^\top \U \rmD \W_\Y^{-1} [\theta_Y], \theta_Y \W_\Y \rangle + \langle \nabla f(\U \Y^\top) \Y \rmD \V_\Y^{-1} [\theta_Y], \theta_U \V_\Y \rangle\\
		& + \langle \rmD \W_\Y\left[\overline{ \grad_\Y \, h_{r}([\U,\Y])}\right], \theta_Y^\top \theta_Y  \rangle /2 + \langle \rmD \V_\Y\left[\overline{ \grad_\Y \, h_{r}([\U,\Y])}\right], \theta_U^\top \theta_U  \rangle /2\\
		 & \overset{ \eqref{eq: FOSP-condition-general3} } = \nabla^2 f(\U \Y^\top)[\U \theta_Y^\top + \theta_U \Y^\top, \U \theta_Y^\top + \theta_U \Y^\top] + 2 \langle  \nabla f(\U \Y^\top), \theta_U \theta_Y^\top \rangle \\
		 & \overset{ \text{Proposition } \ref{prop: general-bijection3}, \eqref{eq: Hessian-con-gradient-general3} } =   \nabla^2 f(\X)[\cL^r_{\U,\Y}(\theta_{(\U,\Y)}),\cL^r_{\U,\Y}(\theta_{(\U,\Y)})] + 2 \langle \nabla f(\X), P_{\U_\perp} \theta_U \Y^\top \V \bSigma^{-1}  \theta_Y^\top  P_{\V_\perp}    \rangle\\
		 & = \Hess f(\X)[\cL^r_{\U,\Y}(\theta_{(\U,\Y)}),\cL^r_{\U,\Y}(\theta_{(\U,\Y)})],
	\end{split}
\end{equation*} where the last equality follows from the expression of $\Hess f(\X)$ in \eqref{eq: embedded-gd-hessian-general} and the definition of $\cL^r_{\U,\Y}$. 

Then, by \eqref{ineq: bijection-spectrum-general3}, \eqref{eq: Hessian-connection-general3} and Theorem \ref{th: hessian-sandwich}, we have $\overline{\Hess \, h_{r}([\U,\Y])}$ has $(p_1+p_2- r)r$ eigenvalues and $\widebar{\lambda}_i(\Hess \, h_{r}([\U,\Y]))$ is sandwiched between $(\sigma_r(\W^{-1}_\Y) \wedge \sigma_r^2(\Y \V_\Y^{-1/2} )) \lambda_i(\Hess f(\X)) $ and $( \sigma_1(\W_\Y^{-1}) \vee \sigma^2_1(\Y\V_\Y^{-1/2})   ) \lambda_i(\Hess f(\X)) $ for $i = 1,\ldots,(p_1+p_2- r)r$.\quad $\blacksquare$

\subsection{Proof of Corollary \ref{coro: landscape connection general case}.}
Here we prove the Riemannian FOSP, SOSP, and strict saddle equivalence of \eqref{eq: general prob} on $\cM_{r}^e$ and $\cM_{r}^{q_1}$, similar proof applies to the equivalence on $\cM_{r}^e$ and $\cM_{r}^{q_2}$ (or $\cM_{r}^{q_3}$). First, by the connection of Riemannian gradients in \eqref{eq: gradient-connect-general1}, the connection of Riemannian FOSPs under two geometries clearly holds. 

Suppose $[\L,\R]$ is a Riemannian SOSP of \eqref{eq: general-opt-problem-quotient-sub1} and let $\X = \L \R^\top$. Given any  $\xi_\X \in T_{\X}\cM^e_{r}$, we have $\Hess f(\X)[\xi_\X, \xi_\X] \overset{ \eqref{eq: Hessian-connection-general1} }=\overline{\Hess \, h_{r}([\L,\R])}[(\cL^{r}_{\L,\R})^{-1}(\xi_\X), (\cL^{r}_{\L,\R})^{-1}(\xi_\X)] \geq 0$, where the inequality is by the SOSP assumption on $[\L,\R]$. Combining the fact $\X$ is a Riemannian FOSP, this shows $\X = \L \R^\top$ is a Riemannian SOSP under the embedded geometry.

Next, let us show the other direction: suppose $\X$ is a Riemannian SOSP under the embedded geometry, then there is a unique $[\L,\R]$ such that $\L \R^\top = \X$ and it is a Riemannian SOSP of \eqref{eq: general-opt-problem-quotient-sub1}. To see this, first the uniqueness of $[\L,\R]$ is guaranteed by the fact $\bar{\ell}: (\L,\R) \in \widebar{\cM}_r^{q_1} \to \L \R^\top \in \cM_r^e$ induces a diffeomorphism between $\cM_r^{q_1}$ and $\cM_r^e$. In addition, we have shown $[\L,\R]$ is a Riemannian FOSP of \eqref{eq: general-opt-problem-quotient-sub1}. Then by \eqref{eq: Hessian-connection-general1}, we have for any $\theta_{(\L,\R)} \in \cH_{(\L,\R)} \widebar{\cM}_{r}^{q_1}$, $\overline{\Hess \, h_{r}([\L,\R])}[\theta_{(\L,\R)}, \theta_{(\L,\R)}]  = \Hess f(\X)[\cL_{\L,\R}^{r}(\theta_{(\L,\R)}),\cL_{\L,\R}^{r}(\theta_{(\L,\R)})] \geq 0$. 

Finally, the strict saddle equivalence follows from the sandwich inequality in Theorem \ref{th: embedded-quotient-connection-general1} and the definition of the strict saddle.
\quad $\blacksquare$

\section{Additional Lemmas} \label{sec: additional-lemmas}
The following lemma provides the connection of horizontal lifts of tangent vectors at different reference points for the quotient geometries we consider.
\begin{lemma} \label{lm: horizontal-lifts-connection}
	\begin{itemize}
		\item[(i)] Suppose $\theta_{[\Y]}$ is a tangent vector to the quotient manifold $\cM_{r+}^{q_1}$ at $[\Y]$. If $\W_\Y = \O \W_{\Y\O} \O^\top$ holds for any $\O \in \bbO_r$, then the horizontal lifts of $\theta_{[\Y]}$ at $\Y$ and $\Y\O$ for $\O \in \bbO_{r}$ are related as $\theta_{\Y \O} = \theta_{\Y} \O$.
		\item[(ii)] Suppose $\theta_{[\U, \B]}$ is a tangent vector to the quotient manifold $\cM_{r+}^{q_2}$ at $[\U,\B]$. The horizontal lifts of $\theta_{[\U,\B]}$ at $(\U, \B)$ and at $(\U \O, \O^\top \B \O)$ are related as $\theta_{(\U \O, \O^\top \B \O)} =[(\theta_U \O)^\top \quad (\O^\top \theta_{B} \O)^\top ]^\top  $ given $\theta_{(\U,\B)} = [\theta_U^\top \quad \theta_B^\top]^\top$.
		\item[(iii)] Suppose $\theta_{[\L, \R]}$ is a tangent vector to the quotient manifold $\cM_{r}^{q_1}$ at $[\L,\R]$. If $\W_{\L,\R} = \M \W_{\L\M,\R \M^{-\top}} \M^\top$ and $\V_{\L,\R} = \M^{-\top } \V_{\L\M,\R \M^{-\top}} \M^{-1}$ hold for any $\M \in \GL(r)$, then the horizontal lifts of $\theta_{[\L,\R]}$ at $(\L, \R)$ and at $(\L \M, \R \M^{-\top})$ are related as $\theta_{(\L \M, \R \M^{-\top})} =[(\theta_L \M)^\top \quad (\theta_{R} \M^{-\top})^\top ]^\top  $ given $\theta_{(\L,\R)} = [\theta_L^\top \quad \theta_R^\top]^\top$.
		\item[(iv)] Suppose $\theta_{[\U, \B, \V]}$ is a tangent vector to the quotient manifold $\cM_{r}^{q_2}$ at $[\U,\B,\V]$. The horizontal lifts of $\theta_{[\U,\B,\V]}$ at $(\U, \B, \V)$ and at $(\U \O, \O^\top \B \O, \V\O)$ are related as $\theta_{(\U \O, \O^\top \B \O, \V \O)} =[(\theta_U \O)^\top \quad (\O^\top \theta_{B} \O)^\top \quad (\theta_V \O)^\top ]^\top  $ given $\theta_{(\U,\B,\V)} = [\theta_U^\top \quad \theta_B^\top \quad \theta_V^\top]^\top$.
		\item[(v)] Suppose $\theta_{[\U, \Y]}$ is a tangent vector to the quotient manifold $\cM_{r}^{q_3}$ at $[\U,\Y]$. The horizontal lifts of $\theta_{[\U,\Y]}$ at $(\U, \Y)$ and at $(\U \O, \Y \O)$ are related as $\theta_{(\U \O, \Y \O)} =[(\theta_U \O)^\top \quad (\theta_{Y} \O)^\top ]^\top  $ given $\theta_{(\U,\Y)} = [\theta_U^\top \quad \theta_Y^\top]^\top$.
	\end{itemize}
\end{lemma}
{\noindent \bf Proof of Lemma \ref{lm: horizontal-lifts-connection}.} Since the proofs for claims (ii)(iii)(iv)(v) are similar to the proof of the first claim, for simplicity, we only present the proof for the first one. Let $l: \cM_{r+}^{q_1} \to \bbR$ be an arbitrary smooth function and define $\bar{l}:= l\circ \pi: \widebar{\cM}^{q_1 }_{r+} \to \bbR$, where $\pi$ is the quotient mapping on the manifold $\cM_{r+}^{q_1}$. Let $z: \Y \to \Y \O$, so we have
\begin{equation} \label{eq: horizonal-lift-connect-eq1}
	\bar{l}(\Y) = \bar{l}(z(\Y)), \quad \forall \Y \in \bbR^{p \times r}_*.
\end{equation}
By taking differential with respect to $\Y$ along direction $\theta_\Y \in \cH_\Y \widebar{\cM}_{r+}^{q_1}$ on both sides of \eqref{eq: horizonal-lift-connect-eq1}, we have
\begin{equation} \label{eq: horizonal-lift-connect-eq2}
	\rmD \bar{l}(\Y)[\theta_\Y] = \rmD \bar{l}(z(\Y))[ \rmD z(\Y) [\theta_\Y] ].
\end{equation} Moreover, by chain rule and the definition of the horizontal lift of a tangent vector, we have
\begin{equation} \label{eq: horizonal-lift-connect-eq3}
	\rmD \bar{l}(\Y) [\theta_\Y] = \rmD l(\pi(\Y))[ \rmD \pi(\Y)[\theta_\Y]  ]  = \rmD l(\pi(\Y)) [ \theta_{[\Y]} ].
\end{equation} In addition,
\begin{equation} \label{eq: horizonal-lift-connect-eq4}
	\rmD z(\Y)[\theta_\Y] = \lim_{t \to 0} (z(\Y + t\theta_\Y) - z(\Y))/t = \theta_\Y \O.
\end{equation}
Then,
\begin{equation*}
\begin{split}
	\rmD l(\pi(\Y\O)) [ \theta_{[\Y]} ] = \rmD l(\pi(\Y)) [ \theta_{[\Y]} ] \overset{\eqref{eq: horizonal-lift-connect-eq3}}= \rmD \bar{l}(\Y) [\theta_\Y]  &\overset{\eqref{eq: horizonal-lift-connect-eq2}}=  \rmD \bar{l}(z(\Y))[ \rmD z(\Y) [\theta_\Y]]\\
	&\overset{\eqref{eq: horizonal-lift-connect-eq4}}=\rmD \bar{l}(z(\Y))[ \theta_\Y \O ] \\
	&\overset{(a)}= \rmD l(\pi(\Y\O)) [ \rmD \pi(\Y \O)[\theta_{\Y} \O] ].
\end{split}
\end{equation*} Here (a) is by chain rule and the definition of $z$. Since the above equation holds for any $l$, this implies $\rmD \pi(\Y \O)[\theta_{\Y} \O]= \theta_{[\Y]}$. Finally, by Lemma \ref{lm: psd-quotient-manifold1-prop}(i), we have $\cH_\Y \widebar{\cM}_{r+}^{q_1} = \{ \theta_\Y': \theta_\Y' = (\U \S + \U_\perp \D) \P^{-\top}, \D \in \bbR^{(p-r) \times r}, \S\P^{-\top}\W_\Y \P^{-1} \in \bbS^{r \times r} \}$, where $\U \in \st(r,p)$ spans the top $r$ eigenspace of $\Y\Y^\top$ and $\P = \U^\top \Y$. Then, it holds that $\cH_{\Y\O} \widebar{\cM}_{r+}^{q_1} = \{ \theta_{\Y\O}': \theta_{\Y\O}' = (\U \S + \U_\perp \D) \P^{-\top} \O, \D \in \bbR^{(p-r) \times r}, \S\P^{-\top}\W_\Y \P^{-1} \in \bbS^{r \times r} \}$. Thus, we have $\theta_\Y \O \in \cH_{\Y\O} \widebar{\cM}_{r+}^{q_1}$ given $\theta_\Y \in \cH_{\Y} \widebar{\cM}_{r+}^{q_1}$ and conclude $\theta_\Y \O$ is the unique horizontal lift of $\theta_{[\Y]}$ at $\Y \O$.
\quad $\blacksquare$

\begin{lemma} \label{lm: quotient-hessian-riemannian-connection-equal} Suppose $f: \cM \to \bbR$ is a smooth function defined on the Riemannian quotient manifold $\cM = \widebar{\cM}/\sim$ and $\widebar{\cM}$ is endowed with the Riemannian metric $\bar{g}_\X$, which induces a Riemannian metric $g_{[\X]}$ on $\cM$. Then $\overline{\Hess f([\X])}[\theta_{[\X]}, \eta_{[\X]}] = \bar{g}_\X\left( \widebar{\nabla}_{\theta_\X} \overline{\grad f }, \eta_\X \right)$ for any $\theta_{[\X]}, \eta_{[\X]} \in T_{[\X]} \cM$.
\end{lemma}
{\noindent \bf Proof of Lemma \ref{lm: quotient-hessian-riemannian-connection-equal}.} Suppose the vertical and horizontal spaces on $\widebar{\cM}$ are $\cV_\X \widebar{\cM}$ and $\cH_\X \widebar{\cM}$ with $\cV_\X \widebar{\cM} \oplus \cH_\X \widebar{\cM} = T_\X \widebar{\cM}$. Let us define $\cV'_\X \widebar{\cM}$ to be the subspace in $T_\X \widebar{\cM}$ that is orthogonal to $\cH_\X \widebar{\cM}$ with respect to $\bar{g}_\X$. In general, $\cV'_\X \widebar{\cM}$ is not equal to $\cV_\X \widebar{\cM}$ unless we pick the horizontal space in a canonical way.

Since $\widebar{\nabla}_{\theta_\X} \overline{\grad f }$ by definition belongs to $T_\X \widebar{\cM}$, we have for $\theta_\X, \eta_\X \in \cH_\X \widebar{\cM}$,
\begin{equation*}
	\begin{split}
		\overline{\Hess f([\X])}[\theta_{[\X]}, \eta_{[\X]}] &= \bar{g}_\X\left( \overline{\Hess f([\X])[\theta_{[\X]}]}, \eta_\X \right) \\
		& \overset{\eqref{eq: quotient-hessian-linear-from} }=\bar{g}_\X\left( P_\X^{\cH}\left( \widebar{\nabla}_{\theta_\X} \overline{\grad f} \right),\eta_\X \right) \\
		& = \bar{g}_\X\left( \widebar{\nabla}_{\theta_\X} \overline{\grad f}- P_\X^{\cV'}\left( \widebar{\nabla}_{\theta_\X} \overline{\grad f} \right),\eta_\X \right)\\
		& \overset{(a)}= \bar{g}_\X\left( \widebar{\nabla}_{\theta_\X} \overline{\grad f }, \eta_\X \right),
	\end{split}
\end{equation*} where (a) is because $\eta_\X \in \cH_\X \widebar{\cM}$ and $P_\X^{\cH}$ and $P_\X^{\cV'}$ denote the projection operators onto $\cH_\X \widebar{\cM}$ and $\cV'_\X \widebar{\cM}$, respectively. \quad $\blacksquare$

\begin{lemma}{\rm(Max-min Theorem for Eigenvalues \cite[Corollary III.1.2]{bhatia2013matrix} )}\label{lm: max-min-theorem} Suppose $\A$ is a Hermitian operator on the Hilbert space $\cH$ with dimension $p$ and inner product $g(\cdot,\cdot)$ and $\A$ has eigenvalues $\lambda_1 \geq \lambda_2 \geq \cdots \geq \lambda_p$. If $\scC_k$ denotes the set of $k$-dimensional subspaces of $\cH$, then $\lambda_k = \max_{C \in \scC_k} \min_{\u \in C, \u \neq 0 } g(\u, \A \u)/g(\u,\u) .$
\end{lemma}

\bibliographystyle{alpha}
\bibliography{reference.bib}

\newcommand{\etalchar}[1]{$^{#1}$}
\begin{thebibliography}{AIDLVH09}

\bibitem[AAM14]{absil2014two}
P-A Absil, Luca Amodei, and Gilles Meyer.
\newblock Two {N}ewton methods on the manifold of fixed-rank matrices endowed
  with {R}iemannian quotient geometries.
\newblock {\em Computational Statistics}, 29(3):569--590, 2014.

\bibitem[AIDLVH09]{absil2009geometric}
P-A Absil, Mariya Ishteva, Lieven De~Lathauwer, and Sabine Van~Huffel.
\newblock A geometric {N}ewton method for {O}ja's vector field.
\newblock {\em Neural computation}, 21(5):1415--1433, 2009.

\bibitem[AMR12]{abraham2012manifolds}
Ralph Abraham, Jerrold~E Marsden, and Tudor Ratiu.
\newblock {\em Manifolds, tensor analysis, and applications}, volume~75.
\newblock Springer Science \& Business Media, 2012.

\bibitem[AMS09]{absil2009optimization}
P-A Absil, Robert Mahony, and Rodolphe Sepulchre.
\newblock {\em Optimization algorithms on matrix manifolds}.
\newblock Princeton University Press, 2009.

\bibitem[AMT13]{absil2013extrinsic}
P-A Absil, Robert Mahony, and Jochen Trumpf.
\newblock An extrinsic look at the {R}iemannian {H}essian.
\newblock In {\em International Conference on Geometric Science of
  Information}, pages 361--368. Springer, 2013.

\bibitem[AS21]{ahn2021riemannian}
Kwangjun Ahn and Felipe Suarez.
\newblock {R}iemannian perspective on matrix factorization.
\newblock {\em arXiv preprint arXiv:2102.00937}, 2021.

\bibitem[BA11]{boumal2011rtrmc}
Nicolas Boumal and P-A Absil.
\newblock Rtrmc: A {R}iemannian trust-region method for low-rank matrix
  completion.
\newblock In {\em Advances in neural information processing systems}, pages
  406--414, 2011.

\bibitem[BAC19]{boumal2019global}
Nicolas Boumal, P-A Absil, and Coralia Cartis.
\newblock Global rates of convergence for nonconvex optimization on manifolds.
\newblock {\em IMA Journal of Numerical Analysis}, 39(1):1--33, 2019.

\bibitem[Bha13]{bhatia2013matrix}
Rajendra Bhatia.
\newblock {\em Matrix analysis}, volume 169.
\newblock Springer Science \& Business Media, 2013.

\bibitem[BKS16]{bhojanapalli2016dropping}
Srinadh Bhojanapalli, Anastasios Kyrillidis, and Sujay Sanghavi.
\newblock Dropping convexity for faster semi-definite optimization.
\newblock In {\em Conference on Learning Theory}, pages 530--582, 2016.

\bibitem[Bou20]{boumal2020introduction}
Nicolas Boumal.
\newblock An introduction to optimization on smooth manifolds.
\newblock {\em http://sma.epfl.ch/~nboumal/\#book}, 2020.

\bibitem[BS10]{bonnabel2010riemannian}
Silv{\`e}re Bonnabel and Rodolphe Sepulchre.
\newblock Riemannian metric and geometric mean for positive semidefinite
  matrices of fixed rank.
\newblock {\em SIAM Journal on Matrix Analysis and Applications},
  31(3):1055--1070, 2010.

\bibitem[CB19]{criscitiello2019efficiently}
Chris Criscitiello and Nicolas Boumal.
\newblock Efficiently escaping saddle points on manifolds.
\newblock {\em Advances in Neural Information Processing Systems},
  32:5987--5997, 2019.

\bibitem[CC18]{chen2018harnessing}
Yudong Chen and Yuejie Chi.
\newblock Harnessing structures in big data via guaranteed low-rank matrix
  estimation: Recent theory and fast algorithms via convex and nonconvex
  optimization.
\newblock {\em IEEE Signal Processing Magazine}, 35(4):14--31, 2018.

\bibitem[CLC19]{chi2019nonconvex}
Yuejie Chi, Yue~M Lu, and Yuxin Chen.
\newblock Nonconvex optimization meets low-rank matrix factorization: An
  overview.
\newblock {\em IEEE Transactions on Signal Processing}, 67(20):5239--5269,
  2019.

\bibitem[CW18]{cai2018exploiting}
Jian-Feng Cai and Ke~Wei.
\newblock Exploiting the structure effectively and efficiently in low-rank
  matrix recovery.
\newblock In {\em Handbook of Numerical Analysis}, volume~19, pages 21--51.
  Elsevier, 2018.

\bibitem[CZ13]{cai2013compressed}
T~Tony Cai and Anru Zhang.
\newblock Compressed sensing and affine rank minimization under restricted
  isometry.
\newblock {\em IEEE Transactions on Signal Processing}, 61(13):3279--3290,
  2013.

\bibitem[DH21]{douik2021low}
Ahmed Douik and Babak Hassibi.
\newblock Low-rank {R}iemannian optimization for graph-based clustering
  applications.
\newblock {\em IEEE Transactions on Pattern Analysis and Machine Intelligence},
  2021.

\bibitem[EAS98]{edelman1998geometry}
Alan Edelman, Tom{\'a}s~A Arias, and Steven~T Smith.
\newblock The geometry of algorithms with orthogonality constraints.
\newblock {\em SIAM journal on Matrix Analysis and Applications},
  20(2):303--353, 1998.

\bibitem[GHJY15]{ge2015escaping}
Rong Ge, Furong Huang, Chi Jin, and Yang Yuan.
\newblock Escaping from saddle points-online stochastic gradient for tensor
  decomposition.
\newblock In {\em Conference on Learning Theory}, pages 797--842. PMLR, 2015.

\bibitem[GM17]{ge2017optimization}
Rong Ge and Tengyu Ma.
\newblock On the optimization landscape of tensor decompositions.
\newblock In {\em Advances in Neural Information Processing Systems}, pages
  3653--3663, 2017.

\bibitem[GP07]{grubivsic2007efficient}
Igor Grubisi{\'c} and Raoul Pietersz.
\newblock Efficient rank reduction of correlation matrices.
\newblock {\em Linear algebra and its applications}, 422(2-3):629--653, 2007.

\bibitem[GS10]{gao2010majorized}
Yan Gao and Defeng Sun.
\newblock A majorized penalty approach for calibrating rank constrained
  correlation matrix problems.
\newblock {\em https://www.polyu.edu.hk/ama/profile/dfsun/MajorPen.pdf}, 2010.

\bibitem[HGZ17]{huang2017solving}
Wen Huang, Kyle~A Gallivan, and Xiangxiong Zhang.
\newblock Solving phaselift by low-rank riemannian optimization methods for
  complex semidefinite constraints.
\newblock {\em SIAM Journal on Scientific Computing}, 39(5):B840--B859, 2017.

\bibitem[HH18]{huang2018blind}
Wen Huang and Paul Hand.
\newblock Blind deconvolution by a steepest descent algorithm on a quotient
  manifold.
\newblock {\em SIAM Journal on Imaging Sciences}, 11(4):2757--2785, 2018.

\bibitem[HLB20]{ha2020equivalence}
Wooseok Ha, Haoyang Liu, and Rina~Foygel Barber.
\newblock An equivalence between critical points for rank constraints versus
  low-rank factorizations.
\newblock {\em SIAM Journal on Optimization}, 30(4):2927--2955, 2020.

\bibitem[HLWY20]{hu2020brief}
Jiang Hu, Xin Liu, Zai-Wen Wen, and Ya-Xiang Yuan.
\newblock A brief introduction to manifold optimization.
\newblock {\em Journal of the Operations Research Society of China},
  8(2):199--248, 2020.

\bibitem[HLZ20]{hou2020fast}
Thomas~Y Hou, Zhenzhen Li, and Ziyun Zhang.
\newblock Fast global convergence for low-rank matrix recovery via {R}iemannian
  gradient descent with random initialization.
\newblock {\em arXiv preprint arXiv:2012.15467}, 2020.

\bibitem[HM12]{helmke2012optimization}
Uwe Helmke and John~B Moore.
\newblock {\em Optimization and dynamical systems}.
\newblock Springer Science \& Business Media, 2012.

\bibitem[JBAS10]{journee2010low}
Michel Journ{\'e}e, Francis Bach, P-A Absil, and Rodolphe Sepulchre.
\newblock Low-rank optimization on the cone of positive semidefinite matrices.
\newblock {\em SIAM Journal on Optimization}, 20(5):2327--2351, 2010.

\bibitem[JMD10]{jain2010guaranteed}
Prateek Jain, Raghu Meka, and Inderjit~S Dhillon.
\newblock Guaranteed rank minimization via singular value projection.
\newblock In {\em Advances in Neural Information Processing Systems}, pages
  937--945, 2010.

\bibitem[JNS13]{jain2013low}
Prateek Jain, Praneeth Netrapalli, and Sujay Sanghavi.
\newblock Low-rank matrix completion using alternating minimization.
\newblock In {\em Proceedings of the 45th Annual ACM Symposium on Theory of
  Computing}, pages 665--674, 2013.

\bibitem[Lee13]{lee2013smooth}
John~M Lee.
\newblock {\em Introduction to Smooth Manifolds}.
\newblock Springer, 2013.

\bibitem[LHLZ20]{luo2020recursive}
Yuetian Luo, Wen Huang, Xudong Li, and Anru~R Zhang.
\newblock Recursive importance sketching for rank constrained least squares:
  Algorithms and high-order convergence.
\newblock {\em arXiv preprint arXiv:2011.08360}, 2020.

\bibitem[LKB22]{levin2022effect}
Eitan Levin, Joe Kileel, and Nicolas Boumal.
\newblock The effect of smooth parametrizations on nonconvex optimization
  landscapes.
\newblock {\em arXiv preprint arXiv:2207.03512}, 2022.

\bibitem[LLA{\etalchar{+}}19]{li2016symmetry}
Xingguo Li, Junwei Lu, Raman Arora, Jarvis Haupt, Han Liu, Zhaoran Wang, and
  Tuo Zhao.
\newblock Symmetry, saddle points, and global optimization landscape of
  nonconvex matrix factorization.
\newblock {\em IEEE Transactions on Information Theory}, 65(6):3489--3514,
  2019.

\bibitem[LLZ21]{luo2021nonconvex}
Yuetian Luo, Xudong Li, and Anru~R Zhang.
\newblock Nonconvex factorization and manifold formulations are almost
  equivalent in low-rank matrix optimization.
\newblock {\em arXiv preprint arXiv:2108.01772}, 2021.

\bibitem[LPP{\etalchar{+}}19]{lee2019first}
Jason~D Lee, Ioannis Panageas, Georgios Piliouras, Max Simchowitz, Michael~I
  Jordan, and Benjamin Recht.
\newblock First-order methods almost always avoid strict saddle points.
\newblock {\em Mathematical programming}, 176(1-2):311--337, 2019.

\bibitem[MA20]{massart2020quotient}
Estelle Massart and P-A Absil.
\newblock Quotient geometry with simple geodesics for the manifold of
  fixed-rank positive-semidefinite matrices.
\newblock {\em SIAM Journal on Matrix Analysis and Applications},
  41(1):171--198, 2020.

\bibitem[MAS12]{mishra2012riemannian}
Bamdev Mishra, K~Adithya Apuroop, and Rodolphe Sepulchre.
\newblock A {R}iemannian geometry for low-rank matrix completion.
\newblock {\em arXiv preprint arXiv:1211.1550}, 2012.

\bibitem[MBS11a]{meyer2011linear}
Gilles Meyer, Silv{\`e}re Bonnabel, and Rodolphe Sepulchre.
\newblock Linear regression under fixed-rank constraints: a {R}iemannian
  approach.
\newblock In {\em Proceedings of the 28th International Conference on Machine
  Learning}, 2011.

\bibitem[MBS11b]{meyer2011regression}
Gilles Meyer, Silv{\`e}re Bonnabel, and Rodolphe Sepulchre.
\newblock Regression on fixed-rank positive semidefinite matrices: A
  {R}iemannian approach.
\newblock {\em The Journal of Machine Learning Research}, 12:593--625, 2011.

\bibitem[Mey11]{meyer2011geometric}
Gilles Meyer.
\newblock Geometric optimization algorithms for linear regression on fixed-rank
  matrices.
\newblock {\em Ph.D. thesis, Universit{\'e} de Li{\`e}ge, Belgique}, 2011.

\bibitem[MMBS13]{mishra2013low}
Bamdev Mishra, Gilles Meyer, Francis Bach, and Rodolphe Sepulchre.
\newblock Low-rank optimization with trace norm penalty.
\newblock {\em SIAM Journal on Optimization}, 23(4):2124--2149, 2013.

\bibitem[MMBS14]{mishra2014fixed}
Bamdev Mishra, Gilles Meyer, Silv{\`e}re Bonnabel, and Rodolphe Sepulchre.
\newblock Fixed-rank matrix factorizations and {R}iemannian low-rank
  optimization.
\newblock {\em Computational Statistics}, 29(3-4):591--621, 2014.

\bibitem[MMS11]{mishra2011low}
Bamdev Mishra, Gilles Meyer, and Rodolphe Sepulchre.
\newblock Low-rank optimization for distance matrix completion.
\newblock In {\em 2011 50th IEEE Conference on Decision and Control and
  European Control Conference}, pages 4455--4460. IEEE, 2011.

\bibitem[MS16]{mishra2016riemannian}
Bamdev Mishra and Rodolphe Sepulchre.
\newblock {R}iemannian preconditioning.
\newblock {\em SIAM Journal on Optimization}, 26(1):635--660, 2016.

\bibitem[MWCC19]{ma2019implicit}
Cong Ma, Kaizheng Wang, Yuejie Chi, and Yuxin Chen.
\newblock Implicit regularization in nonconvex statistical estimation: Gradient
  descent converges linearly for phase retrieval, matrix completion, and blind
  deconvolution.
\newblock {\em Foundations of Computational Mathematics}, pages 1--182, 2019.

\bibitem[MZL19]{maunu2019well}
Tyler Maunu, Teng Zhang, and Gilad Lerman.
\newblock A well-tempered landscape for non-convex robust subspace recovery.
\newblock {\em Journal of Machine Learning Research}, 20:1--59, 2019.

\bibitem[NS12]{ngo2012scaled}
Thanh Ngo and Yousef Saad.
\newblock Scaled gradients on {G}rassmann manifolds for matrix completion.
\newblock {\em Advances in neural information processing systems},
  25:1412--1420, 2012.

\bibitem[Pet06]{petersen2006riemannian}
Peter Petersen.
\newblock {\em Riemannian geometry}, volume 171.
\newblock Springer, 2006.

\bibitem[RFP10]{recht2010guaranteed}
Benjamin Recht, Maryam Fazel, and Pablo~A Parrilo.
\newblock Guaranteed minimum-rank solutions of linear matrix equations via
  nuclear norm minimization.
\newblock {\em SIAM review}, 52(3):471--501, 2010.

\bibitem[SFF19]{sun2019escaping}
Yue Sun, Nicolas Flammarion, and Maryam Fazel.
\newblock Escaping from saddle points on {R}iemannian manifolds.
\newblock In {\em Advances in Neural Information Processing Systems},
  volume~32, pages 7276--7286, 2019.

\bibitem[SL15]{sun2015guaranteed}
Ruoyu Sun and Zhi-Quan Luo.
\newblock Guaranteed matrix completion via nonconvex factorization.
\newblock In {\em Foundations of Computer Science (FOCS), 2015 IEEE 56th Annual
  Symposium on}, pages 270--289. IEEE, 2015.

\bibitem[SQW18]{sun2018geometric}
Ju~Sun, Qing Qu, and John Wright.
\newblock A geometric analysis of phase retrieval.
\newblock {\em Foundations of Computational Mathematics}, 18(5):1131--1198,
  2018.

\bibitem[SWC12]{shalit2012online}
Uri Shalit, Daphna Weinshall, and Gal Chechik.
\newblock Online learning in the embedded manifold of low-rank matrices.
\newblock {\em Journal of Machine Learning Research}, 13:429--458, 2012.

\bibitem[TBS{\etalchar{+}}16]{tu2016low}
Stephen Tu, Ross Boczar, Max Simchowitz, Mahdi Soltanolkotabi, and Benjamin
  Recht.
\newblock Low-rank solutions of linear matrix equations via {P}rocrustes flow.
\newblock In {\em International Conference on Machine Learning}, pages
  964--973, 2016.

\bibitem[TMC21]{tong2020accelerating}
Tian Tong, Cong Ma, and Yuejie Chi.
\newblock Accelerating ill-conditioned low-rank matrix estimation via scaled
  gradient descent.
\newblock {\em Journal of Machine Learning Research}, 22(150):1--63, 2021.

\bibitem[UV20]{uschmajew2018critical}
Andr{\'e} Uschmajew and Bart Vandereycken.
\newblock On critical points of quadratic low-rank matrix optimization
  problems.
\newblock {\em IMA Journal of Numerical Analysis}, 40(4):2626--2651, 2020.

\bibitem[Van13]{vandereycken2013low}
Bart Vandereycken.
\newblock Low-rank matrix completion by {R}iemannian optimization.
\newblock {\em SIAM Journal on Optimization}, 23(2):1214--1236, 2013.

\bibitem[VAV13]{vandereycken2013riemannian}
Bart Vandereycken, P-A Absil, and Stefan Vandewalle.
\newblock A {R}iemannian geometry with complete geodesics for the set of
  positive semidefinite matrices of fixed rank.
\newblock {\em IMA Journal of Numerical Analysis}, 33(2):481--514, 2013.

\bibitem[VV10]{vandereycken2010riemannian}
Bart Vandereycken and Stefan Vandewalle.
\newblock A {R}iemannian optimization approach for computing low-rank solutions
  of {L}yapunov equations.
\newblock {\em SIAM Journal on Matrix Analysis and Applications},
  31(5):2553--2579, 2010.

\bibitem[WCCL16]{wei2016guarantees}
Ke~Wei, Jian-Feng Cai, Tony~F Chan, and Shingyu Leung.
\newblock Guarantees of {R}iemannian optimization for low rank matrix recovery.
\newblock {\em SIAM Journal on Matrix Analysis and Applications},
  37(3):1198--1222, 2016.

\bibitem[ZFZ21]{zhang2021preconditioned}
Jialun Zhang, Salar Fattahi, and Richard~Y Zhang.
\newblock Preconditioned gradient descent for over-parameterized nonconvex
  matrix factorization.
\newblock {\em Advances in Neural Information Processing Systems},
  34:5985--5996, 2021.

\bibitem[ZHVZ22]{zheng2022riemannian}
Shixin Zheng, Wen Huang, Bart Vandereycken, and Xiangxiong Zhang.
\newblock Riemannian optimization using three different metrics for {H}ermitian
  {PSD} fixed-rank constraints: an extended version.
\newblock {\em arXiv preprint arXiv:2204.07830}, 2022.

\bibitem[ZKHC21]{zhuo2021computational}
Jiacheng Zhuo, Jeongyeol Kwon, Nhat Ho, and Constantine Caramanis.
\newblock On the computational and statistical complexity of over-parameterized
  matrix sensing.
\newblock {\em arXiv preprint arXiv:2102.02756}, 2021.

\bibitem[ZQW20]{zhang2020symmetry}
Yuqian Zhang, Qing Qu, and John Wright.
\newblock From symmetry to geometry: Tractable nonconvex problems.
\newblock {\em arXiv preprint arXiv:2007.06753}, 2020.

\end{thebibliography}

\end{document}